\documentclass[alpha-refs]{wiley-article}

\usepackage{color}
\usepackage{ulem}

\usepackage{siunitx}

\RequirePackage[OT1]{fontenc}
\RequirePackage[colorlinks,citecolor=blue,urlcolor=blue]{hyperref}
\usepackage[linesnumbered, ruled,vlined]{algorithm2e}
\usepackage{url}
  
\usepackage{graphicx}

\usepackage{tikz}
\usetikzlibrary{calc,fadings,decorations.pathreplacing, angles,quotes, shadings}
\usepackage{ragged2e}
\usetikzlibrary{arrows}
\tikzstyle{block}=[draw opacity=1,line width=1.4cm]
\usepackage{pstricks,pst-node,pst-text,pst-3d}
\usepackage{amsmath}
\usepackage{latexsym}
\usepackage{amsmath}
\usepackage{graphicx}
\usepackage{epic}
\usepackage{eepic}
\usepackage{caption}
\usepackage{subcaption}
\usepackage{bbm}
\usepackage{stackengine}
\usepackage{paralist} 

\usepackage{xr}
\definecolor{AliceBlue}{rgb}{0.94,0.97,1.00}
\definecolor{AntiqueWhite1}{rgb}{1.00,0.94,0.86}
\definecolor{AntiqueWhite2}{rgb}{0.93,0.87,0.80}
\definecolor{AntiqueWhite3}{rgb}{0.80,0.75,0.69}
\definecolor{AntiqueWhite4}{rgb}{0.55,0.51,0.47}
\definecolor{AntiqueWhite}{rgb}{0.98,0.92,0.84}
\definecolor{BlanchedAlmond}{rgb}{1.00,0.92,0.80}
\definecolor{BlueViolet}{rgb}{0.54,0.17,0.89}
\definecolor{CadetBlue1}{rgb}{0.60,0.96,1.00}
\definecolor{CadetBlue2}{rgb}{0.56,0.90,0.93}
\definecolor{CadetBlue3}{rgb}{0.48,0.77,0.80}
\definecolor{CadetBlue4}{rgb}{0.33,0.53,0.55}
\definecolor{CadetBlue}{rgb}{0.37,0.62,0.63}
\definecolor{CornflowerBlue}{rgb}{0.39,0.58,0.93}
\definecolor{DarkBlue}{rgb}{0.00,0.00,0.55}
\definecolor{DarkCyan}{rgb}{0.00,0.55,0.55}
\definecolor{DarkGoldenrod1}{rgb}{1.00,0.73,0.06}
\definecolor{DarkGoldenrod2}{rgb}{0.93,0.68,0.05}
\definecolor{DarkGoldenrod3}{rgb}{0.80,0.58,0.05}
\definecolor{DarkGoldenrod4}{rgb}{0.55,0.40,0.03}
\definecolor{DarkGoldenrod}{rgb}{0.72,0.53,0.04}
\definecolor{DarkGray}{rgb}{0.66,0.66,0.66}
\definecolor{DarkGreen}{rgb}{0.00,0.39,0.00}
\definecolor{DarkGrey}{rgb}{0.66,0.66,0.66}
\definecolor{DarkKhaki}{rgb}{0.74,0.72,0.42}
\definecolor{DarkMagenta}{rgb}{0.55,0.00,0.55}
\definecolor{DarkOliveGreen1}{rgb}{0.79,1.00,0.44}
\definecolor{DarkOliveGreen2}{rgb}{0.74,0.93,0.41}
\definecolor{DarkOliveGreen3}{rgb}{0.64,0.80,0.35}
\definecolor{DarkOliveGreen4}{rgb}{0.43,0.55,0.24}
\definecolor{DarkOliveGreen}{rgb}{0.33,0.42,0.18}
\definecolor{DarkOrange1}{rgb}{1.00,0.50,0.00}
\definecolor{DarkOrange2}{rgb}{0.93,0.46,0.00}
\definecolor{DarkOrange3}{rgb}{0.80,0.40,0.00}
\definecolor{DarkOrange4}{rgb}{0.55,0.27,0.00}
\definecolor{DarkOrange}{rgb}{1.00,0.55,0.00}
\definecolor{DarkOrchid1}{rgb}{0.75,0.24,1.00}
\definecolor{DarkOrchid2}{rgb}{0.70,0.23,0.93}
\definecolor{DarkOrchid3}{rgb}{0.60,0.20,0.80}
\definecolor{DarkOrchid4}{rgb}{0.41,0.13,0.55}
\definecolor{DarkOrchid}{rgb}{0.60,0.20,0.80}
\definecolor{DarkRed}{rgb}{0.55,0.00,0.00}
\definecolor{DarkSalmon}{rgb}{0.91,0.59,0.48}
\definecolor{DarkSeaGreen1}{rgb}{0.76,1.00,0.76}
\definecolor{DarkSeaGreen2}{rgb}{0.71,0.93,0.71}
\definecolor{DarkSeaGreen3}{rgb}{0.61,0.80,0.61}
\definecolor{DarkSeaGreen4}{rgb}{0.41,0.55,0.41}
\definecolor{DarkSeaGreen}{rgb}{0.56,0.74,0.56}
\definecolor{DarkSlateBlue}{rgb}{0.28,0.24,0.55}
\definecolor{DarkSlateGray1}{rgb}{0.59,1.00,1.00}
\definecolor{DarkSlateGray2}{rgb}{0.55,0.93,0.93}
\definecolor{DarkSlateGray3}{rgb}{0.47,0.80,0.80}
\definecolor{DarkSlateGray4}{rgb}{0.32,0.55,0.55}
\definecolor{DarkSlateGray}{rgb}{0.18,0.31,0.31}
\definecolor{DarkSlateGrey}{rgb}{0.18,0.31,0.31}
\definecolor{DarkTurquoise}{rgb}{0.00,0.81,0.82}
\definecolor{DarkViolet}{rgb}{0.58,0.00,0.83}
\definecolor{DeepPink1}{rgb}{1.00,0.08,0.58}
\definecolor{DeepPink2}{rgb}{0.93,0.07,0.54}
\definecolor{DeepPink3}{rgb}{0.80,0.06,0.46}
\definecolor{DeepPink4}{rgb}{0.55,0.04,0.31}
\definecolor{DeepPink}{rgb}{1.00,0.08,0.58}
\definecolor{DeepSkyBlue1}{rgb}{0.00,0.75,1.00}
\definecolor{DeepSkyBlue2}{rgb}{0.00,0.70,0.93}
\definecolor{DeepSkyBlue3}{rgb}{0.00,0.60,0.80}
\definecolor{DeepSkyBlue4}{rgb}{0.00,0.41,0.55}
\definecolor{DeepSkyBlue}{rgb}{0.00,0.75,1.00}
\definecolor{DimGray}{rgb}{0.41,0.41,0.41}
\definecolor{DimGrey}{rgb}{0.41,0.41,0.41}
\definecolor{DodgerBlue1}{rgb}{0.12,0.56,1.00}
\definecolor{DodgerBlue2}{rgb}{0.11,0.53,0.93}
\definecolor{DodgerBlue3}{rgb}{0.09,0.45,0.80}
\definecolor{DodgerBlue4}{rgb}{0.06,0.31,0.55}
\definecolor{DodgerBlue}{rgb}{0.12,0.56,1.00}
\definecolor{FloralWhite}{rgb}{1.00,0.98,0.94}
\definecolor{ForestGreen}{rgb}{0.13,0.55,0.13}
\definecolor{GhostWhite}{rgb}{0.97,0.97,1.00}
\definecolor{GreenYellow}{rgb}{0.68,1.00,0.18}
\definecolor{HotPink1}{rgb}{1.00,0.43,0.71}
\definecolor{HotPink2}{rgb}{0.93,0.42,0.65}
\definecolor{HotPink3}{rgb}{0.80,0.38,0.56}
\definecolor{HotPink4}{rgb}{0.55,0.23,0.38}
\definecolor{HotPink}{rgb}{1.00,0.41,0.71}
\definecolor{IndianRed1}{rgb}{1.00,0.42,0.42}
\definecolor{IndianRed2}{rgb}{0.93,0.39,0.39}
\definecolor{IndianRed3}{rgb}{0.80,0.33,0.33}
\definecolor{IndianRed4}{rgb}{0.55,0.23,0.23}
\definecolor{IndianRed}{rgb}{0.80,0.36,0.36}
\definecolor{LavenderBlush1}{rgb}{1.00,0.94,0.96}
\definecolor{LavenderBlush2}{rgb}{0.93,0.88,0.90}
\definecolor{LavenderBlush3}{rgb}{0.80,0.76,0.77}
\definecolor{LavenderBlush4}{rgb}{0.55,0.51,0.53}
\definecolor{LavenderBlush}{rgb}{1.00,0.94,0.96}
\definecolor{LawnGreen}{rgb}{0.49,0.99,0.00}
\definecolor{LemonChiffon1}{rgb}{1.00,0.98,0.80}
\definecolor{LemonChiffon2}{rgb}{0.93,0.91,0.75}
\definecolor{LemonChiffon3}{rgb}{0.80,0.79,0.65}
\definecolor{LemonChiffon4}{rgb}{0.55,0.54,0.44}
\definecolor{LemonChiffon}{rgb}{1.00,0.98,0.80}
\definecolor{LightBlue1}{rgb}{0.75,0.94,1.00}
\definecolor{LightBlue2}{rgb}{0.70,0.87,0.93}
\definecolor{LightBlue3}{rgb}{0.60,0.75,0.80}
\definecolor{LightBlue4}{rgb}{0.41,0.51,0.55}
\definecolor{LightBlue}{rgb}{0.68,0.85,0.90}
\definecolor{LightCoral}{rgb}{0.94,0.50,0.50}
\definecolor{LightCyan1}{rgb}{0.88,1.00,1.00}
\definecolor{LightCyan2}{rgb}{0.82,0.93,0.93}
\definecolor{LightCyan3}{rgb}{0.71,0.80,0.80}
\definecolor{LightCyan4}{rgb}{0.48,0.55,0.55}
\definecolor{LightCyan}{rgb}{0.88,1.00,1.00}
\definecolor{LightGoldenrod1}{rgb}{1.00,0.93,0.55}
\definecolor{LightGoldenrod2}{rgb}{0.93,0.86,0.51}
\definecolor{LightGoldenrod3}{rgb}{0.80,0.75,0.44}
\definecolor{LightGoldenrod4}{rgb}{0.55,0.51,0.30}
\definecolor{LightGoldenrodYellow}{rgb}{0.98,0.98,0.82}
\definecolor{LightGoldenrod}{rgb}{0.93,0.87,0.51}
\definecolor{LightGray}{rgb}{0.83,0.83,0.83}
\definecolor{LightGreen}{rgb}{0.56,0.93,0.56}
\definecolor{LightGrey}{rgb}{0.83,0.83,0.83}
\definecolor{LightPink1}{rgb}{1.00,0.68,0.73}
\definecolor{LightPink2}{rgb}{0.93,0.64,0.68}
\definecolor{LightPink3}{rgb}{0.80,0.55,0.58}
\definecolor{LightPink4}{rgb}{0.55,0.37,0.40}
\definecolor{LightPink}{rgb}{1.00,0.71,0.76}
\definecolor{LightSalmon1}{rgb}{1.00,0.63,0.48}
\definecolor{LightSalmon2}{rgb}{0.93,0.58,0.45}
\definecolor{LightSalmon3}{rgb}{0.80,0.51,0.38}
\definecolor{LightSalmon4}{rgb}{0.55,0.34,0.26}
\definecolor{LightSalmon}{rgb}{1.00,0.63,0.48}
\definecolor{LightSeaGreen}{rgb}{0.13,0.70,0.67}
\definecolor{LightSkyBlue1}{rgb}{0.69,0.89,1.00}
\definecolor{LightSkyBlue2}{rgb}{0.64,0.83,0.93}
\definecolor{LightSkyBlue3}{rgb}{0.55,0.71,0.80}
\definecolor{LightSkyBlue4}{rgb}{0.38,0.48,0.55}
\definecolor{LightSkyBlue}{rgb}{0.53,0.81,0.98}
\definecolor{LightSlateBlue}{rgb}{0.52,0.44,1.00}
\definecolor{LightSlateGray}{rgb}{0.47,0.53,0.60}
\definecolor{LightSlateGrey}{rgb}{0.47,0.53,0.60}
\definecolor{LightSteelBlue1}{rgb}{0.79,0.88,1.00}
\definecolor{LightSteelBlue2}{rgb}{0.74,0.82,0.93}
\definecolor{LightSteelBlue3}{rgb}{0.64,0.71,0.80}
\definecolor{LightSteelBlue4}{rgb}{0.43,0.48,0.55}
\definecolor{LightSteelBlue}{rgb}{0.69,0.77,0.87}
\definecolor{LightYellow1}{rgb}{1.00,1.00,0.88}
\definecolor{LightYellow2}{rgb}{0.93,0.93,0.82}
\definecolor{LightYellow3}{rgb}{0.80,0.80,0.71}
\definecolor{LightYellow4}{rgb}{0.55,0.55,0.48}
\definecolor{LightYellow}{rgb}{1.00,1.00,0.88}
\definecolor{LimeGreen}{rgb}{0.20,0.80,0.20}
\definecolor{MediumAquamarine}{rgb}{0.40,0.80,0.67}
\definecolor{MediumBlue}{rgb}{0.00,0.00,0.80}
\definecolor{MediumOrchid1}{rgb}{0.88,0.40,1.00}
\definecolor{MediumOrchid2}{rgb}{0.82,0.37,0.93}
\definecolor{MediumOrchid3}{rgb}{0.71,0.32,0.80}
\definecolor{MediumOrchid4}{rgb}{0.48,0.22,0.55}
\definecolor{MediumOrchid}{rgb}{0.73,0.33,0.83}
\definecolor{MediumPurple1}{rgb}{0.67,0.51,1.00}
\definecolor{MediumPurple2}{rgb}{0.62,0.47,0.93}
\definecolor{MediumPurple3}{rgb}{0.54,0.41,0.80}
\definecolor{MediumPurple4}{rgb}{0.36,0.28,0.55}
\definecolor{MediumPurple}{rgb}{0.58,0.44,0.86}
\definecolor{MediumSeaGreen}{rgb}{0.24,0.70,0.44}
\definecolor{MediumSlateBlue}{rgb}{0.48,0.41,0.93}
\definecolor{MediumSpringGreen}{rgb}{0.00,0.98,0.60}
\definecolor{MediumTurquoise}{rgb}{0.28,0.82,0.80}
\definecolor{MediumVioletRed}{rgb}{0.78,0.08,0.52}
\definecolor{MidnightBlue}{rgb}{0.10,0.10,0.44}
\definecolor{MintCream}{rgb}{0.96,1.00,0.98}
\definecolor{MistyRose1}{rgb}{1.00,0.89,0.88}
\definecolor{MistyRose2}{rgb}{0.93,0.84,0.82}
\definecolor{MistyRose3}{rgb}{0.80,0.72,0.71}
\definecolor{MistyRose4}{rgb}{0.55,0.49,0.48}
\definecolor{MistyRose}{rgb}{1.00,0.89,0.88}
\definecolor{NavajoWhite1}{rgb}{1.00,0.87,0.68}
\definecolor{NavajoWhite2}{rgb}{0.93,0.81,0.63}
\definecolor{NavajoWhite3}{rgb}{0.80,0.70,0.55}
\definecolor{NavajoWhite4}{rgb}{0.55,0.47,0.37}
\definecolor{NavajoWhite}{rgb}{1.00,0.87,0.68}
\definecolor{NavyBlue}{rgb}{0.00,0.00,0.50}
\definecolor{OldLace}{rgb}{0.99,0.96,0.90}
\definecolor{OliveDrab1}{rgb}{0.75,1.00,0.24}
\definecolor{OliveDrab2}{rgb}{0.70,0.93,0.23}
\definecolor{OliveDrab3}{rgb}{0.60,0.80,0.20}
\definecolor{OliveDrab4}{rgb}{0.41,0.55,0.13}
\definecolor{OliveDrab}{rgb}{0.42,0.56,0.14}
\definecolor{OrangeRed1}{rgb}{1.00,0.27,0.00}
\definecolor{OrangeRed2}{rgb}{0.93,0.25,0.00}
\definecolor{OrangeRed3}{rgb}{0.80,0.22,0.00}
\definecolor{OrangeRed4}{rgb}{0.55,0.15,0.00}
\definecolor{OrangeRed}{rgb}{1.00,0.27,0.00}
\definecolor{PaleGoldenrod}{rgb}{0.93,0.91,0.67}
\definecolor{PaleGreen1}{rgb}{0.60,1.00,0.60}
\definecolor{PaleGreen2}{rgb}{0.56,0.93,0.56}
\definecolor{PaleGreen3}{rgb}{0.49,0.80,0.49}
\definecolor{PaleGreen4}{rgb}{0.33,0.55,0.33}
\definecolor{PaleGreen}{rgb}{0.60,0.98,0.60}
\definecolor{PaleTurquoise1}{rgb}{0.73,1.00,1.00}
\definecolor{PaleTurquoise2}{rgb}{0.68,0.93,0.93}
\definecolor{PaleTurquoise3}{rgb}{0.59,0.80,0.80}
\definecolor{PaleTurquoise4}{rgb}{0.40,0.55,0.55}
\definecolor{PaleTurquoise}{rgb}{0.69,0.93,0.93}
\definecolor{PaleVioletRed1}{rgb}{1.00,0.51,0.67}
\definecolor{PaleVioletRed2}{rgb}{0.93,0.47,0.62}
\definecolor{PaleVioletRed3}{rgb}{0.80,0.41,0.54}
\definecolor{PaleVioletRed4}{rgb}{0.55,0.28,0.36}
\definecolor{PaleVioletRed}{rgb}{0.86,0.44,0.58}
\definecolor{PapayaWhip}{rgb}{1.00,0.94,0.84}
\definecolor{PeachPuff1}{rgb}{1.00,0.85,0.73}
\definecolor{PeachPuff2}{rgb}{0.93,0.80,0.68}
\definecolor{PeachPuff3}{rgb}{0.80,0.69,0.58}
\definecolor{PeachPuff4}{rgb}{0.55,0.47,0.40}
\definecolor{PeachPuff}{rgb}{1.00,0.85,0.73}
\definecolor{PowderBlue}{rgb}{0.69,0.88,0.90}
\definecolor{RosyBrown1}{rgb}{1.00,0.76,0.76}
\definecolor{RosyBrown2}{rgb}{0.93,0.71,0.71}
\definecolor{RosyBrown3}{rgb}{0.80,0.61,0.61}
\definecolor{RosyBrown4}{rgb}{0.55,0.41,0.41}
\definecolor{RosyBrown}{rgb}{0.74,0.56,0.56}
\definecolor{RoyalBlue1}{rgb}{0.28,0.46,1.00}
\definecolor{RoyalBlue2}{rgb}{0.26,0.43,0.93}
\definecolor{RoyalBlue3}{rgb}{0.23,0.37,0.80}
\definecolor{RoyalBlue4}{rgb}{0.15,0.25,0.55}
\definecolor{RoyalBlue}{rgb}{0.25,0.41,0.88}
\definecolor{SaddleBrown}{rgb}{0.55,0.27,0.07}
\definecolor{SandyBrown}{rgb}{0.96,0.64,0.38}
\definecolor{SeaGreen1}{rgb}{0.33,1.00,0.62}
\definecolor{SeaGreen2}{rgb}{0.31,0.93,0.58}
\definecolor{SeaGreen3}{rgb}{0.26,0.80,0.50}
\definecolor{SeaGreen4}{rgb}{0.18,0.55,0.34}
\definecolor{SeaGreen}{rgb}{0.18,0.55,0.34}
\definecolor{SkyBlue1}{rgb}{0.53,0.81,1.00}
\definecolor{SkyBlue2}{rgb}{0.49,0.75,0.93}
\definecolor{SkyBlue3}{rgb}{0.42,0.65,0.80}
\definecolor{SkyBlue4}{rgb}{0.29,0.44,0.55}
\definecolor{SkyBlue}{rgb}{0.53,0.81,0.92}
\definecolor{SlateBlue1}{rgb}{0.51,0.44,1.00}
\definecolor{SlateBlue2}{rgb}{0.48,0.40,0.93}
\definecolor{SlateBlue3}{rgb}{0.41,0.35,0.80}
\definecolor{SlateBlue4}{rgb}{0.28,0.24,0.55}
\definecolor{SlateBlue}{rgb}{0.42,0.35,0.80}
\definecolor{SlateGray1}{rgb}{0.78,0.89,1.00}
\definecolor{SlateGray2}{rgb}{0.73,0.83,0.93}
\definecolor{SlateGray3}{rgb}{0.62,0.71,0.80}
\definecolor{SlateGray4}{rgb}{0.42,0.48,0.55}
\definecolor{SlateGray}{rgb}{0.44,0.50,0.56}
\definecolor{SlateGrey}{rgb}{0.44,0.50,0.56}
\definecolor{SpringGreen1}{rgb}{0.00,1.00,0.50}
\definecolor{SpringGreen2}{rgb}{0.00,0.93,0.46}
\definecolor{SpringGreen3}{rgb}{0.00,0.80,0.40}
\definecolor{SpringGreen4}{rgb}{0.00,0.55,0.27}
\definecolor{SpringGreen}{rgb}{0.00,1.00,0.50}
\definecolor{SteelBlue1}{rgb}{0.39,0.72,1.00}
\definecolor{SteelBlue2}{rgb}{0.36,0.67,0.93}
\definecolor{SteelBlue3}{rgb}{0.31,0.58,0.80}
\definecolor{SteelBlue4}{rgb}{0.21,0.39,0.55}
\definecolor{SteelBlue}{rgb}{0.27,0.51,0.71}
\definecolor{VioletRed1}{rgb}{1.00,0.24,0.59}
\definecolor{VioletRed2}{rgb}{0.93,0.23,0.55}
\definecolor{VioletRed3}{rgb}{0.80,0.20,0.47}
\definecolor{VioletRed4}{rgb}{0.55,0.13,0.32}
\definecolor{VioletRed}{rgb}{0.82,0.13,0.56}
\definecolor{WhiteSmoke}{rgb}{0.96,0.96,0.96}
\definecolor{YellowGreen}{rgb}{0.60,0.80,0.20}
\definecolor{aliceblue}{rgb}{0.94,0.97,1.00}
\definecolor{antiquewhite}{rgb}{0.98,0.92,0.84}
\definecolor{aquamarine1}{rgb}{0.50,1.00,0.83}
\definecolor{aquamarine2}{rgb}{0.46,0.93,0.78}
\definecolor{aquamarine3}{rgb}{0.40,0.80,0.67}
\definecolor{aquamarine4}{rgb}{0.27,0.55,0.45}
\definecolor{aquamarine}{rgb}{0.50,1.00,0.83}
\definecolor{azure1}{rgb}{0.94,1.00,1.00}
\definecolor{azure2}{rgb}{0.88,0.93,0.93}
\definecolor{azure3}{rgb}{0.76,0.80,0.80}
\definecolor{azure4}{rgb}{0.51,0.55,0.55}
\definecolor{azure}{rgb}{0.94,1.00,1.00}
\definecolor{beige}{rgb}{0.96,0.96,0.86}
\definecolor{bisque1}{rgb}{1.00,0.89,0.77}
\definecolor{bisque2}{rgb}{0.93,0.84,0.72}
\definecolor{bisque3}{rgb}{0.80,0.72,0.62}
\definecolor{bisque4}{rgb}{0.55,0.49,0.42}
\definecolor{bisque}{rgb}{1.00,0.89,0.77}
\definecolor{black}{rgb}{0.00,0.00,0.00}
\definecolor{blanchedalmond}{rgb}{1.00,0.92,0.80}
\definecolor{blue1}{rgb}{0.00,0.00,1.00}
\definecolor{blue2}{rgb}{0.00,0.00,0.93}
\definecolor{blue3}{rgb}{0.00,0.00,0.80}
\definecolor{blue4}{rgb}{0.00,0.00,0.55}
\definecolor{blueviolet}{rgb}{0.54,0.17,0.89}
\definecolor{blue}{rgb}{0.00,0.00,1.00}
\definecolor{brown1}{rgb}{1.00,0.25,0.25}
\definecolor{brown2}{rgb}{0.93,0.23,0.23}
\definecolor{brown3}{rgb}{0.80,0.20,0.20}
\definecolor{brown4}{rgb}{0.55,0.14,0.14}
\definecolor{brown}{rgb}{0.65,0.16,0.16}
\definecolor{burlywood1}{rgb}{1.00,0.83,0.61}
\definecolor{burlywood2}{rgb}{0.93,0.77,0.57}
\definecolor{burlywood3}{rgb}{0.80,0.67,0.49}
\definecolor{burlywood4}{rgb}{0.55,0.45,0.33}
\definecolor{burlywood}{rgb}{0.87,0.72,0.53}
\definecolor{cadetblue}{rgb}{0.37,0.62,0.63}
\definecolor{chartreuse1}{rgb}{0.50,1.00,0.00}
\definecolor{chartreuse2}{rgb}{0.46,0.93,0.00}
\definecolor{chartreuse3}{rgb}{0.40,0.80,0.00}
\definecolor{chartreuse4}{rgb}{0.27,0.55,0.00}
\definecolor{chartreuse}{rgb}{0.50,1.00,0.00}
\definecolor{chocolate1}{rgb}{1.00,0.50,0.14}
\definecolor{chocolate2}{rgb}{0.93,0.46,0.13}
\definecolor{chocolate3}{rgb}{0.80,0.40,0.11}
\definecolor{chocolate4}{rgb}{0.55,0.27,0.07}
\definecolor{chocolate}{rgb}{0.82,0.41,0.12}
\definecolor{coral1}{rgb}{1.00,0.45,0.34}
\definecolor{coral2}{rgb}{0.93,0.42,0.31}
\definecolor{coral3}{rgb}{0.80,0.36,0.27}
\definecolor{coral4}{rgb}{0.55,0.24,0.18}
\definecolor{coral}{rgb}{1.00,0.50,0.31}
\definecolor{cornflowerblue}{rgb}{0.39,0.58,0.93}
\definecolor{cornsilk1}{rgb}{1.00,0.97,0.86}
\definecolor{cornsilk2}{rgb}{0.93,0.91,0.80}
\definecolor{cornsilk3}{rgb}{0.80,0.78,0.69}
\definecolor{cornsilk4}{rgb}{0.55,0.53,0.47}
\definecolor{cornsilk}{rgb}{1.00,0.97,0.86}
\definecolor{cyan1}{rgb}{0.00,1.00,1.00}
\definecolor{cyan2}{rgb}{0.00,0.93,0.93}
\definecolor{cyan3}{rgb}{0.00,0.80,0.80}
\definecolor{cyan4}{rgb}{0.00,0.55,0.55}
\definecolor{cyan}{rgb}{0.00,1.00,1.00}
\definecolor{darkblue}{rgb}{0.00,0.00,0.55}
\definecolor{darkcyan}{rgb}{0.00,0.55,0.55}
\definecolor{darkgoldenrod}{rgb}{0.72,0.53,0.04}
\definecolor{darkgray}{rgb}{0.66,0.66,0.66}
\definecolor{darkgreen}{rgb}{0.00,0.39,0.00}
\definecolor{darkgrey}{rgb}{0.66,0.66,0.66}
\definecolor{darkkhaki}{rgb}{0.74,0.72,0.42}
\definecolor{darkmagenta}{rgb}{0.55,0.00,0.55}
\definecolor{darkolive}{rgb}{0.33,0.42,0.18}
\definecolor{darkorange}{rgb}{1.00,0.55,0.00}
\definecolor{darkorchid}{rgb}{0.60,0.20,0.80}
\definecolor{darkred}{rgb}{0.55,0.00,0.00}
\definecolor{darksalmon}{rgb}{0.91,0.59,0.48}
\definecolor{darksea}{rgb}{0.56,0.74,0.56}
\definecolor{darkslate}{rgb}{0.18,0.31,0.31}
\definecolor{darkslate}{rgb}{0.18,0.31,0.31}
\definecolor{darkslate}{rgb}{0.28,0.24,0.55}
\definecolor{darkturquoise}{rgb}{0.00,0.81,0.82}
\definecolor{darkviolet}{rgb}{0.58,0.00,0.83}
\definecolor{deeppink}{rgb}{1.00,0.08,0.58}
\definecolor{deepsky}{rgb}{0.00,0.75,1.00}
\definecolor{dimgray}{rgb}{0.41,0.41,0.41}
\definecolor{dimgrey}{rgb}{0.41,0.41,0.41}
\definecolor{dodgerblue}{rgb}{0.12,0.56,1.00}
\definecolor{firebrick1}{rgb}{1.00,0.19,0.19}
\definecolor{firebrick2}{rgb}{0.93,0.17,0.17}
\definecolor{firebrick3}{rgb}{0.80,0.15,0.15}
\definecolor{firebrick4}{rgb}{0.55,0.10,0.10}
\definecolor{firebrick}{rgb}{0.70,0.13,0.13}
\definecolor{floralwhite}{rgb}{1.00,0.98,0.94}
\definecolor{forestgreen}{rgb}{0.13,0.55,0.13}
\definecolor{gainsboro}{rgb}{0.86,0.86,0.86}
\definecolor{ghostwhite}{rgb}{0.97,0.97,1.00}
\definecolor{gold1}{rgb}{1.00,0.84,0.00}
\definecolor{gold2}{rgb}{0.93,0.79,0.00}
\definecolor{gold3}{rgb}{0.80,0.68,0.00}
\definecolor{gold4}{rgb}{0.55,0.46,0.00}
\definecolor{goldenrod1}{rgb}{1.00,0.76,0.15}
\definecolor{goldenrod2}{rgb}{0.93,0.71,0.13}
\definecolor{goldenrod3}{rgb}{0.80,0.61,0.11}
\definecolor{goldenrod4}{rgb}{0.55,0.41,0.08}
\definecolor{goldenrod}{rgb}{0.85,0.65,0.13}
\definecolor{gold}{rgb}{1.00,0.84,0.00}
\definecolor{gray0}{rgb}{0.00,0.00,0.00}
\definecolor{gray100}{rgb}{1.00,1.00,1.00}
\definecolor{gray10}{rgb}{0.10,0.10,0.10}
\definecolor{gray11}{rgb}{0.11,0.11,0.11}
\definecolor{gray12}{rgb}{0.12,0.12,0.12}
\definecolor{gray13}{rgb}{0.13,0.13,0.13}
\definecolor{gray14}{rgb}{0.14,0.14,0.14}
\definecolor{gray15}{rgb}{0.15,0.15,0.15}
\definecolor{gray16}{rgb}{0.16,0.16,0.16}
\definecolor{gray17}{rgb}{0.17,0.17,0.17}
\definecolor{gray18}{rgb}{0.18,0.18,0.18}
\definecolor{gray19}{rgb}{0.19,0.19,0.19}
\definecolor{gray1}{rgb}{0.01,0.01,0.01}
\definecolor{gray20}{rgb}{0.20,0.20,0.20}
\definecolor{gray21}{rgb}{0.21,0.21,0.21}
\definecolor{gray22}{rgb}{0.22,0.22,0.22}
\definecolor{gray23}{rgb}{0.23,0.23,0.23}
\definecolor{gray24}{rgb}{0.24,0.24,0.24}
\definecolor{gray25}{rgb}{0.25,0.25,0.25}
\definecolor{gray26}{rgb}{0.26,0.26,0.26}
\definecolor{gray27}{rgb}{0.27,0.27,0.27}
\definecolor{gray28}{rgb}{0.28,0.28,0.28}
\definecolor{gray29}{rgb}{0.29,0.29,0.29}
\definecolor{gray2}{rgb}{0.02,0.02,0.02}
\definecolor{gray30}{rgb}{0.30,0.30,0.30}
\definecolor{gray31}{rgb}{0.31,0.31,0.31}
\definecolor{gray32}{rgb}{0.32,0.32,0.32}
\definecolor{gray33}{rgb}{0.33,0.33,0.33}
\definecolor{gray34}{rgb}{0.34,0.34,0.34}
\definecolor{gray35}{rgb}{0.35,0.35,0.35}
\definecolor{gray36}{rgb}{0.36,0.36,0.36}
\definecolor{gray37}{rgb}{0.37,0.37,0.37}
\definecolor{gray38}{rgb}{0.38,0.38,0.38}
\definecolor{gray39}{rgb}{0.39,0.39,0.39}
\definecolor{gray3}{rgb}{0.03,0.03,0.03}
\definecolor{gray40}{rgb}{0.40,0.40,0.40}
\definecolor{gray41}{rgb}{0.41,0.41,0.41}
\definecolor{gray42}{rgb}{0.42,0.42,0.42}
\definecolor{gray43}{rgb}{0.43,0.43,0.43}
\definecolor{gray44}{rgb}{0.44,0.44,0.44}
\definecolor{gray45}{rgb}{0.45,0.45,0.45}
\definecolor{gray46}{rgb}{0.46,0.46,0.46}
\definecolor{gray47}{rgb}{0.47,0.47,0.47}
\definecolor{gray48}{rgb}{0.48,0.48,0.48}
\definecolor{gray49}{rgb}{0.49,0.49,0.49}
\definecolor{gray4}{rgb}{0.04,0.04,0.04}
\definecolor{gray50}{rgb}{0.50,0.50,0.50}
\definecolor{gray51}{rgb}{0.51,0.51,0.51}
\definecolor{gray52}{rgb}{0.52,0.52,0.52}
\definecolor{gray53}{rgb}{0.53,0.53,0.53}
\definecolor{gray54}{rgb}{0.54,0.54,0.54}
\definecolor{gray55}{rgb}{0.55,0.55,0.55}
\definecolor{gray56}{rgb}{0.56,0.56,0.56}
\definecolor{gray57}{rgb}{0.57,0.57,0.57}
\definecolor{gray58}{rgb}{0.58,0.58,0.58}
\definecolor{gray59}{rgb}{0.59,0.59,0.59}
\definecolor{gray5}{rgb}{0.05,0.05,0.05}
\definecolor{gray60}{rgb}{0.60,0.60,0.60}
\definecolor{gray61}{rgb}{0.61,0.61,0.61}
\definecolor{gray62}{rgb}{0.62,0.62,0.62}
\definecolor{gray63}{rgb}{0.63,0.63,0.63}
\definecolor{gray64}{rgb}{0.64,0.64,0.64}
\definecolor{gray65}{rgb}{0.65,0.65,0.65}
\definecolor{gray66}{rgb}{0.66,0.66,0.66}
\definecolor{gray67}{rgb}{0.67,0.67,0.67}
\definecolor{gray68}{rgb}{0.68,0.68,0.68}
\definecolor{gray69}{rgb}{0.69,0.69,0.69}
\definecolor{gray6}{rgb}{0.06,0.06,0.06}
\definecolor{gray70}{rgb}{0.70,0.70,0.70}
\definecolor{gray71}{rgb}{0.71,0.71,0.71}
\definecolor{gray72}{rgb}{0.72,0.72,0.72}
\definecolor{gray73}{rgb}{0.73,0.73,0.73}
\definecolor{gray74}{rgb}{0.74,0.74,0.74}
\definecolor{gray75}{rgb}{0.75,0.75,0.75}
\definecolor{gray76}{rgb}{0.76,0.76,0.76}
\definecolor{gray77}{rgb}{0.77,0.77,0.77}
\definecolor{gray78}{rgb}{0.78,0.78,0.78}
\definecolor{gray79}{rgb}{0.79,0.79,0.79}
\definecolor{gray7}{rgb}{0.07,0.07,0.07}
\definecolor{gray80}{rgb}{0.80,0.80,0.80}
\definecolor{gray81}{rgb}{0.81,0.81,0.81}
\definecolor{gray82}{rgb}{0.82,0.82,0.82}
\definecolor{gray83}{rgb}{0.83,0.83,0.83}
\definecolor{gray84}{rgb}{0.84,0.84,0.84}
\definecolor{gray85}{rgb}{0.85,0.85,0.85}
\definecolor{gray86}{rgb}{0.86,0.86,0.86}
\definecolor{gray87}{rgb}{0.87,0.87,0.87}
\definecolor{gray88}{rgb}{0.88,0.88,0.88}
\definecolor{gray89}{rgb}{0.89,0.89,0.89}
\definecolor{gray8}{rgb}{0.08,0.08,0.08}
\definecolor{gray90}{rgb}{0.90,0.90,0.90}
\definecolor{gray91}{rgb}{0.91,0.91,0.91}
\definecolor{gray92}{rgb}{0.92,0.92,0.92}
\definecolor{gray93}{rgb}{0.93,0.93,0.93}
\definecolor{gray94}{rgb}{0.94,0.94,0.94}
\definecolor{gray95}{rgb}{0.95,0.95,0.95}
\definecolor{gray96}{rgb}{0.96,0.96,0.96}
\definecolor{gray97}{rgb}{0.97,0.97,0.97}
\definecolor{gray98}{rgb}{0.98,0.98,0.98}
\definecolor{gray99}{rgb}{0.99,0.99,0.99}
\definecolor{gray9}{rgb}{0.09,0.09,0.09}
\definecolor{gray}{rgb}{0.75,0.75,0.75}
\definecolor{green1}{rgb}{0.00,1.00,0.00}
\definecolor{green2}{rgb}{0.00,0.93,0.00}
\definecolor{green3}{rgb}{0.00,0.80,0.00}
\definecolor{green4}{rgb}{0.00,0.55,0.00}
\definecolor{greenyellow}{rgb}{0.68,1.00,0.18}
\definecolor{green}{rgb}{0.00,1.00,0.00}
\definecolor{grey0}{rgb}{0.00,0.00,0.00}
\definecolor{grey100}{rgb}{1.00,1.00,1.00}
\definecolor{grey10}{rgb}{0.10,0.10,0.10}
\definecolor{grey11}{rgb}{0.11,0.11,0.11}
\definecolor{grey12}{rgb}{0.12,0.12,0.12}
\definecolor{grey13}{rgb}{0.13,0.13,0.13}
\definecolor{grey14}{rgb}{0.14,0.14,0.14}
\definecolor{grey15}{rgb}{0.15,0.15,0.15}
\definecolor{grey16}{rgb}{0.16,0.16,0.16}
\definecolor{grey17}{rgb}{0.17,0.17,0.17}
\definecolor{grey18}{rgb}{0.18,0.18,0.18}
\definecolor{grey19}{rgb}{0.19,0.19,0.19}
\definecolor{grey1}{rgb}{0.01,0.01,0.01}
\definecolor{grey20}{rgb}{0.20,0.20,0.20}
\definecolor{grey21}{rgb}{0.21,0.21,0.21}
\definecolor{grey22}{rgb}{0.22,0.22,0.22}
\definecolor{grey23}{rgb}{0.23,0.23,0.23}
\definecolor{grey24}{rgb}{0.24,0.24,0.24}
\definecolor{grey25}{rgb}{0.25,0.25,0.25}
\definecolor{grey26}{rgb}{0.26,0.26,0.26}
\definecolor{grey27}{rgb}{0.27,0.27,0.27}
\definecolor{grey28}{rgb}{0.28,0.28,0.28}
\definecolor{grey29}{rgb}{0.29,0.29,0.29}
\definecolor{grey2}{rgb}{0.02,0.02,0.02}
\definecolor{grey30}{rgb}{0.30,0.30,0.30}
\definecolor{grey31}{rgb}{0.31,0.31,0.31}
\definecolor{grey32}{rgb}{0.32,0.32,0.32}
\definecolor{grey33}{rgb}{0.33,0.33,0.33}
\definecolor{grey34}{rgb}{0.34,0.34,0.34}
\definecolor{grey35}{rgb}{0.35,0.35,0.35}
\definecolor{grey36}{rgb}{0.36,0.36,0.36}
\definecolor{grey37}{rgb}{0.37,0.37,0.37}
\definecolor{grey38}{rgb}{0.38,0.38,0.38}
\definecolor{grey39}{rgb}{0.39,0.39,0.39}
\definecolor{grey3}{rgb}{0.03,0.03,0.03}
\definecolor{grey40}{rgb}{0.40,0.40,0.40}
\definecolor{grey41}{rgb}{0.41,0.41,0.41}
\definecolor{grey42}{rgb}{0.42,0.42,0.42}
\definecolor{grey43}{rgb}{0.43,0.43,0.43}
\definecolor{grey44}{rgb}{0.44,0.44,0.44}
\definecolor{grey45}{rgb}{0.45,0.45,0.45}
\definecolor{grey46}{rgb}{0.46,0.46,0.46}
\definecolor{grey47}{rgb}{0.47,0.47,0.47}
\definecolor{grey48}{rgb}{0.48,0.48,0.48}
\definecolor{grey49}{rgb}{0.49,0.49,0.49}
\definecolor{grey4}{rgb}{0.04,0.04,0.04}
\definecolor{grey50}{rgb}{0.50,0.50,0.50}
\definecolor{grey51}{rgb}{0.51,0.51,0.51}
\definecolor{grey52}{rgb}{0.52,0.52,0.52}
\definecolor{grey53}{rgb}{0.53,0.53,0.53}
\definecolor{grey54}{rgb}{0.54,0.54,0.54}
\definecolor{grey55}{rgb}{0.55,0.55,0.55}
\definecolor{grey56}{rgb}{0.56,0.56,0.56}
\definecolor{grey57}{rgb}{0.57,0.57,0.57}
\definecolor{grey58}{rgb}{0.58,0.58,0.58}
\definecolor{grey59}{rgb}{0.59,0.59,0.59}
\definecolor{grey5}{rgb}{0.05,0.05,0.05}
\definecolor{grey60}{rgb}{0.60,0.60,0.60}
\definecolor{grey61}{rgb}{0.61,0.61,0.61}
\definecolor{grey62}{rgb}{0.62,0.62,0.62}
\definecolor{grey63}{rgb}{0.63,0.63,0.63}
\definecolor{grey64}{rgb}{0.64,0.64,0.64}
\definecolor{grey65}{rgb}{0.65,0.65,0.65}
\definecolor{grey66}{rgb}{0.66,0.66,0.66}
\definecolor{grey67}{rgb}{0.67,0.67,0.67}
\definecolor{grey68}{rgb}{0.68,0.68,0.68}
\definecolor{grey69}{rgb}{0.69,0.69,0.69}
\definecolor{grey6}{rgb}{0.06,0.06,0.06}
\definecolor{grey70}{rgb}{0.70,0.70,0.70}
\definecolor{grey71}{rgb}{0.71,0.71,0.71}
\definecolor{grey72}{rgb}{0.72,0.72,0.72}
\definecolor{grey73}{rgb}{0.73,0.73,0.73}
\definecolor{grey74}{rgb}{0.74,0.74,0.74}
\definecolor{grey75}{rgb}{0.75,0.75,0.75}
\definecolor{grey76}{rgb}{0.76,0.76,0.76}
\definecolor{grey77}{rgb}{0.77,0.77,0.77}
\definecolor{grey78}{rgb}{0.78,0.78,0.78}
\definecolor{grey79}{rgb}{0.79,0.79,0.79}
\definecolor{grey7}{rgb}{0.07,0.07,0.07}
\definecolor{grey80}{rgb}{0.80,0.80,0.80}
\definecolor{grey81}{rgb}{0.81,0.81,0.81}
\definecolor{grey82}{rgb}{0.82,0.82,0.82}
\definecolor{grey83}{rgb}{0.83,0.83,0.83}
\definecolor{grey84}{rgb}{0.84,0.84,0.84}
\definecolor{grey85}{rgb}{0.85,0.85,0.85}
\definecolor{grey86}{rgb}{0.86,0.86,0.86}
\definecolor{grey87}{rgb}{0.87,0.87,0.87}
\definecolor{grey88}{rgb}{0.88,0.88,0.88}
\definecolor{grey89}{rgb}{0.89,0.89,0.89}
\definecolor{grey8}{rgb}{0.08,0.08,0.08}
\definecolor{grey90}{rgb}{0.90,0.90,0.90}
\definecolor{grey91}{rgb}{0.91,0.91,0.91}
\definecolor{grey92}{rgb}{0.92,0.92,0.92}
\definecolor{grey93}{rgb}{0.93,0.93,0.93}
\definecolor{grey94}{rgb}{0.94,0.94,0.94}
\definecolor{grey95}{rgb}{0.95,0.95,0.95}
\definecolor{grey96}{rgb}{0.96,0.96,0.96}
\definecolor{grey97}{rgb}{0.97,0.97,0.97}
\definecolor{grey98}{rgb}{0.98,0.98,0.98}
\definecolor{grey99}{rgb}{0.99,0.99,0.99}
\definecolor{grey9}{rgb}{0.09,0.09,0.09}
\definecolor{grey}{rgb}{0.75,0.75,0.75}
\definecolor{honeydew1}{rgb}{0.94,1.00,0.94}
\definecolor{honeydew2}{rgb}{0.88,0.93,0.88}
\definecolor{honeydew3}{rgb}{0.76,0.80,0.76}
\definecolor{honeydew4}{rgb}{0.51,0.55,0.51}
\definecolor{honeydew}{rgb}{0.94,1.00,0.94}
\definecolor{hotpink}{rgb}{1.00,0.41,0.71}
\definecolor{indianred}{rgb}{0.80,0.36,0.36}
\definecolor{ivory1}{rgb}{1.00,1.00,0.94}
\definecolor{ivory2}{rgb}{0.93,0.93,0.88}
\definecolor{ivory3}{rgb}{0.80,0.80,0.76}
\definecolor{ivory4}{rgb}{0.55,0.55,0.51}
\definecolor{ivory}{rgb}{1.00,1.00,0.94}
\definecolor{khaki1}{rgb}{1.00,0.96,0.56}
\definecolor{khaki2}{rgb}{0.93,0.90,0.52}
\definecolor{khaki3}{rgb}{0.80,0.78,0.45}
\definecolor{khaki4}{rgb}{0.55,0.53,0.31}
\definecolor{khaki}{rgb}{0.94,0.90,0.55}
\definecolor{lavenderblush}{rgb}{1.00,0.94,0.96}
\definecolor{lavender}{rgb}{0.90,0.90,0.98}
\definecolor{lawngreen}{rgb}{0.49,0.99,0.00}
\definecolor{lemonchiffon}{rgb}{1.00,0.98,0.80}
\definecolor{lightblue}{rgb}{0.68,0.85,0.90}
\definecolor{lightcoral}{rgb}{0.94,0.50,0.50}
\definecolor{lightcyan}{rgb}{0.88,1.00,1.00}
\definecolor{lightgoldenrod}{rgb}{0.93,0.87,0.51}
\definecolor{lightgoldenrod}{rgb}{0.98,0.98,0.82}
\definecolor{lightgray}{rgb}{0.83,0.83,0.83}
\definecolor{lightgreen}{rgb}{0.56,0.93,0.56}
\definecolor{lightgrey}{rgb}{0.83,0.83,0.83}
\definecolor{lightpink}{rgb}{1.00,0.71,0.76}
\definecolor{lightsalmon}{rgb}{1.00,0.63,0.48}
\definecolor{lightsea}{rgb}{0.13,0.70,0.67}
\definecolor{lightsky}{rgb}{0.53,0.81,0.98}
\definecolor{lightslate}{rgb}{0.47,0.53,0.60}
\definecolor{lightslate}{rgb}{0.47,0.53,0.60}
\definecolor{lightslate}{rgb}{0.52,0.44,1.00}
\definecolor{lightsteel}{rgb}{0.69,0.77,0.87}
\definecolor{lightyellow}{rgb}{1.00,1.00,0.88}
\definecolor{limegreen}{rgb}{0.20,0.80,0.20}
\definecolor{linen}{rgb}{0.98,0.94,0.90}
\definecolor{magenta1}{rgb}{1.00,0.00,1.00}
\definecolor{magenta2}{rgb}{0.93,0.00,0.93}
\definecolor{magenta3}{rgb}{0.80,0.00,0.80}
\definecolor{magenta4}{rgb}{0.55,0.00,0.55}
\definecolor{magenta}{rgb}{1.00,0.00,1.00}
\definecolor{maroon1}{rgb}{1.00,0.20,0.70}
\definecolor{maroon2}{rgb}{0.93,0.19,0.65}
\definecolor{maroon3}{rgb}{0.80,0.16,0.56}
\definecolor{maroon4}{rgb}{0.55,0.11,0.38}
\definecolor{maroon}{rgb}{0.69,0.19,0.38}
\definecolor{mediumaquamarine}{rgb}{0.40,0.80,0.67}
\definecolor{mediumblue}{rgb}{0.00,0.00,0.80}
\definecolor{mediumorchid}{rgb}{0.73,0.33,0.83}
\definecolor{mediumpurple}{rgb}{0.58,0.44,0.86}
\definecolor{mediumsea}{rgb}{0.24,0.70,0.44}
\definecolor{mediumslate}{rgb}{0.48,0.41,0.93}
\definecolor{mediumspring}{rgb}{0.00,0.98,0.60}
\definecolor{mediumturquoise}{rgb}{0.28,0.82,0.80}
\definecolor{mediumviolet}{rgb}{0.78,0.08,0.52}
\definecolor{midnightblue}{rgb}{0.10,0.10,0.44}
\definecolor{mintcream}{rgb}{0.96,1.00,0.98}
\definecolor{mistyrose}{rgb}{1.00,0.89,0.88}
\definecolor{moccasin}{rgb}{1.00,0.89,0.71}
\definecolor{navajowhite}{rgb}{1.00,0.87,0.68}
\definecolor{navyblue}{rgb}{0.00,0.00,0.50}
\definecolor{navy}{rgb}{0.00,0.00,0.50}
\definecolor{oldlace}{rgb}{0.99,0.96,0.90}
\definecolor{olivedrab}{rgb}{0.42,0.56,0.14}
\definecolor{orange1}{rgb}{1.00,0.65,0.00}
\definecolor{orange2}{rgb}{0.93,0.60,0.00}
\definecolor{orange3}{rgb}{0.80,0.52,0.00}
\definecolor{orange4}{rgb}{0.55,0.35,0.00}
\definecolor{orangered}{rgb}{1.00,0.27,0.00}
\definecolor{orange}{rgb}{1.00,0.65,0.00}
\definecolor{orchid1}{rgb}{1.00,0.51,0.98}
\definecolor{orchid2}{rgb}{0.93,0.48,0.91}
\definecolor{orchid3}{rgb}{0.80,0.41,0.79}
\definecolor{orchid4}{rgb}{0.55,0.28,0.54}
\definecolor{orchid}{rgb}{0.85,0.44,0.84}
\definecolor{palegoldenrod}{rgb}{0.93,0.91,0.67}
\definecolor{palegreen}{rgb}{0.60,0.98,0.60}
\definecolor{paleturquoise}{rgb}{0.69,0.93,0.93}
\definecolor{paleviolet}{rgb}{0.86,0.44,0.58}
\definecolor{papayawhip}{rgb}{1.00,0.94,0.84}
\definecolor{peachpuff}{rgb}{1.00,0.85,0.73}
\definecolor{peru}{rgb}{0.80,0.52,0.25}
\definecolor{pink1}{rgb}{1.00,0.71,0.77}
\definecolor{pink2}{rgb}{0.93,0.66,0.72}
\definecolor{pink3}{rgb}{0.80,0.57,0.62}
\definecolor{pink4}{rgb}{0.55,0.39,0.42}
\definecolor{pink}{rgb}{1.00,0.75,0.80}
\definecolor{plum1}{rgb}{1.00,0.73,1.00}
\definecolor{plum2}{rgb}{0.93,0.68,0.93}
\definecolor{plum3}{rgb}{0.80,0.59,0.80}
\definecolor{plum4}{rgb}{0.55,0.40,0.55}
\definecolor{plum}{rgb}{0.87,0.63,0.87}
\definecolor{powderblue}{rgb}{0.69,0.88,0.90}
\definecolor{purple1}{rgb}{0.61,0.19,1.00}
\definecolor{purple2}{rgb}{0.57,0.17,0.93}
\definecolor{purple3}{rgb}{0.49,0.15,0.80}
\definecolor{purple4}{rgb}{0.33,0.10,0.55}
\definecolor{purple}{rgb}{0.63,0.13,0.94}
\definecolor{red1}{rgb}{1.00,0.00,0.00}
\definecolor{red2}{rgb}{0.93,0.00,0.00}
\definecolor{red3}{rgb}{0.80,0.00,0.00}
\definecolor{red4}{rgb}{0.55,0.00,0.00}
\definecolor{red}{rgb}{1.00,0.00,0.00}
\definecolor{rosybrown}{rgb}{0.74,0.56,0.56}
\definecolor{royalblue}{rgb}{0.25,0.41,0.88}
\definecolor{saddlebrown}{rgb}{0.55,0.27,0.07}
\definecolor{salmon1}{rgb}{1.00,0.55,0.41}
\definecolor{salmon2}{rgb}{0.93,0.51,0.38}
\definecolor{salmon3}{rgb}{0.80,0.44,0.33}
\definecolor{salmon4}{rgb}{0.55,0.30,0.22}
\definecolor{salmon}{rgb}{0.98,0.50,0.45}
\definecolor{sandybrown}{rgb}{0.96,0.64,0.38}
\definecolor{seagreen}{rgb}{0.18,0.55,0.34}
\definecolor{seashell1}{rgb}{1.00,0.96,0.93}
\definecolor{seashell2}{rgb}{0.93,0.90,0.87}
\definecolor{seashell3}{rgb}{0.80,0.77,0.75}
\definecolor{seashell4}{rgb}{0.55,0.53,0.51}
\definecolor{seashell}{rgb}{1.00,0.96,0.93}
\definecolor{sienna1}{rgb}{1.00,0.51,0.28}
\definecolor{sienna2}{rgb}{0.93,0.47,0.26}
\definecolor{sienna3}{rgb}{0.80,0.41,0.22}
\definecolor{sienna4}{rgb}{0.55,0.28,0.15}
\definecolor{sienna}{rgb}{0.63,0.32,0.18}
\definecolor{skyblue}{rgb}{0.53,0.81,0.92}
\definecolor{slateblue}{rgb}{0.42,0.35,0.80}
\definecolor{slategray}{rgb}{0.44,0.50,0.56}
\definecolor{slategrey}{rgb}{0.44,0.50,0.56}
\definecolor{snow1}{rgb}{1.00,0.98,0.98}
\definecolor{snow2}{rgb}{0.93,0.91,0.91}
\definecolor{snow3}{rgb}{0.80,0.79,0.79}
\definecolor{snow4}{rgb}{0.55,0.54,0.54}
\definecolor{snow}{rgb}{1.00,0.98,0.98}
\definecolor{springgreen}{rgb}{0.00,1.00,0.50}
\definecolor{steelblue}{rgb}{0.27,0.51,0.71}
\definecolor{tan1}{rgb}{1.00,0.65,0.31}
\definecolor{tan2}{rgb}{0.93,0.60,0.29}
\definecolor{tan3}{rgb}{0.80,0.52,0.25}
\definecolor{tan4}{rgb}{0.55,0.35,0.17}
\definecolor{tan}{rgb}{0.82,0.71,0.55}
\definecolor{thistle1}{rgb}{1.00,0.88,1.00}
\definecolor{thistle2}{rgb}{0.93,0.82,0.93}
\definecolor{thistle3}{rgb}{0.80,0.71,0.80}
\definecolor{thistle4}{rgb}{0.55,0.48,0.55}
\definecolor{thistle}{rgb}{0.85,0.75,0.85}
\definecolor{tomato1}{rgb}{1.00,0.39,0.28}
\definecolor{tomato2}{rgb}{0.93,0.36,0.26}
\definecolor{tomato3}{rgb}{0.80,0.31,0.22}
\definecolor{tomato4}{rgb}{0.55,0.21,0.15}
\definecolor{tomato}{rgb}{1.00,0.39,0.28}
\definecolor{turquoise1}{rgb}{0.00,0.96,1.00}
\definecolor{turquoise2}{rgb}{0.00,0.90,0.93}
\definecolor{turquoise3}{rgb}{0.00,0.77,0.80}
\definecolor{turquoise4}{rgb}{0.00,0.53,0.55}
\definecolor{turquoise}{rgb}{0.25,0.88,0.82}
\definecolor{violetred}{rgb}{0.82,0.13,0.56}
\definecolor{violet}{rgb}{0.93,0.51,0.93}
\definecolor{wheat1}{rgb}{1.00,0.91,0.73}
\definecolor{wheat2}{rgb}{0.93,0.85,0.68}
\definecolor{wheat3}{rgb}{0.80,0.73,0.59}
\definecolor{wheat4}{rgb}{0.55,0.49,0.40}
\definecolor{wheat}{rgb}{0.96,0.87,0.70}
\definecolor{whitesmoke}{rgb}{0.96,0.96,0.96}
\definecolor{white}{rgb}{1.00,1.00,1.00}
\definecolor{yellow1}{rgb}{1.00,1.00,0.00}
\definecolor{yellow2}{rgb}{0.93,0.93,0.00}
\definecolor{yellow3}{rgb}{0.80,0.80,0.00}
\definecolor{yellow4}{rgb}{0.55,0.55,0.00}
\definecolor{yellowgreen}{rgb}{0.60,0.80,0.20}
\definecolor{yellow}{rgb}{1.00,1.00,0.00}

\usepackage[psamsfonts]{amsfonts} 

\usepackage{pdflscape}

\usepackage{pgfplots}
\usepackage{dsfont}  

\usepackage{mathtools}

\newcommand\pgfmathsinandcos[3]{%
  \pgfmathsetmacro#1{sin(#3)}%
  \pgfmathsetmacro#2{cos(#3)}%
}

\newcommand\LongitudePlane[3][current plane]{%
  \pgfmathsinandcos\sinEl\cosEl{#2} 
  \pgfmathsinandcos\sint\cost{#3} 
  \tikzset{#1/.estyle={cm={\cost,\sint*\sinEl,0,\cosEl,(0,0)}}}
}
\newcommand\LatitudePlane[3][current plane]{%
  \pgfmathsinandcos\sinEl\cosEl{#2} 
  \pgfmathsinandcos\sint\cost{#3} 
  \pgfmathsetmacro\yshift{\cosEl*\sint}
  \tikzset{#1/.estyle={cm={\cost,0,0,\cost*\sinEl,(0,\yshift)}}} %
}
\newcommand\DrawLongitudeCircle[2][1]{
  \LongitudePlane{\angEl}{#2}
  \tikzset{current plane/.prefix style={scale=#1}}
  \pgfmathsetmacro\angVis{atan(sin(#2)*cos(\angEl)/sin(\angEl))} %
  \draw[current plane, blue] (\angVis:1) arc (\angVis:\angVis+180:1);
  \draw[current plane,dashed] (\angVis-180:1) arc (\angVis-180:\angVis:1);
}

\newcommand\DrawLatitudeCircle[2][1]{
  \LatitudePlane{\angEl}{#2}
  \tikzset{current plane/.prefix style={scale=#1}}
  \pgfmathsetmacro\sinVis{sin(#2)/cos(#2)*sin(\angEl)/cos(\angEl)}
  \pgfmathsetmacro\angVis{asin(min(1,max(\sinVis,-1)))}
  \draw[current plane, orange] (\angVis:1) arc (\angVis:-\angVis-180:1);
  \draw[current plane,dashed] (180-\angVis:1) arc (180-\angVis:\angVis:1);
}

\newcommand{\mockalph}[1]{}

\newcommand{\mbf}[1]{\mbox{\boldmath$#1$\unboldmath}}

\newcommand{\dd}{{\mathrm d}}

\newcommand{\cqfd}{\hfill $\square$}

\newcommand{\R}{\mathbb R}

\newcommand{\n}{^{(n)}}

\newcommand{\Xb}{\mathbf{X}}
\newcommand{\Sb}{\mathbf{S}}
\newcommand{\sbb}{\mathbf{s}}

\newcommand{\Jb}{\mathbf{J}}

\newcommand{\Qb}{\mathbf{Q}}

\newcommand{\Zb}{\mathbf{Z}}

\newcommand{\Fb}{\mathbf{F}}

\newcommand{\ub}{\ensuremath{\mathbf{u}}}

\newcommand{\vb}{\ensuremath{\mathbf{v}}}
\newcommand{\xb}{\ensuremath{\mathbf{x}}}
\newcommand{\yb}{\ensuremath{\mathbf{y}}}
\newcommand{\zb}{\ensuremath{\mathbf{z}}}

\newcommand{\Db}{\ensuremath{\mathbf{D}}}
\newcommand{\Ub}{\ensuremath{\mathbf{U}}}

\newcommand{\Yb}{\ensuremath{\mathbf{Y}}}

\newcommand{\Gb}{\ensuremath{\mathbf{G}}}

\newcommand{\thetab}{{\pmb \theta}}

\newcommand{\Gamb}{{\pmb \Gamma}}

\newcommand{\bth}{{\boldsymbol \theta}}

\newcommand{\tbomega}{{\scriptstyle{\boldsymbol{\omega}}}}
\newcommand{\tbtau}{{\scriptstyle{\boldsymbol{\tau}}}}
\newcommand{\tbnu}{{\scriptstyle{\boldsymbol{\nu}}}}

\DeclareMathAlphabet\mathbfcal{OMS}{cmsy}{b}{n}

\def\bth{\mbox{\boldmath$\theta$}}

\def\bgam{\mbox{\boldmath$\gamma$}}
\def\btau{\mbox{\boldmath$\tau$}}
\def\bnu{\mbox{\boldmath$\nu$}}
\def\bDelta{\mbox{\boldmath$\Delta$}}

\def\bGamma{\mbox{\boldmath$\Gamma$}}

\def\bOmega{\mbox{\boldmath$\Omega$}}

\def\bvp{\mbox{\boldmath$\varphi$}}

\newcommand{\pr}{^{\prime}}

\newcommand{\thetatM}{{\widehat\thetab}^{\phantom n}_{\rm M}{\phantom\theta}\hspace{-3.9mm}^{{(n)}}}
\newcommand{\thetat}{{\widehat\thetab}^{\phantom n}{\phantom\theta}\hspace{-3.0mm}^{{(n)}}}
\newcommand{\thetatn}{{\widehat\thetab}^{\phantom n}{\phantom\theta}\hspace{-2.5mm}^{{(n)}}}

\newcommand{\thetatsub}{{\,\widehat\thetab}^{\; \phantom n}{\phantom\theta}\hspace{-3.0mm}^{{(n)}} }

\newcommand{\ny}{n\rightarrow\infty}

\title{Nonparametric Measure-Transportation-Based Methods for Directional Data}

\author[1\authfn{1}]{Hallin M.}
\author[2\authfn{2}]{Liu H.}
\author[1\authfn{3}]{Verdebout T.}

\affil[1]{ECARES and D\'{e}partement de Math\'{e}matique,
Universit\'{e} libre de Bruxelles, Brussels, Belgium}
\affil[2]{International Institute of Finance, School of Management, University of Science and Technology of China, Hefei, 230026, Anhui Province, P.R.~China}

\corraddress{Marc Hallin, Universit\'{e} libre de Bruxelles, ECARES and D\'{e}partement de Math\'{e}matique,  B-1050 Brussels, Belgium}
\corremail{mhallin@ulb.ac.be}

\fundinginfo{$^*$Marc Hallin gratefully acknowledges the support of the COST (European Cooperation in Science and Technology) Action HiTEc CA21163 and the Czech Science Foundation grants GA\v{C}R22036365  and GA24-100788. 

$^\dagger$Hang Liu's research is supported by the USTC grants WK2040000055 and YD2040002016. 

$^\ddagger$Thomas Verdebout is supported by an ARC grant of the ULB and the Communaut\'e Fran{\c c}aise de Belgique, a Projet de Recherche (PDR) grant of the FNRS, and the Fonds Thelam of the Fondation Roi Baudouin. }

\runningauthor{Hallin, Liu, and Verdebout}

\begin{document}

\maketitle

\begin{abstract}

This paper proposes various nonparametric tools  based on measure transportation for directional data. We use optimal transports to define  new notions of  distribution and quantile  functions on the hypersphere, with meaningful quantile contours and regions and  closed-form formulas  under the classical assumption of rotational symmetry. The empirical versions of our distribution functions enjoy the expected Glivenko-Cantelli property of traditional distribution functions. They provide fully distribution-free concepts of ranks and signs and define data-driven systems of (curvilinear) parallels and (hyper)meridians.  Based on this, we also construct a universally consistent test of uniformity and a class of  fully distribution-free and universally consistent tests for directional MANOVA which, in simulations,  outperform all their existing competitors.  A real-data example involving the analysis of sunspots concludes the paper.    

\keywords{directional statistics, directional quantiles, ranks, signs, optimal transport, directional GOF, directional MANOVA}
\end{abstract}


\section{Introduction}

 Directional data analysis is dealing with random directions---taking values on  circles and  (hyper)spheres or, more generally, on manifolds. Directional data can be found in a variety of fields, including astronomy \citep{Marinucci2008, Marinucci2011}, environmetrics \citep{Garcia-Portugues2014, Ameijeiras2018, Kume2018}, biology and medicine \citep{Dryden2005, Hamelryck2006, Dortet-Bernadet2008}, to cite only a few. Directions generally are represented as points on the unit hypersphere~$ \mathcal{S}^{d-1}\coloneqq  \{ \textbf{x} \in {\mathbb R}^d:  \|\textbf{x} \|^2\coloneqq  \xb^\top \xb =1 \}$ in~${\mathbb R}^d$. The special nature of their sample spaces gives directional statistics a specific flavor. The models of interest are of the standard statistical types (one-sample and  multi-sample location, analysis of variance, regression, etc.) and inference is carried out according to the usual principles, but their   implementation takes  special forms and requires distinctive techniques. A comprehensive exposition of the theoretical background and inference methods for circular and spherical data 
  can be found in \cite{Mardia1999a} or \cite{Jammalamadaka2001}. More recent advances are presented in  \cite{Ley2017a, Ley2018}, while statistical tools for data on more general spaces are studied in \cite{Chikuse2003} and \cite{BB08, BB12}.

The most popular parametric model in directional statistics, which can be traced back to the first decades of the~20th century, is the von Mises-Fisher model, characterized by the family of von Mises distributions, with densities    (with respect to the classical surface area measure on $ \mathcal{S}^{d-1}$) of the form\vspace{-1mm}
\begin{equation} \label{vMFdens}
\zb \in \mathcal{S}^{d-1}  \mapsto c_\kappa {\rm exp}( \kappa \zb^\top \thetab),\vspace{-1mm}
\end{equation}
where $\thetab \in \mathcal{S}^{d-1}$ plays the role of a location parameter, $\kappa \in\R^+$, which drives the probability mass in the vicinity of~$\thetab$, is the so-called {\it concentration parameter}, $c_\kappa$ is a normalizing constant, and $\zb^\top$ denotes the transpose of $\zb$. The von Mises distribution is often seen as the `` Gaussian directional distribution'' due to the fact that  the maximum likelihood estimator of $\thetab$ is the normalized sample average. 

The von Mises distributions belong to the more general class of rotationally symmetric distributions, which contains all distributions  with densities (still with respect to the surface area measure on $\mathcal{S}^{d-1}$) of the form\vspace{-1mm}
\begin{equation} \label{Rotsymfirst}
\zb \in \mathcal{S}^{d-1} \mapsto c_{f} f(\zb^\top \thetab), \vspace{-1mm}
\end{equation}
where $\thetab \in \mathcal{S}^{d-1}$, $f$ is some positive {\it  angular function}, and $c_f$ is a norming constant. 
 The projection $({\bf I}_p- \thetab \thetab^\top) \Zb$ of a rotationally symmetric random vector $\Zb$ onto the tangent space (to $ \mathcal{S}^{d-1}$) at $\thetab$ has a spherical distribution. Inference for the parameters of rotationally symmetric distributions has  been considered recently  in \cite{Christie2015}, \cite{Kanika2015}, and \cite{Paindaveine2020b, Paindaveine2020}.  Extensions of rotational symmetry  yielding, after projection  onto the tangent space  (to $ \mathcal{S}^{d-1}$) at $\thetab$,  elliptically symmetric rather than spherical distributions  have been proposed in \cite{Kent1982},  \cite{Scealy2019}, and \cite{Garcia-Portugues2020}; see also \cite{Kume2013}, \cite{Kume2018}, and \cite{Kent2018a}.  
 
 {Rank-based methods (including MANOVA) for directional data have been proposed in \cite{Ley2013, Ley2017} and \cite{Verdebout2017}}. While the concepts of ranks considered there enjoy several attractive features, their distribution-freeness, hence also their applicability, is limited to the class of  rotationally symmetric distributions. Based on sample projections onto the mean direction, the quantiles studied in \cite{Ley2014} are  canonical for rotationally symmetric distributions on $\mathcal{S}^{d-1}$ but fail to be informative for general distributions on $\mathcal{S}^{d-1}$. 
 
 The assumptions of rotational symmetry  and their elliptical extensions are pervasive in the literature on directional data and   can be compared, in this respect, to the assumptions of spherical and  elliptical  symmetry in traditional multivariate analysis. Just as the latter, they drastically simplify the mathematical structure of inference problems but also very severely restrict their validity: in most applications, rotational symmetry, indeed, is extremely unlikely to hold.  The situation in the directional context, thus, is quite comparable to the situation in the general multivariate context: the absence of a canonical ordering, at first sight, precludes the canonical definition of essential order-related statistical tools as ranks, signs, and quantiles unless one is ready to make a very strong assumption of rotational symmetry (or its ``elliptical extensions") in the directional case, of spherical or elliptical symmetry in the general multivariate case. Based on measure transportation ideas and novel concepts of multivariate distribution and quantile functions in ${\mathbb R}^d$, new non- and semi-parametric methods have been proposed very recently in \cite{chernoetal17} and \cite{Hallin2021}, which have demonstrated their efficiency in a variety of  multivariate inference problems: see \cite{deb2021efficiency}, \cite{HMordantS}, \cite{GhosalSenAOS}, \cite{HHH,HLH2022,HLH2022+}, \cite{shi2019distribution, shi2022}, \cite{DGM22}, \cite{Bodi2021}.  \cite{Rigollet2022} use optimal transport to perform causal inference, while machine learning methods based on optimal transport have been proposed, e.g., in \cite{kolouri2017}. We refer to~\cite{MHAnnR} for a nontechnical survey.  

The objective of this paper is to develop a measure-transportation-based approach for directional data similar to the one adopted in the above references, all of which  are dealing with observations in ${\mathbb R}^d$. Measure transportation on hyperspheres and more general manifolds, however,  is not a cosmetic adaptation of the Euclidean case, as it involves non-Euclidean topologies and  metrics. As a consequence, the measure-transportation-based concepts and methods developed here do not follow along the same lines as  in \cite{chernoetal17} and \cite{Hallin2021}. Below, we mainly build upon \cite{McCann2001}, where the celebrated  results of \cite{McCann1995} on measure-preserving monotone mappings in~${\mathbb R}^d$ are extended to Riemanniann manifolds, to propose directional concepts of distribution and quantile functions inducing a distribution-specific system of curvilinear parallels and (hyper)meridians. These directional distribution functions can be seen as canonical transformations in the sense of \cite{Jupp2020}. 
The empirical counterparts of these   concepts yield directional ranks, signs, and empirical quantiles enjoying the   properties one is expecting from such notions (in particular, distribution-freeness of ranks and signs and quantile regions with prescribed probability contents irrespective of the underlying distribution); empirical distribution functions, moreover, are shown to satisfy a Glivenko-Cantelli consistency property. 

The motivation for these definitions are the same as in ${\mathbb R}$ and ${\mathbb R}^d$: conducting distribution-free inference based on ranks and signs, constructing quantile contours and (closed, connected, and nested) quantile regions, 
detecting directional outliers, defining directional values at risk,  performing directional quantile regression, etc.  Earlier attemps have been made with similar objectives, which are only partially successful: the notion of ranks studied in  \cite{Ley2013, Ley2017} and \cite{Verdebout2017}, indeed, is not enjoying  distribution-freeness unless the underlying distribution is rotationally symmetric; similarly,  the quantile contours proposed in \cite{Ley2014} are satisfactory under rotational symmetry only.

In order to show that our measure-transportation-based concepts, contrary to   earlier proposals, are achieving these objectives, we consider two very classical problems in directional inference:
\begin{compactenum}
\item[(i)] (testing for directional uniformity) the classical problem of testing uniformity on $\mathcal{S}^{d-1}$. This GOF problem is probably the oldest and most fundamental one in directional statistics, and can be traced back to the discussion by \cite{Bernoulli1735} on  whether the closeness of the orbital planes of various planets arose by chance or not. To cite only a few recent works, tests of uniformity in the context of noisy directional data were studied in \cite{Lacour2014} and \cite{Kim2016} and, in a high-dimensional context, in \cite{Cai2012}, \cite{Cai2013}, and \cite{Cutting2017}. Projection-based tests were proposed by \cite{Cuesta-Albertos2009} and \cite{Garcia-Portugues2020b}. Tests of uniformity were used in \cite{Garcia-Portugues2020} in the problem of testing for rotational symmetry. Inference for location in the vicinity of   uniformity  was considered in \cite{Paindaveine2017}. We refer to \cite{Garcia-Portugues2020a} for a recent review of this topic. 
\item[(ii)] (directional MANOVA) MANOVA on the hypersphere---testing the hypothesis of no-treatment effect, where the treatment can impact   location, concentration, or skewness, is another very  classical problem in directional inference. Due to the absence of distribution-free   tools, the traditional approach typically reduces to  a pseudo-von Mises MANOVA procedure---see for instance the recent contributions by  \cite{figueiredo2006}, \cite{Ley2017}, \cite{sengupta2020}, and \cite{kulkarni2022}, and the references therein; a major weakness of  pseudo-von Mises methods, however, is that their validity, which requires the classical Fisher consistency property, is guaranteed under rotational symmetry only.  
\end{compactenum}
In both cases, our directional measure-transportation-based concepts   provide natural nonparametric and fully distri\-bution-free solutions which (in simulations) perform equally well as some of their competitors under rotational symmetry and, most importantly, outperform them all in the non-rotationally-symmetric case: 
\begin{compactenum}
\item[(i)]  (testing for directional uniformity)  the directional distribution function for the uniform over $\mathcal{S}^{d-1}$ is the identity function; a (fully distribution-free) Cram\' er-von Mises-type test of uniformity on $\mathcal{S}^{d-1}$ based on our concept of empirical directional distribution function is shown to be {\it universally} consistent while outperforming (in simulations)  its traditional  competitors (the Rayleigh,  Bingham,   Ajne,   Gin\'e, and Bakshaev tests), as well as the more recent  projection-based tests by \cite{Garcia-Portugues2020b}. In sharp contrast with the universal consistency of our measure-transportation-based Cram\' er-von Mises test, these competitors, moreover, typically have undetectable {\it blind spots}. Also note that our measure-transportation-based Cram\' er-von Mises tests straightforwardly extend from testing uniformity  to the general GOF problem;  
%
%
\item[(ii)] (directional MANOVA) the class of fully distribution-free MANOVA tests we are proposing is based on our concepts of directional ranks and signs; as in the traditional case of univariate rank tests, it involves    the choice of a score function, which can target various types of alternatives. 
 We provide the asymptotic ((non)central chi-square) distributions of the test statistics under the null and under local alternatives  and show via simulations  how our tests outperform their  classical  pseudo-von Mises competitors (the validity of which, moreover, is limited to the rotationally symmetric case).  
\end{compactenum}

The paper is organized as follows. In Section \ref{sec:Wasser}, we provide the relevant mathematical notions and briefly discuss the concept of measure transportation on Riemannian manifolds and hyperspheres. We define the population versions~$\bf F$ and $\bf Q$ of our directional distribution  and quantile functions in Section~\ref{sec.PopQuan}, then  derive their explicit forms in the rotationally symmetric case. In Section~\ref{sec:emp}, we introduce the empirical version~${\bf F}\n$ of $\bf F$ and establish its   Glivenko-Cantelli asymptotic behavior by adapting the results of \cite{Hallin2021} to the directional setup. Section~\ref{numerillsec}   illustrates the construction of empirical quantile contours in  simulated data. We also introduce new ranks and signs that are adapted the directional setup. Our quantile contours, ranks and signs require a construction which is specific to the non-Euclidean context considered here (and  therefore differs from the classical $\R^d$ case). 
  Under rotational symmetry, however, our quantile contours coincide with those  of \cite{Ley2014}, which are canonical for  rotationally symmetric distributions; some distributional and equivariance properties our signs and ranks are discussed in Supplementary Material A.    In Section~\ref{sec:gof},   we  construct  totally new Cram\' er-von Mises-type tests for  goodness-of-fit based on~${\bf F}\n$, establish their universal consistency against fixed alternatives,  and compare their performance with that of their many competitors. In Section~\ref{MANOVAsec}, we propose a class of distribution-free MANOVA testing procedures and study their asymptotic distributions under the null and under local alternatives using the classical technical tools of \cite{HHH, HLH2022, HLH2022+}; we   also establish their consistency   under non-local alternatives. A real-data application is presented in Section~\ref{sec:real}, where we implement our MANOVA procedure on sunspots data. Finally, Supplementary Material A,~B,~C, and D are collecting (i) the basic properties of our ranks and signs, (ii) additional simulation results, (iii) an illustration of our quantile contours on a protein dataset, and (iv) the proofs of the various results provided in the paper.


\section{Measure transportation on Riemannian manifolds} \label{sec:Wasser}
\setcounter{equation}{0}

Throughout the paper, we  denote by~$d(\yb, \zb)\coloneqq\left\vert\arccos (\yb^\top \zb)\right\vert$ the geodesic distance  and by $c(\yb, \zb)\coloneqq d^2(\yb, \zb)/2$ the squared Riemannian distance between two points $\yb$ and $\zb$ on the unit hypersphere~$ \mathcal{S}^
{d-1}$;  both $c$ and $d$ are continuous and bounded. When  equipped with the geodesic distance,~$ \mathcal{S}^
{d-1}$ is a  separable complete metric space, hence a Polish metric space,   with Borel $\sigma$-field~${\cal B}^{d-1}$. Denote by $\lambda^{d-1}$ the  {\it  surface area measure} over $ \mathcal{S}^
{d-1}$.  Throughout, measurability tacitly is understood with respect to~${\cal B}^{d-1}$ and  ${\rm P}$ is assumed to be in the family~$\mathfrak{P}_d$  of~$\lambda^
{d-1}$-absolutely continuous distributions with densities bounded away from zero and infinity.

Measure transportation on Polish spaces and, more particularly, on Riemannian mani\-folds, is a well-studied subject; below, we mainly follow \cite{McCann2001}, where the celebrated  results of \cite{McCann1995} on measure-preserving monotone mappings in~${\mathbb R}^d$ are extended to Riemanniann manifolds, \cite{AmbrosioPratelli2003}, and \cite{SCHACHERMAYERTEICHMANN2008}. Let ${\rm P}$ and~${\rm Q}$ denote two probability measures on $ \mathcal{S}^{d-1}$.    Monge's optimal transport problem on ${\cal S}^{d-1} $ consists in minimizing, among the set  ${\mathcal S}({\rm P}, {\rm Q})$ of all  measurable transport maps ${\bf S}\!:  \mathcal{S}^{d-1}\! \to\!  \mathcal{S}^{d-1}$  such\linebreak  that~$
 (\Sb \# {\rm P})(V)\coloneqq 
{\rm P}(\Sb^{-1}(V)) 
= {\rm Q}({V})$   for all~$V\!\in\!{\cal B}^{d-1}$ (in the measure transportation notation and termino-\linebreak logy,~$\Sb\!\in\!{\mathcal S}({\rm P}, {\rm Q})$ is   {\it  pushing  ${\rm P}$ forward   
 to}~${\rm Q}$), the transportation~cost
\begin{equation} \label{cost}
C_{\rm M}({\bf S})\coloneqq  \int_{{\cal S}^{d-1}} c\big(\zb, {\bf S}(\zb)\big) \dd {\rm P}(\zb)= {\rm E}_{\rm P}\big[c(\Zb, {\bf S}(\Zb))\big].
\end{equation}
Due to the fact that ${\mathcal S}({\rm P}, {\rm Q})$ is not convex, Monge's problem is an uneasy one.

%

Closely related to Monge's problem 
 is   the so-called Kantorovich problem. 
 Denote by~$\Gamma({\rm P}, {\rm Q})$ the set of all probability measures $\gamma$ on $ \mathcal{S}^{d-1}\! \times~\!\mathcal{S}^{d-1}\!$ with ``marginals" ${\rm P}$ and~${\rm Q}$;  an element $\gamma$ of $\Gamma({\rm P}, {\rm Q})$ is called a {\it transport plan}. 
The Kantorovich  problem on~${\cal S}^{d-1} $  consists in minimizing, among all  transport  plans~$\gamma \in\Gamma({\rm P}, {\rm Q})$, the expectation\vspace{-1mm}
  \begin{equation}\label{costK}  C_{\rm K}(\gamma)\coloneqq \int_{ \mathcal{S}^{d-1} \times  \mathcal{S}^{d-1}} c(\yb, \zb) \dd\gamma(\yb, \zb)\vspace{-1mm}
\end{equation}
of the transportation cost $c(\Yb,\Zb)$ under $(\Yb,\Zb)\sim\gamma$.
Unlike ${\mathcal S}({\rm P}, {\rm Q})$, the set~$\Gamma({\rm P}, {\rm Q})$ is a convex subset of a Banach space;  since $c$ is bounded,  the existence of a solution  is guaranteed. 

The Kantorovich problem is a relaxation of Monge's. Letting $({\bf I}_{\rm d}~\!\times~\!{\bf S}){\bf Z}\coloneqq ~\!({\bf Z}, {\bf S}({\bf Z}))$, a  transport map~${\bf S} \in {\mathcal S}({\rm P}, {\rm Q})$ indeed  induces a transport plan $\gamma=({\bf I}_{\rm d} \times {\bf S}) \# {\rm P}\in \Gamma({\rm P}, {\rm Q})$ which is concentrated on the graph of $\Sb$ (i.e., such that~$\gamma\big(({\bf Z}, {\bf S}({\bf Z}))\in\big\{(\zb,\Sb(\zb)): \zb\in \mathcal{S}^{d-1}\}\big)~\!\!=~\!\!1$) and satisfying $C_{\rm K}(\gamma)=C_{\rm M}(\Sb)$. As we shall see (Proposition~\ref{keylemma} below), the solution of Kantorovich's problem is precisely of the form $({\bf I}_{\rm d} \times {\bf F}) \# {\rm P}$ where~$ {\bf F}$, thus, is a solution of Monge's problem, and\vspace{-1mm}  $$
C_{\rm M}(\Fb)=\min_{\Sb\in{\cal S}({\rm P},{\rm Q})}C_{\rm M}(\Sb) = \min_{\gamma\in{\Gamma}({\rm P},{\rm Q})}C_{\rm K}(\gamma) 
= C_{\rm K}(({\bf I}_{\rm d} \times {\bf F}) \# {\rm P})
,\vspace{-1mm}$$ 
where $\min_{\gamma\in{\Gamma}({\rm P},{\rm Q})}C^{1/2}_{\rm K}(\gamma) \eqqcolon{\cal W}_2({\rm P},{\rm Q})$ is the Wasserstein distance between~${\rm P}$ and ${\rm Q}$.

Before providing  more precise statements about the solutions of the Monge and Kantorovich problems,   let us introduce some notation from differential geometry. Recall  that  the {\it tangent space} (tangent to $ \mathcal{S}^{d-1}$) at  $\xb \in  \mathcal{S}^{d-1}$ is the $(d-1)$-dimensional linear subspace 
${\cal T}_\xb \mathcal{S}^{d-1}\coloneqq  \{ \zb \in \R^d, \zb^\top \xb=0 \}.$
 Letting $\xb \in  \mathcal{S}^{d-1}$ and ${\bf v} \in {\cal T}_\xb \mathcal{S}^{d-1}$, the directional derivative at $\xb \in   \mathcal{S}^{d-1}$ of a smooth function $\varsigma:  \mathcal{S}^{d-1} \to \R$   in   direction~${\bf v}$ is defined as 
${v}({\varsigma})\coloneqq  (\varsigma \circ \nu)^\prime (0),$
 where~$\nu~\!:~\![0,1]~\!\to~\!\mathcal{S}^{d-1}$  
 is a differentiable path   such that~$\nu (0)= \xb$ and~$\nu\pr (0)= {\bf v}$.   The gra\-dient~$ \nabla \varsigma(\xb)$ of~$\varsigma$ at~$\xb \in  \mathcal{S}^{d-1}$  then is defined as the vector in~${\cal T}_\xb \mathcal{S}^{d-1}$ such that
 ${\bf v}^\top \nabla \varsigma(\xb)= {v}({\varsigma})$
 for all ${\bf v} \in {\cal T}_\xb \mathcal{S}^{d-1}$ and $\xb \in  \mathcal{S}^{d-1}$. The exponential map at a point $\xb \in~\!\mathcal{S}^{d-1}$ is a map from the tangent space ${\cal T}_\xb \mathcal{S}^{d-1}$ to~$ \mathcal{S}^{d-1}$. More precisely, {denoting by~$\nu_{\bf v}(t)$  the unique geodesic with tangent vector ${\bf v} \in {\cal T}_\xb \mathcal{S}^{d-1}$ running through $\xb \in \mathcal{S}^{d-1}$ ($\nu_{\bf v}(0)=\xb$), then~$\exp_{\xb} ({\bf v})\coloneqq  \nu_{\bf v} (1)$.} Finally, consider the following concepts of {\it $c$-cyclical monotonicity} and {\it $c$-convexity}. 
\begin{definition}
A subset $\Omega \in {\cal S}^{d-1} \times {\cal S}^{d-1}$ is called {{\it $c$-cyclically monotone}} ($c$ the squared Riemannian distance on ${\cal S}^{d-1}$)  if, denoting by $\Sigma(k)$  the set of permutations of~$\{1,\ldots,k\}$, for all $k \in {\mathbb N}$, all  $\sigma \in \Sigma(k)$, and all $(\xb_1, \yb_1), \ldots, (\xb_k, \yb_k)\in \Omega$, \vspace{-1mm}
$$\sum_{i=1}^k c(\xb_i, \yb_i) \leq \sum_{i=1}^k c(\xb_{\sigma(i)}, \yb_i).$$
\end{definition}
\begin{definition}
Two functions $\psi$ and $\phi$ from  $\mathcal{S}^{d-1}$ to $\R$ such that  \begin{equation}\label{cconv}
\psi(\xb)=\!\! \sup_{\yb \in   \mathcal{S}^{d-1}} \{ \phi(\yb) -c(\xb, \yb) \} \text{ and }\ 
 \phi(\yb)=\!\! \inf_{\xb \in   \mathcal{S}^{d-1}} \{ \psi(\xb) + c(\xb, \yb) \} \end{equation}
are called {\it $c$-convex} and {\it $c$-concave}, respectively;  call $\phi$ (resp. $\psi$) the {\it  c-transform} of $\psi$ (resp. $\phi$).
\end{definition}
Since the squared Riemannian distance $c$ satisfies the conditions of Theorem 10.26(i) of \cite{Villani2009}, the functions~$\psi$ and $\phi$ in \eqref{cconv} are a.e.\ differentiable. The following theorem summarizes, for the Monge and Kantorovich problems with transportation cost the squared Riemannian distance   on ${\cal S}^{d-1}$, some of the results contained in Chapter 5 of \cite{Villani2009} (for the necessity part of (i), first established by \cite{Rusch1996}), Theorem~1 of \cite{SCHACHERMAYERTEICHMANN2008} (for the sufficiency part of (i); see also \cite{Pratelli2008}),   Theorems~8 and~9 and Corollary~10 of  \cite{McCann2001} (for (ii), (iii), and (iv)). The assumptions made in these references are automatically satisfied here, since~${\cal S}^{d-1}$ equipped with the geodesic distance is a Polish metric space and $c$ is a continuous bounded cost function. 
%
%
%
\begin{proposition} \label{keylemma} Let $\rm P \in \mathfrak{P}_d$ and $\rm Q$ denote two probability measures on~$\mathcal{S}^{d-1}$. Then,  
\begin{compactenum}
  \item[(i)] a transport plan $\gamma\in\Gamma({\rm P},{\rm Q})$ is the solution (minimizing $C_{\rm K}$  in \eqref{costK} over $\Gamma ({\rm P}, {\rm Q})$) of the Kantorovich problem if and only if it is supported on a    $c$-cycli\-cally monotone subset of~$ {\cal S}^{d-1} \times {\cal S}^{d-1}$; 
 \item[(ii)] this solution of the Kantorovich problem exists, is unique, and is of  the   form~$({\bf I}_{\rm d} \times {\bf F}) \# {\rm P}$ where ${\bf F}\in{\mathcal S}({\rm P},{\rm Q})$ is the~$\rm P$-a.s.\ unique solution of the corresponding Monge problem (minimizing  $C_{\rm M}$   in \eqref{cost}  over ${\mathcal S} ({\rm P}, {\rm Q})$); 
 \item[(iii)] there exist $c$-concave differentiable mappings 
   $\psi$  from~$\mathcal{S}^{d-1}$ to $\R$ such that
  ${\bf F}(\xb)= \exp_{\xb}(-\nabla \psi (\xb))$,  
   $\lambda^{d-1}\text{-a.e.}$
   \end{compactenum}
 If, moreover, $\rm Q \in \mathfrak{P}_d$, then
 \begin{compactenum}
 \item[(iv)]  ${\bf Q}(\xb)\coloneqq  \exp_{\xb}(-\nabla \phi (\xb))$, with $\phi$ the $c$-transform of $\psi$, belongs to ${\mathcal S}({\rm Q}, {\rm P})$, is the~$\rm Q$-a.s. unique minimizer of $C_{\rm M}$ in~\eqref{cost} over ${\mathcal S} ({\rm Q}, {\rm P})$  (i.e., the $\rm Q$-a.s.\  unique solution of the  corresponding  Monge problem),  and  satisfies 
 $${\bf Q}({\bf F}(\xb)) = \xb\quad \rm P\text{-a.s.\quad and \quad} {\bf F}({\bf Q}(\xb)) = \xb\quad \rm Q\text{-a.s.} $$
   \end{compactenum}
\end{proposition}

Proposition \ref{keylemma} entails the existence of a mapping ${\bf F}$   pushing  any probability measure~${\rm P} \in \mathfrak{P}_d$ 
 forward to any other probability measure ${\rm Q} \in \mathfrak{P}_d$ 
  and  minimizing, among all such mappings, the transportation cost \eqref{cost}.  In the next Section, we will use this optimality result   to define  directional distribution and quantile functions.

Our   results all require  densities that are bounded away from zero and $\infty$. It would be interesting to see if our results and our methods could be extended  (i) to densities supported on a subset of the hypersphere and (ii) to other Riemanian manifolds. Encouraging results on the uniform consistency of empirical transports under (i) have been obtained  in \cite{delbarrio2020regularity},  \cite{delbarrio2023regularity}, and \cite{segers2023}. All these results are limited to the Euclidean case, though, and extending them to the hypersphere is beyond the scope of this paper. As for (ii),  whether   our methodology adapts to more general Riemannian manifolds, and how,  is a challenging question which, definitely, should be left for future research.


\section{Directional distribution and  quantile functions}\label{sec.PopQuan} \setcounter{equation}{0}
\subsection{Definitions and some basic properties}\label{31sec}
The concepts of directional distribution and quantile functions we are proposing are based on the optimal transport~${\bf F}$ from a generic probability measure ${\rm P} \in \mathfrak{P}_d$ to the uniform probability measure~${\rm P}^{\bf U}$ on $ \mathcal{S}^{d-1}$.  The idea is inspired by the so-called {\it center-outward distribution}  and {\it quantile functions}  ${\Fb_\pm}$ and~${\Qb_\pm}\coloneqq \Fb_\pm^{-1}$ introduced in 
\cite{Hallin2021} for distributions on $\R ^d$, which  involve optimal transports to the spherical uniform over the unit ball. 
 On the hypersphere, the definition of the  {\it directional distribution function} of a random vector $\Zb\sim{\rm P}^\Zb \in \mathfrak{P}_d$ with values on ${\cal S}^{d-1}$ is 
  as follows (
  we indifferently use the expressions {\it directional distribution and quantile functions of }$\Zb$ and {\it directional distribution and quantile functions of}  ${\rm P}^\Zb$).
 \begin{definition} \label{defdistribution}
 Call  {\it (directional) distribution function of} ${\rm P}^\Zb \in \mathfrak{P}_d$ the a.s.~unique optimal transport map $\Fb$ from $ \mathcal{S}^{d-1}$ to $ \mathcal{S}^{d-1}$  such that $\Fb\# {\rm P}^\Zb={\rm P}^{\bf U}$. 
  \end{definition}
  This definition   is very natural, since   $\Fb$  is such that, for $\Zb \sim {\rm P}^\Zb$,  $\Fb(\Zb) \sim {\rm P}^{\bf U}$---so that $\bf F$ has the same probability integral  transformation flavor as   classical univariate distribution functions. We then have the following important property  (see Section~\ref{Sec.Proofs} of the  Supplementary Material for a proof).
\begin{proposition}\label{ProHom}  The distribution  function $\Fb$ of $\Zb\sim{\rm P}^\Zb \in \mathfrak{P}_d$ is a homeomorphism from $ \mathcal{S}^{d-1}$ to $ \mathcal{S}^{d-1}$.
\end{proposition}
 It directly follows from Proposition~\ref{ProHom} that the inverse~$\Qb\coloneqq \Fb^{-1}$ of $\Fb$ is well defined, which yields the following definition of a {\it directional quantile function}.
 \begin{definition} Call  $\Qb\coloneqq \Fb^{-1}$ the {\it directional quantile function of} ${\rm P}^\Zb\in \mathfrak{P}_d$.\vspace{-2mm}
  \end{definition}
  

In $\R^d$,  the {\it center-outward distribution function} $\Fb_\pm$  was defined as an optimal transport to the spherical uniform~${\rm U}_d$ over the unit ball of $\R^d$. For  ${\rm U}_d$, nested balls of the form~$\tau{\mathcal S}^{d-1}$,  centered at the origin and with ${\rm U}_d$-probability contents~$\tau$, naturally play  the role of quantile regions of order $\tau\in [0,1)$; the origin thus naturally qualifies as the median of~${\rm U}_d$. The center-outward quantile regions then are obtained as the images by~$\Qb_\pm\coloneqq \Fb_\pm^{-1}$ of these balls. 

The basic idea is quite similar  in $ \mathcal{S}^{d-1}$.  If, however,  nested regions with ${\rm P}^{\bf U}$-probability contents $\tau\in [0,1]$   are to be defined on~$ \mathcal{S}^{d-1}$, a central point or {\it pole}~$\thetab\in\mathcal{S}^{d-1}$ has to be chosen, playing the role of a directional median for~${\rm P}^{\bf Z}$.   Since the uniform~${\rm P}^{\bf U}$   automatically enjoys rotational symmetry with respect to the image ${\bf F}(\thetab)$ of that pole, the nested regions playing the role of ${\rm P}^{\bf U}$-quantile regions should be invariant with respect to rotations with axis ${\bf F}(\thetab)$: hence, the collection of spherical caps with ${\rm P}^{\bf U}$-probability contents $\tau$  centered at~${\bf F}(\thetab)$  naturally qualifies as the family of (nested) quantile regions of orders~$\tau\in~\!\![0,1]$ of  ${\rm P}^{\bf U}$. 

The choice of a directional median or pole $\thetab$, of course, should be ${\rm P}^\Zb$-specific. The literature on directional data provides several concepts of {\it directional medians}, the most convenient of which is the so-called {\it  Fr\'echet mean} of $\Zb\sim{\rm P}^\Zb$, defined as \vspace{-1mm}
 \begin{equation}\label{def.Frechet}
\thetab_{\rm Fr}\coloneqq  {\rm argmin}_{\zb \in  \mathcal{S}^{d-1}} {\rm E}_{{\rm P}^\Zb}[c(\Zb, \zb)]
\end{equation}
(recall that $c(\zb_1, \zb_2)$ denotes the squared Riemannian distance between $\zb_1$ and $ \zb_2$). The Fr\' echet mean always exists but is not necessarily unique. In   case \eqref{def.Frechet} has multiple solutions,  one can randomly select~$\thetab$ as one of them.

%
  Another possible choice 
   can be based on the transport map $\Fb$ (equivalently, on the quantile function $\Qb$) itself. For each $\zb \in {\cal S}^{d-1}$,  consider the hyper-hemisphere~${\cal H}_{\zb}\coloneqq \{ {\bf u} \in {\cal S}^{d-1}: \ub^\top {\bf F}(\zb) \geq 0\}$ centered at ${\bf F}(\zb)$. Since~${\bf F}$  pushes~${\rm P}^\Zb$  forward to ${\rm P}^\Ub\!$, we have that ${\rm P}^\Ub({\cal H}_{\zb})={\rm P}^\Zb({\bf Q}({\cal H}_{\zb}))=1/2$. In general,  however, ${\bf Q}({\cal H}_{\zb})
   $ is not a hyper-hemisphere and~$1/2={\rm P}^\Ub({\cal H}_{\zb}) \neq {\rm P}^\Ub({\bf Q}({\cal H}_{\zb}))$. A transport-based concept of median for ${\rm P}^\Zb$ is thus\vspace{-1mm}
\begin{equation} \label{newtheta}
\thetab_{\rm Tr}\coloneqq  {\rm argmin}_{\zb \in {\cal S}^{d-1}} {\rm P}^\Ub({\bf Q}({\cal H}_{\zb})).
\end{equation}
%
 This  directional median or pole $\thetab_{\rm Tr}$ indeed is such that the area of the image by ${\bf Q}$ of the hyper-hemisphere centered at~${\bf F}(\thetab_{\rm Tr})$ is minimal 
  among  the areas of all 
 images~${\bf Q}({\cal H}_{\zb})$ of  hyper-hemispheres, indicating a concentration of the ${\rm P}^\Zb$ probability mass around~$\thetab_{\rm Tr}$.  Such a pole always exists. Again, it may not be unique and one can randomly select $\thetab$ within the set of points satisfying \eqref{newtheta}.

Assume that a pole $\thetab_{\rm M}$, say, has been selected  for ${\rm P}^\Zb$ ($\thetab_{\rm Fr}$ or  $\thetab_{\rm Tr}$ being possible choices). 
The spherical cap with~${\rm P}^{\bf U}$-probability  $\tau$ centered at~${\bf F}(\thetab_{\rm M})$ is\vspace{-1mm}
\[
{\mathbb C}_\tau^{\Ub} \coloneqq 
\left\{{\bf u}\in {\cal S}^{d-1} : F_{*}\left({\bf u}^\top {\bf F}(\thetab_{\rm M}) \right)\geq 1-\tau
\right\} \quad 0\leq\tau \leq 1,\vspace{-1mm}
\]
with boundary\vspace{-1mm}
\begin{equation}\label{paraldef}
{\mathcal C}_\tau^{\Ub} \coloneqq 
\left\{{\bf u}\in {\cal S}^{d-1} : F_{*}\left({\bf u}^\top {\bf F}(\thetab_{\rm M}) \right)= 1-\tau
\right\} \quad 0\leq\tau \leq 1,\vspace{-1mm}
\end{equation}
{(a $(d-2)$-sphere)} where \vspace{-1mm}
\begin{equation}\label{eq:unif}
F_{*}(u)\coloneqq  {\int_{-1}^u (1-s^2)^{(d-3)/2} \; \dd s}{\Big/}{\int_{-1}^1 (1-s^2)^{(d-3)/2} \;  \dd s}, \quad -1 \leq u \leq 1.
\end{equation}
Indeed, it is easy to see  that  $u\mapsto F_{*}(u)$ is the distribution function of  
${\bf U}^\top{\bf F}(\thetab_{\rm M})$ where~${\bf U}\sim {{\rm P}^\Ub}$. Hence,  $F_*({\bf u}^\top{\bf F}(\thetab_{\rm M}))$ is the ${\rm P}^\Ub$-probability of the ${\bf F}(\thetab_{\rm M})$-centered spherical cap running through $\ub$ and  a measure of $\ub$'s latitude, ranging from~0 (for~$\ub =- {\bf F}(\thetab_{\rm M})$) to 1 (for $\ub = {\bf F}(\thetab_{\rm M})$) and scaled in such a way that the latitude of~${\bf U}\sim {{\rm P}^\Ub}$ is uniform over $[0,1]$. Therefore, for the hypersphere $ {\cal S}^{d-1}$ {equipped with the uniform distribution~$ {{\rm P}^\Ub}$} and the pole ${\bf F}(\thetab_{\rm M})$, the contour ${\mathcal C}_\tau^{\Ub}$ plays the role of a {\it parallel} of order $\tau$ (${\mathcal C}_{1/2}^{\Ub}$, thus, is an {\it equatorial hypersphere}).

Accordingly,   define the quantile contour and quantile region of order $\tau$ of $\Zb \sim {\rm P}^\Zb$  as the images 
\begin{equation}\label{qcontdef}
{\mathcal C}_\tau \coloneqq {\Qb}\left({\mathcal C}_\tau^{\Ub}\right)
=
 \left\{\zb\in \mathcal{S}^{d-1} : F_{*} \left((\Fb (\zb)^\top \Fb(\thetab_{\rm M}))\right) =1-\tau
\right\}
\vspace{-2mm}\end{equation}
and   \vspace{-2mm}
\begin{equation}\label{qregdef}
{\mathbb C}_\tau \coloneqq
 {\Qb}\left({\mathbb C}_\tau^{\Ub}\right)
=
 \left\{\zb\in \mathcal{S}^{d-1} : F_{*} \left((\Fb  (\zb))^\top \Fb(\thetab_{\rm M})\right) \geq 1- \tau
\right\},\vspace{-2mm}\end{equation}
 by $\bf Q$ of ${\mathcal C}_\tau^{\Ub}$ and ${\mathbb C}_\tau^{\Ub}$, respectively. Since $\Qb$ is a measure-preserving transformation pushing ${\rm P}^\Ub$ forward to ${\rm P}^\Zb$, the~${\rm P}^\Zb$ probability content of ${\mathbb C}_\tau$ is $\tau$. 
 The following proposition summarizes the main properties of $\Fb$, $\Qb$,   the quantile contours~$\mathcal C_\tau$ and quantile regions~$\mathbb C_\tau$. They all follow from the definitions and the continuity of $\Fb$ and~$\Qb$; details are left to the reader.
\begin{proposition} \label{characterize}
Let $\Zb\sim{\rm P}^\Zb  \in \mathfrak{P}_d$ have distribution and quantile functions~$\Fb$ and~$\Qb$, respectively. Then, 
\begin{compactenum}
\item[(i)] $\Fb$ entirely characterizes ${\rm P}^\Zb$,  $\Fb(\Zb)\sim\Fb\#{\rm P}^\Zb=~\!{\rm P}^\Ub$, and  $F_{*} \left((\Fb (\Zb))^\top \Fb(\thetab_{\rm M})\right) \sim~\!{\rm U}_{[0,1]}$;
\item[(ii)] $\Qb$ entirely characterizes ${\rm P}^\Zb$ and $\Qb(\Ub)\sim\Qb\#{\rm P}^\Ub=~\!{\rm P}^\Zb$
;
\item[(iii)] the quantile contours $\mathcal C_\tau$, $\tau\in[0,1]$ are continuous;  the quantile regions~${\mathbb C}_\tau$ are   closed, connected, and nested; their intersec\-tion~$\bigcap_{\tau\in[0,1]} \mathbb{C}_\tau$ is the {\rm directional median}~$\thetab_{\rm M}$;
\item[(iv)] the probability content ${\rm P}^\Zb\left(\mathbb C_\tau
\right)$ of ${\mathbb C}_\tau$, $\tau\in[0,1]$,  is $\tau$, irrespective of $\rm P ^\Zb\in  \mathfrak{P}_d$. \vspace{-2mm}
\end{compactenum}
\end{proposition}

%
%
%
%
%
%

\subsection{A natural distribution-specific coordinate system}\label{coordsec} The quantile function $\Qb$ of~${\rm P}^\Zb$ actually creates over $\mathcal{S}^{d-1}$ a coordinate system with latitudes and (hyper)longitudes adapted to the distribution of $\Zb$. 

\begin{figure}[t!]
\centering\vspace{-3mm}
\begin{subfigure}[b]{0.4\textwidth} 
\begin{tikzpicture}
\def\R{2.5} 
\def\angEl{35} 
\def\angAz{-105} 
\def\angPhi{-40} 
\def\angBeta{19} 

\pgfmathsetmacro\H{\R*cos(\angEl)} 
\tikzset{xyplane/.estyle={cm={cos(\angAz),sin(\angAz)*sin(\angEl),-sin(\angAz), cos(\angAz)*sin(\angEl),(0,-\H)}}}
\fill[ball color=lightgreen] (0,0) circle (\R); 
\draw (0,0) circle (\R);

\coordinate (O) at (0,0);
\LatitudePlane[equator]{\angEl}{0};

\draw[->] (0,0) -- (0,\R) node[above] {\scriptsize  ${\bf F}(\thetab_{\rm M})$};
\draw[->] (0,0) -- (.6*\R,.6*\R) node[above] {\scriptsize ${\bf F}({\bf z})$};
\draw[->] (0,0) -- (0,.77*\R) node[below=.45cm, left=-.01cm, , inner sep=0.1ex] {\scriptsize {$(({\bf F} (\thetab_{\rm M}))^\top {\bf F}(\zb)){\bf F}(\thetab_{\rm M})$}};
\draw[dashed] (.6*\R,.6*\R) -- (0,.77*\R);
\draw[->] (0,0) -- (.6*\R, -.2*\R) node[below left] {};
\draw[dashed] (.6*\R,.6*\R) -- (.6*\R, -.2*\R);
\draw[darkorange,->] (0,0) -- (1.45*.6*\R, -1.45*.2*\R) node[above=.8cm, right=-1.3cm, inner sep=0.1ex] {\small  ${\bf S}_{{\bf F}(\thetab_{\rm M})} {\bf F}(\zb)$};

\draw[equator,->] (0,0) circle (\R);
\end{tikzpicture}
\end{subfigure}
\hspace{10mm}
\begin{subfigure}[b]{0.4\textwidth}
\centering
\begin{tikzpicture} 
\label{Figquant1}

\def\R{2.5} 
\def\angEl{35} 
\def\angAz{-105} 
\def\angPhi{-40} 
\def\angBeta{19} 

\pgfmathsetmacro\H{\R*cos(\angEl)} 
\tikzset{xyplane/.estyle={cm={cos(\angAz),sin(\angAz)*sin(\angEl),-sin(\angAz),
                              cos(\angAz)*sin(\angEl),(0,-\H)}}}

\fill[ball color=lightgreen] (0,0) circle (\R); 
\draw (0,0) circle (\R);

\coordinate (O) at (0,0);
\LatitudePlane[equator]{\angEl}{0};

\draw[->] (0,0) -- (0,\R) node[above] {{\scriptsize  ${{\pmb \theta}_{{\rm M}}}$}};
\draw[->] (0,0) -- (.6*\R,.6*\R) node[above] {{\scriptsize  ${\bf z}$}};
\draw[->] (0,0) -- (.85*\R,.05*\R) node[above=.4cm, left=-.1cm] {{\scriptsize  ${\bf F}({\bf z})$}};
\draw[->] (0,0) -- (0,.8*\R) node[below=.0cm, left=.1cm] {\scriptsize  {$(\zb^\top \thetab_{\rm M} ) \thetab_{\rm M}$}};
\draw[->] (0,0) -- (0,.35*\R) node[below=.0cm, left=.1cm] {\scriptsize  {$F_f^*(\zb^\top \thetab_{\rm M} ) \thetab_{\rm M}$}};
\draw[dashed] (.6*\R,.6*\R) -- (0,.8*\R);
\draw[dashed] (.85*\R,.05*\R) -- (0,.35*\R);
\draw[dashed] (.85*\R,.05*\R) -- (.85*\R, -.3*\R);
\draw[->] (0,0) -- (.6*\R, -.2*\R) node[below left] {};
\draw[dashed] (.6*\R,.6*\R) -- (.6*\R, -.2*\R);
\draw[darkorange,->] (0,0) -- (1.45*.6*\R, -1.45*.2*\R) node[above=.5cm, left=.05cm] {$ {\bf S}_{{\pmb \theta}_{\rm M}}(\zb)$};

\draw[equator,->] (0,0) circle (\R);
\end{tikzpicture}
\end{subfigure} 

\caption{ Left: the tangent-normal decomposition of   ${\bf F}$ 
  with respect to ${\bf F}(\thetab_{\rm M})$. Right: in the rotationally symmetric case, the optimal transport $\Fb$  reduces to an optimal univariate transport acting on the projection~${\bf z}^\top \thetab_{\text{\rm M}}$ of~$\zb$ along the rotation axis $\thetab_{\text{\rm M}}$.}\label{tangentfig}
\end{figure}
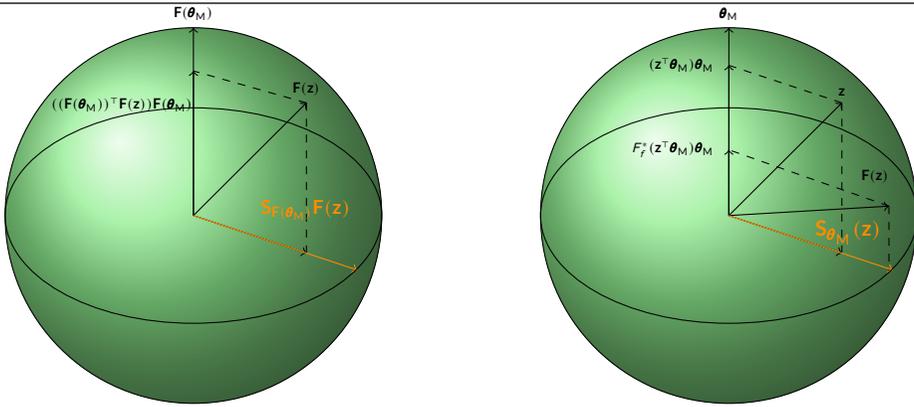

The usual latitude-(hyper)longitude coordinate system is based on the classical   {\it  tangent-normal} decomposition,  with respect to a pole ${\bf F}(\thetab_{\rm M})$, say, of a point   $\ub\in \mathcal{S}^{d-1}$  into the sum of two mutually orthogonal terms
  \begin{equation}\label{tangenteq}
{\bf u}  = ({\bf u}^\top  {\bf F}(\thetab_{\rm M})) {\bf F}(\thetab_{\rm M})+ ({\bf I}_p- {\bf F}(\thetab_{\rm M})({\bf F} (\thetab_{\rm M}))^\top){\bf u} 
=
 ({\bf u}^\top  {\bf F}(\thetab_{\rm M})){\bf F}(\thetab_{\rm M})+ \sqrt{1- ({\bf u}^\top  {\bf F}(\thetab_{\rm M}))^2} {\bf S}_{ {\bf F}(\thetab_{\rm M})\!}({\bf u} )
\end{equation}
(see the left panel of Figure~\ref{tangentfig} for~$d=3$) where ${\bf u}^\top  {\bf F}(\thetab_{\rm M})$, being constant over the parallel $\mathcal C ^\Ub_{\tau }$ ($\tau =1-F_*({\bf u}^\top  {\bf F}(\thetab_{\rm M}))$) defined in Section~\ref{31sec}, is a latitude  while the unit vector (a {\it directional sign}) 
\begin{equation} \label{sign}
\mathbf{S}_{{\bf F}(\thetab_{\rm M})}({\ub})\coloneqq  {(\ub-(\ub^\top {\bf F}(\thetab_{\rm M})){\bf F}(\thetab_{\rm M}))
   }/{\|\ub-(\ub^\top {\bf F}(\thetab_{\rm M})) {\bf F}(\thetab_{\rm M}) \|}
 \end{equation}   
(with the convention  ${\boldsymbol 0}/0={\boldsymbol 0}$ for $\ub=\pm  {\bf F}(\thetab_{\rm M})$)  with values  on the equatorial hypersphere ${\mathcal C}_{1/2}^{\Ub}$---a hyperlongitude, thus---characterizes hypermeridians 
$
{\mathcal M}_{\bf s}^\Ub\coloneqq \left\{ \ub \in \mathcal{S}^{d-1}: \mathbf{S}_{{\bf F}(\thetab_{\rm M})}({\ub})={\bf s}
\right\}$, ${\bf s}\in{\mathcal C}_{1/2}^{\Ub}
$. 

This parallel/hypermeridian system is well adapted to ${\mathcal S}^{d-1}$ when equipped with the  surface area measure 
or the uniform measure ${\rm P}^\Ub$:  parallels,   for instance, then coincide with ${\rm P}^\Ub$'s quantile contours $\mathcal C ^\Ub_{\tau }$, while the hyper\-longi\-tudes~$\mathbf{S}_{{\bf F}(\thetab_{\rm M})}({\Ub})$, for $\Ub\sim{\rm P}^\Ub$, are uniform over the equatorial hypersphere. Its image by $\Qb$ is more natural under the probability measure ${\rm P}^\Zb$: the same properties, indeed, then hold for the ``curvilinear''    parallels $\Qb(\mathcal C ^\Ub_{\tau })={\mathcal C}_\tau$ and the  ``curvilinear''  hypermeridians \vspace{-2mm}
\begin{equation}\label{meriddef}
{\mathcal M}_{\bf s}\coloneqq \Qb({\mathcal M}_{\bf s}^\Ub)\coloneqq \left\{\zb \in \mathcal{S}^{d-1}:  \mathbf{S}_{{\bf F}(\thetab_{\rm M})}(\Fb(\zb)
)={\bf s}
\right\},\quad {\bf s}\in{\mathcal C}_{1/2}^{\Ub}.
 \vspace{-4mm}\end{equation}
%
\subsection{Rotational symmetry}
Optimal transports seldom admit closed-form expressions. In the particular case of rotational symmetry, however, explicit expressions for~${\bf F}$ and $\Qb$ are possible.

Consider  the rotationally symmetric variable $\Zb\sim {\rm P}^{\Zb}$ with uniquely defined axis $\pm\,\thetab$; all sensible directional\linebreak  medians~$\thetab _{\rm M}$ (including $\thetab_{\rm Fr}$ and~$\thetab_{\rm Tr}$) then are  lying  along that axis, and  ${\rm P}^{\Zb}$ has density  \eqref{Rotsymfirst} with  $\thetab=\thetab_{\rm M}$: \linebreak write ${\rm P}^{\Zb}\eqqcolon {\rm P}_{\thetab_{\rm M}, f}$. Denote~by  
$$F_{f}(u)\coloneqq  {\int_{-1}^u f(s) (1-s^2)^{(p-3)/2} \; \dd s}{\big /}{\int_{-1}^1 f(s) (1-s^2)^{(p-3)/2} \; \dd s}, \quad -1 \leq u \leq 1,$$
 the distribution function  of  $\Zb^\top \thetab_{\rm M}$, by   $Q_{f}\coloneqq F_{f}^{-1}$ its  quantile function. Under rotational symmetry, we have that 
\begin{equation} \label{sign}
\mathbf{S}_{\thetab_{\rm M}}(\mathbf{\Zb}) =  {(\Zb-(\Zb^\top \thetab_{\rm M})\thetab_{\rm M})
   }/{\|\Zb-(\Zb^\top \thetab_{\rm M}) \thetab_{\rm M} \|}
 \end{equation}   
 (with the convention that ${\boldsymbol 0}/0={\boldsymbol 0}$)   is uniformly distributed over the {\it equatorial hypersphere} 
  \begin{equation}\label{eqdef}
 \mathcal{S}^
{d-2}_{\thetab^\perp_{\rm M}}\coloneqq \left\{ \ub \in \R^d:   \ub^\top \thetab_{\rm M}=0, \ub^\top \ub=1\right\}
\end{equation}
 (see, e.g., \cite{Paindaveine2017}). 
Actually, $\mathbf{S}_{\thetab_{\rm M}}(\mathbf{\zb})$ is the point in the equatorial sphere with the same {\it hyperlongitudes} as~$\zb$ and can be interpreted as the {\it directional sign} of $\zb$. Since $\Fb(\Zb)$, by definition, is uniformly distributed over~$ \mathcal{S}^{d-1}$, 
 $(\Fb (\Zb))^\top \thetab_{\text{\rm M}}$ has distribution function $F_*$ as defined in  \eqref{eq:unif}
and  quantile function $Q_{*}=F_*^{-1}$. We then have an explicit  form for the directional distribution function ${\bf F}$ of~$\Zb\sim {\rm P}_{\thetab_{\rm M}, f}$ (see Section~\ref{Sec.Proofs} of the  Supplementary Material  for a proof).
\begin{proposition}\label{PropSym}
  Let $\Zb$ have rotationally symmetric distribution ${\rm P}_{\thetab_{\rm M}, f}$. Then,   letting $F_{{f}}^*(u)\coloneqq Q_* (F_f(u))$, 
  \begin{equation}
    \label{eq:2}
    \mathbf{F}(\mathbf{z}) = F_{f}^* ({\bf z}^\top \thetab_{\rm M}) \thetab_{\rm M} +
    \sqrt{1 - (F^*_{{f}}({\bf z}^\top \thetab_{\rm M}))^2} \,\mathbf{S}_{\thetab_{\rm M}}(\mathbf{z}).
  \end{equation} 
\end{proposition}

Proposition \ref{PropSym} tells us that, in the rotationally symmetric case, the optimal transport $\Fb$ essentially reduces to the optimal univariate transport, the function $F_{{f}}^*$ acting on the projection~${\bf z}^\top \thetab_{\text{\rm M}}$ of $\zb$ onto the  axis $\thetab_{\text{\rm M}}$. This is illustrated in Figure \ref{tangentfig} (right panel). A consequence of Proposition \ref{PropSym} is that, in the rotationally symmetric case, our quantile regions coincide with the canonical quantile regions of \cite{Ley2014}.  

\section{Empirical distribution and quantile functions, ranks,\\  and signs} \label{sec:emp}\vspace{-3mm}\setcounter{equation}{0}

So far, we have been dealing with  population concepts of distribution and quantile functions. If a statistical analysis is to be performed, we need empirical versions of the same; these  involve  concepts of ranks and signs, which are completely new in the context of directional observations.\vspace{-3mm}

\subsection{Empirical directional distribution function}\label{sec:emp1}

Denoting   by  $\Zb_1^{(n)}, \ldots, \Zb_n^{(n)}$ a sample  of $n$  i.i.d.\ observations with 
 distribution ${\rm P}^\Zb$, consider  a ``regular'' $n$-point
\linebreak  grid~$\mathfrak{G}\n=~\!\{{\scriptstyle{\mathfrak{G}}}\n_1,\ldots,{\scriptstyle{\mathfrak{G}}}\n_n\}$ over the unit hypersphere ${\mathcal S}^{d-1}$. In this section, the only property required from 
 $\mathfrak{G}\n$  is that the sequence of uniform distributions over the~$n$-tuple ${\scriptstyle{\mathfrak{G}}}\n_1,\ldots,{\scriptstyle{\mathfrak{G}}}\n_n$ of gridpoints converges weakly to the  uniform distribution ${\rm P}^{\bf U}$ over ${\mathcal S}^{d-1}$ as~$n\to\infty$. 

 Since $\Fb$ minimizes the transportation cost \eqref{cost}, a plugin estimator 
  is obtained as the solution 
 of an optimal coupling problem between the observations and the grid~$\mathfrak{G}\n$. More precisely, let $\mathcal{T}\n$ denote the set of all permutations of the inte\-gers~$\{1,\ldots,n\}$: each permutation $T\in~\!\mathcal{T}\n$ defines a   bi\-jection~$\Zb\n_i\mapsto {\scriptstyle{\mathfrak{G}}}\n_{T(i)}$  between~$\{\Zb_1^{(n)}, \ldots, \Zb_n^{(n)}\}$ and~$\mathfrak{G}\n$.  The {\it empirical directional distribution function} $\Fb^{(n)}$ is then defined as the mapping 
\begin{equation}\label{Fndef}
\Fb^{(n)}: \Zb^{(n)}\coloneqq  (\Zb_1^{(n)}, \ldots, \Zb_n^{(n)}) \mapsto (\Fb^{(n)}(\Zb_1^{(n)}), \ldots, \Fb^{(n)}(\Zb_n^{(n)}))
\end{equation}
satisfying (with $c(\zb_1,\zb_2)$ the squared Riemannian distance)
\begin{equation}\label{eq.optCouple}
\sum_{i=1}^n c(\Zb_i^{(n)} , \Fb^{(n)}(\Zb_i^{(n)})) = \underset{T\in \mathcal{T}\n }{\min}  \sum_{i=1}^n  c(\Zb_i^{(n)} , {\scriptstyle{\mathfrak{G}}}\n_{T(i)}). 
\end{equation}
The terminology {\it empirical directional distribution function} is justified by the following Glivenko-Cantelli result (see Section~\ref{Sec.Proofs} of the  Supplementary Material  for  a  proof).
\begin{proposition}\label{Gliv}
{\rm (Glivenko-Cantelli).} Let $\Zb_1^{(n)}, \ldots, \Zb_n^{(n)}$ be i.i.d.\ with distribution~${\rm P}^{\Zb}$ over~$ \mathcal{S}^{d-1}$. Then, provided that the sequence of uniform discrete distributions over the $n$ gridpoints $\{{\scriptstyle{\mathfrak{G}}}\n_1,\ldots,{\scriptstyle{\mathfrak{G}}}\n_n\}$ of ${\mathfrak{G}}\n$ converges weakly,  as $n\to\infty$,  to the uniform distribution~${\rm P}^{\Ub}$ over~$ \mathcal{S}^{d-1}$, 
\begin{equation}\label{GC}
{\rm max}_{1 \leq i \leq n} \big\Vert \Fb^{(n)}(\Zb_i^{(n)})- \Fb(\Zb_i^{(n)})\big\Vert \longrightarrow 0\qquad\text{almost surely as $\ny$. }
\end{equation}
\end{proposition}
If a consistent estimation of $\Fb$ is the   objective, the empirical distribution $\Fb^{(n)}$ defined in \eqref{Fndef}, in view of \eqref{GC},   offers a perfect solution. If   empirical counterparts of the quantile contours and regions $\mathcal C_\tau$ and $\mathbb C_\tau$, hence of the parallels ${\mathcal C}_\tau$ and hyperme\-ridians~${\mathcal M}_{\bf s}$ (see \eqref{paraldef} and \eqref{meriddef}) of Section~\ref{coordsec} are to be constructed, or if ranks and signs generating a maximal ancillary sigma-field  are to be defined, however, more structure is required from the grid~$\mathfrak{G}\n$.

\subsection{Directional ranks, signs, and empirical quantiles}\label{ranksandsigns}

In this subsection,  we show how imposing some additional finite-$n$ structure on the grid~$\mathfrak{G}\n$  yields  empirical versions   of $\Fb$ with natural data-driven concepts of 
\begin{compactenum}
\item[(i)] empirical quantile contours and regions consistently estimating the actual ones,
\item[(ii)] a data-driven coordinate system of empirical parallels and hypermeridians adapting to the underlying distribution ${\rm P}^\Zb$ of the observations, and 
\item[(iii)] distribution-free 
  signs and ranks,   paving the way to a theory of rank-based inference for directional data with unspecified density. 
\end{compactenum}

All these concepts involve   a pole $\thetab_{\rm M}$ (the population Fr\' echet mean, for instance) or    a consistent estimator $\thetatM$ 
 thereof  (e.g., the empirical Fr\' echet mean). That 
  $\thetatM$ in turn will define the pole $\thetat\!$ of a  structured grid~$\mathfrak{G}\n({\thetat}\!)$ 
  to be used in the estimation of $\Fb$. 
%
 Accordingly, we are proceeding  in two steps: a first step yielding  an empirical pole $\thetat\!$, then  a second step estimating $\Fb$ on the basis of    the grid~$\mathfrak{G}^{(n)}(\thetat)$ exploiting the role of $\thetat$ as an empirical pole.\smallskip

\noindent{\bf Step 1.} Construct an empirical version~$\thetat$ of the population quantity~${\bf F}(\thetab_{\rm M})$. Let~$\thetatM$ be a consistent estimator of~$\thetab_{\rm M}$ and consider an $n$-point grid $\mathfrak{G}\n_0$ satisfying the assumptions of Proposition~\ref{Gliv};  denoting by~${\bf F}\n_0(\Zb_i^{(n)})$ the resulting estimator of ${\bf F}(\Zb_i^{(n)})$, choose 
$\thetat\coloneqq {\bf F}\n_0(\Zb_i^{(n)})$ where $\Zb_i^{(n)}$ is the sample point closest to $\thetatM$ in the~$d(\cdot, \cdot)$  metric:  Proposition~\ref{Gliv} and the continuity of $\Fb$ imply that $\thetat$ a.s.\ converges to~${\bf F}(\thetab_{\rm M})$.

\noindent{\bf Step 2.} (2a) Construct a further regular grid $\mathfrak{G}\n({\thetat})$ over 
${\cal S}^{d-1}$. Factorizing $n$ into $n = n_R n_S+n_0$ where $n_R, n_S, n_0 \in \mathbb{N}$ and $0 \leq n_0 < \min\{n_R, n_S\}$, 
define that grid $\mathfrak{G}\n({\thetat})$   as the product of two independent grids:
  \begin{compactenum}
  \item[(i)] a reference grid $\mathfrak{S}^{(n_S)}\coloneqq \{{\sbb_1,\ldots,\sbb_{n_S}}\}$ over ${\mathcal S}^{d-2}$; again, this grid should be as uniform as possible but the only requirement is the weak convergence, as $n_S\to~\!\infty$, of  the uniform distribution over $\mathfrak{S}^{(n_S)}$ to the uniform over~${\mathcal S}^{d-2}$. Note that for~$d=~\!3$,  a fully regular grid $\sbb_1,\ldots,\sbb_{n_S}$ is obtained by dividing the unit circle ${\cal S}^1$ into $n_S$ equal parts (see Figure~\ref{coord});
  \item[(ii)] a grid of $n_R$ points over the unit interval, of the form ${i}/{(n_R+1)}$, $i=1, \ldots, n_R$.\vspace{2mm}
  \end{compactenum}

  (2b)   Import the grid $\mathfrak{S}^{(n_S)}$ of step (2a) to the equatorial space defined by the pole~$\thetat$ computed in step 1.  More precisely, construct the grid $\mathfrak{G}\n({\thetat})$ that consists in~$n_0$ copies of  $\thetat$ (if $n_0 \neq 0$) and the $n_R n_S$ points~${\scriptstyle{\mathfrak{G}}}\n_{ij}({\thetat})$  such that  \vspace{-3mm}
 $$
1 - F_{*}(({\scriptstyle{\mathfrak{G}}}\n_{ij}({\thetat}))^\top \thetat)
 = \frac{i}{n_R+1}
\quad 
\text{ and } \quad 
\mathbf{S}_{\,\thetatsub}
({\scriptstyle{\mathfrak{G}}}\n_{ij}({\thetatsub})) = \Gamb_{\!\thetatsub} {{\bf s}}_j,\quad i=1,\ldots,n_R,\, j=1,\ldots,n_S
\vspace{-1mm}$$ 
where $\Gamb_\thetab$ denotes a $d \times (d-1)$ semi-orthogonal matrix such that $\Gamb_\thetab \Gamb_\thetab^\top= {\bf I}_d - \thetab \thetab^\top$ and~$\Gamb_\thetab^\top\Gamb_\thetab= {\bf I}_{d-1}$; the columns of~$(\thetat,\Gamb_{\!\thetatsub})$, thus, constitute an orthonormal coordinate system of ${\mathbb R}^d$,  the columns of $\Gamb_{\!\thetatsub}$ an arbitrary orthonormal coordinate system of the equatorial hyperplane determined by $\thetat\!$. Note that, due to its dependence on $\thetat\!$, this grid~$\mathfrak{G}\n({\thetat})$~is~random. \vspace{1.5mm}
 \vspace{0mm}

\noindent (2c) 
 Denote by~$\Fb^{(n)}(\Zb_1^{(n)}), \ldots, \Fb^{(n)}(\Zb_n^{(n)})$ the solutions of the optimal coupling problem~\eqref{eq.optCouple} based on this second grid $\mathfrak{G}\n({\thetat})$.
  \smallskip\smallskip

  The  (directional) empirical sign ${\bf S}_i\n$ and the (directional) rank $R_i\n$ of   $\Zb_i\n$ then are  naturally defined as\vspace{-1.5mm} 
\begin{equation}\label{signFemp}
 {\bf S}_i\n \coloneqq {\bf S}_{\,\thetatsub}
 ({\bf F}\n(\Zb_i\n))\vspace{-4.5mm}
\end{equation}
and  \vspace{-3.5mm}
\begin{equation}\label{rankFemp}
R\n_i  = R\n(\Zb_i\n)\coloneqq  (n_R+1)\left [1 - F_{*}(({\bf F}\n(\Zb_i\n))^\top\thetat
 ))\right], \quad i=1,\ldots,n\vspace{-2.5mm}
\end{equation}
(with values in $\{1,\ldots,n_R\}$), respectively, provided that ${\bf F}\n(\Zb_i\n)\neq\thetat$; for $\Zb_i\n$ such that ${\bf F}\n(\Zb_i\n)=\thetat$, let~$ {\bf S}_i\n\coloneqq {\boldsymbol 0} $ 
 (the pole has no specific longitude) and  $R\n_i  \coloneqq 0 \vspace{0mm}$.

The empirical versions of the quantile contours and regions 
\eqref{qcontdef} and~\eqref{qregdef} are the collections of observations 
\begin{equation}\label{def.EmpQuanCont}
\mathcal{C}\n_{j/(n_R+1)} \! \coloneqq \! \left\{\Zb_i^{(n)}\! \!: R_i\n\! = j\right\}\quad\!\!\!\text{and}\!\quad\!\! 
\mathbb{C}\n_{j/(n_R+1)} \!\coloneqq  \!\left\{\Zb_i^{(n)}\!\! : R_i\n \leq j\right\},\quad\! j=1,\ldots, n_R,
\end{equation}
respectively. An empirical quantile contour of order $j/(n_R+1)$ (the empirical parallel of order~$j/(n_R+1)$)  
thus consists of the $n_S$ sample points with given rank~$j$. Similarly, the empirical hypermeridian  with longitude ${\bf s}$ consists of the~$n_R$ observations with given  sign (hyperlongitude) ${\bf s}$:
\begin{equation}\label{def.EmpMerid}{\mathcal M}\n_{{\bf s}} 
 \!\coloneqq  \!\left\{\Zb_i^{(n)}\! : \ \Fb\n(\Zb_i^{(n)})\neq \thetat\text{ and }  {\bf S}_i\n = {\bf s}\right\},
 \quad {{\bf s}}\in\{ \Gamb_{\thetatsub} {{\bf s}}_1, \ldots,  \Gamb_{\thetatsub} {{\bf s}}_{n_S}\}.
\end{equation}

Because of the data-driven choice of $\thetat$, the grid $\mathfrak{G}\n({\thetat})$ is random. 
If  $n_R\to\infty$ and $n_S\to\infty$,  the assumptions in  Proposition~\ref{Gliv} are satisfied a.s., and the Glivenko-Cantelli result \eqref{GC} holds a.s.~as well. Similar consistency properties then follow for quantile contours and regions, parallels, and hypermeridians. In Supplementary Material A, we establish and discuss several properties---distribution-freeness, ancillarity, and equivariance---of directional signs and ranks.

  \begin{center}
\begin{figure}[ht!]
\centering
\begin{subfigure}[b]{0.4\textwidth}\hspace{-20mm}
\begin{tikzpicture}[>=latex]
\def\R{3.5}
 \draw (0,0) circle (\R);
\foreach \ang in {0,...,31} {
  \draw [blue] (0,0) -- (\ang * 180 / 16:4);
}
\node[below=0pt] at (\R+1.5, 0.3) {${\bf s}_1=(1,0)^\top$};
\node[below=0pt] at (\R+.9, 1) {\ ${\bf s}_2$};
\node[below=0pt] at (\R+1, -.6) { ${\bf s}_{n_S}$};
\draw[dashed] (\R+.8,1.4) arc (\R+.8: 35: 1.4);
\draw[dashed] (\R+.4,-2) arc (-40: 5: 1);\hspace{-10mm}

\draw (\R+1,0) node {$\bullet$};
\draw (\R+0.95 ,0.67) node {$\bullet$};
\draw (\R + 0.95,-0.677) node {$\bullet$};
\draw (\R+3.9, 2.2)--(\R+4.1,2.2) node[right=2mm] {\scriptsize $1/(n_R+1)$};
\draw (\R+3.9, 1.9)--(\R+4.1, 1.9) node[right=2mm] {\scriptsize $2/(n_R+1)$};
\draw (\R+3.9, 1.6)--(\R+4.1, 1.6);
\draw (\R+3.9, 1)--(\R+4.1, 1);
\draw (\R+3.9, .7)--(\R+4.1, .7);
\draw (\R+3.9, 1.3)--(\R+4.1, 1.3) node[above= 2mm, right=2mm] {\scriptsize $\vdots$};
\draw (\R+3.9, -2.2)--(\R+4.1, -2.2) node[right=2mm] {\scriptsize $n_R/(n_R+1)$};
\draw[|-|,thick, orange] (\R+4, 2.5) node[above=1mm]{0} --(\R+4, -2.5)node[below=1mm]{1};

\end{tikzpicture} 
\end{subfigure}
\\ 
\begin{subfigure}[b]{0.4\textwidth}\hspace{-41mm}
\begin{tikzpicture} 
\label{Figquant3}

\def\R{3.5} 
\def\angEl{35} 
\def\angAz{-105} 
\def\angPhi{-40} 
\def\angBeta{19} 

\pgfmathsetmacro\H{\R*cos(\angEl)} 
\tikzset{xyplane/.estyle={cm={cos(\angAz),sin(\angAz)*sin(\angEl),-sin(\angAz),
                              cos(\angAz)*sin(\angEl),(0,-\H)}}}

\fill[ball color=lightgreen] (0,0) circle (\R); 

\coordinate (O) at (0,0);
\LatitudePlane[equator]{\angEl}{0};

\coordinate (N) at (0,\R);

\foreach \t in {20, 40, 60, 80, 100 , 120, 140, 160, 180} { \DrawLongitudeCircle[\R]{\t}}
\foreach \t in {-40, -30, -20, -10, 10,20, 30, 40, 50, 60, 70} { \DrawLatitudeCircle[\R]{\t}}


\node[above=3pt] at (N) {${\bf F}\n(\widehat{\thetab}_{\rm M})$};


\draw[->] (-1.71, -1.78) -- (-\R, 3.2) node[left=5pt] {\scriptsize grid of $n_S$ points  determining};
\node[left=2pt, below=0pt] at (-\R-1.5, 3.1) {\scriptsize  the meridians (the longitudes)};
\draw[->] (.44*\R, .56*\R ) -- (.9*\R, .8*\R) node[above right] {\scriptsize grid of $n_R$ points  determining};
\node[below=8pt] at (\R+1.0, \R-.3) {\scriptsize   the parallels (the latitudes)};

\draw[equator,->] (0,0) circle (\R);
\end{tikzpicture}
\end{subfigure}
\caption{\small The grid used to define signs and ranks on ${\cal S}^{d-1}$ for $d=3$. The final grid   is obtained as the product of a    reference grid~$\mathfrak{S}^{(n_S)}\coloneqq \{{\sbb_1,\ldots,\sbb_{n_S}}\}$ over ${\mathcal S}^{d-2}$ (here,~$n_S$ equispaced points on the circle ${\mathcal S}^{1}$) and $n_R$ equispaced  points on the~unit interval. \vspace{-4mm}
\label{coord}}
\end{figure}
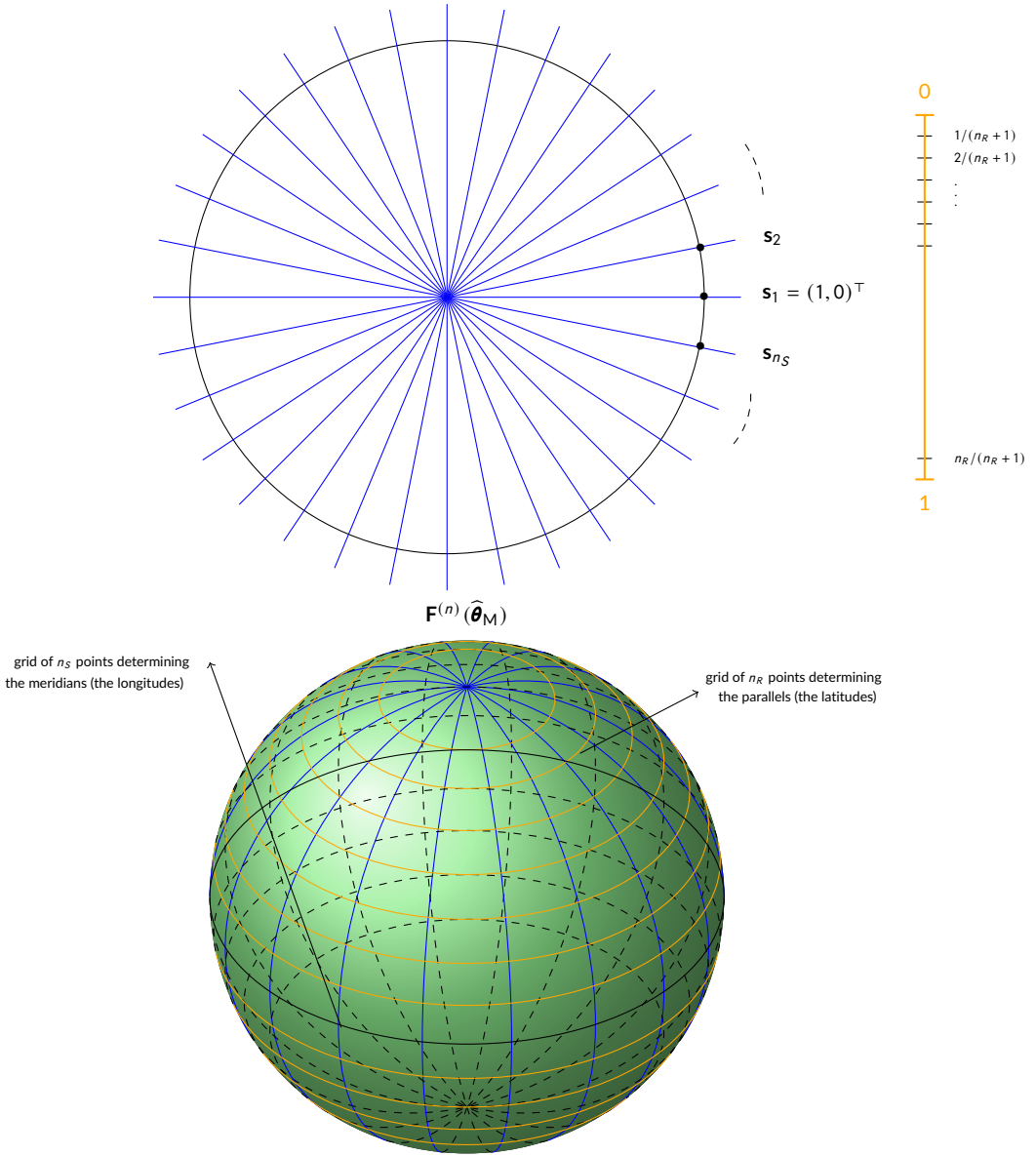
\end{center}

\vspace{-8mm}

\subsection{Numerical illustration: simulations}\label{numerillsec}
 To illustrate the concept of quantile contours in~\eqref{def.EmpQuanCont}, we generated sequences of $n = 2001$ i.i.d.\ unit random vectors from three types of distributions on $ \mathcal{S}^2$: 
\begin{compactenum}
\item[(i)] the von Mises-Fisher (vMF) distribution (see \eqref{vMFdens}) with concentration $\kappa = 10$ and location $\thetab = (0, 0, 1)^\top$. Below,~$\mathcal{M}_{d}(\thetab, \kappa)$ denotes the von Mises distribution on ${\cal S}^{d-1}$ with location $\thetab$ and concentration $\kappa$;
\item[(ii)] the {tangent vMF distribution} as defined  in \cite{Garcia-Portugues2020}. The  {tangent vMF distribution} with location~$\thetab$, angular function~$G$, skewness direction ${\mbf \mu}$, and skewness intensity $\kappa$  is the distribution of \vspace{-2mm}
\begin{equation} \label{tangentvMF}
  \Zb \coloneqq V\thetab + \sqrt{1 - V^2} {\boldsymbol \Gamma}_{\thetab} {\bf U},\vspace{-2mm} \end{equation}
where  $\Gamb_\thetab$ is the $d \times (d-1)$ semi-orthogonal matrix described in Part (2b) of Step~2 of the construction (Section~\ref{ranksandsigns}) of the grid,   $V$ 
 (an absolutely continuous scalar random variable with values in [-1,1]) and ${\bf U} \sim \mathcal{M}_{2}({\mbf \mu}, \kappa)$ are mutually independent; in the simulation, we set  
 $\thetab = (0, 0, 1)^\top$,  ${\mbf \mu} = (0.7, \sqrt{0.51})^\top$, $\kappa = 10$, and $V= 2\widetilde{V} -1$ with~$\widetilde{V}\sim\text{\rm Beta}(2, 8)$; 
\item[(iii)] a mixture of two vMF distributions---the distribution of
 $I[U \leq 0.3] \Zb_1 + I[U > 0.3] \Zb_2,$ where $U \sim {\rm U}[0, 1]$,~$I[\: ]$ denotes the indicator function,  and $\Zb_1 \sim \mathcal{M}_3(\thetab_1, \kappa_1)$ and  $\Zb_2 \sim \mathcal{M}_3(\thetab_2, \kappa_2)$   are mutually independent. In the simulation, we set~$\kappa_1 = 20$, $\thetab_1 = (0, -0.5, \sqrt{0.75})^\top$, $\kappa_2 = 20$, and $\thetab_2 = (0, 0.5, \sqrt{0.75})^\top$.
\end{compactenum}

\begin{figure}[h!]
     \centering
     \begin{subfigure}[b]{0.29\textwidth}
         \hspace{3mm}\vspace{-4.8mm}
         \includegraphics[trim=8 15 25 8,clip, width=1\textwidth]{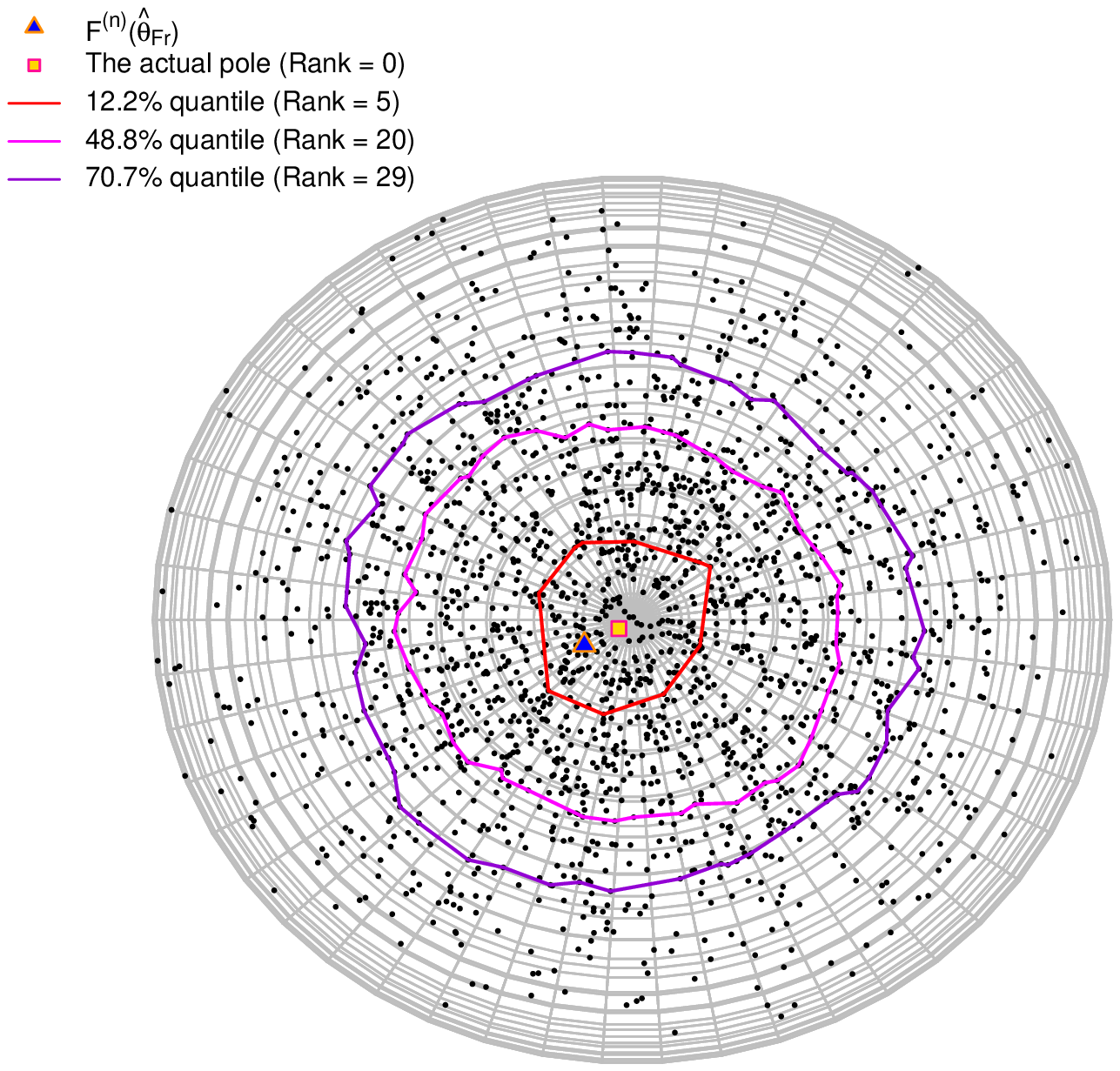} 
     \end{subfigure}
     \hfill    
 \begin{subfigure}[b]{0.28\textwidth}
         \hspace{-13mm}
       \vspace{-13.8mm}
         \includegraphics[trim=20 15 25 12,clip, width=1.4\textwidth]{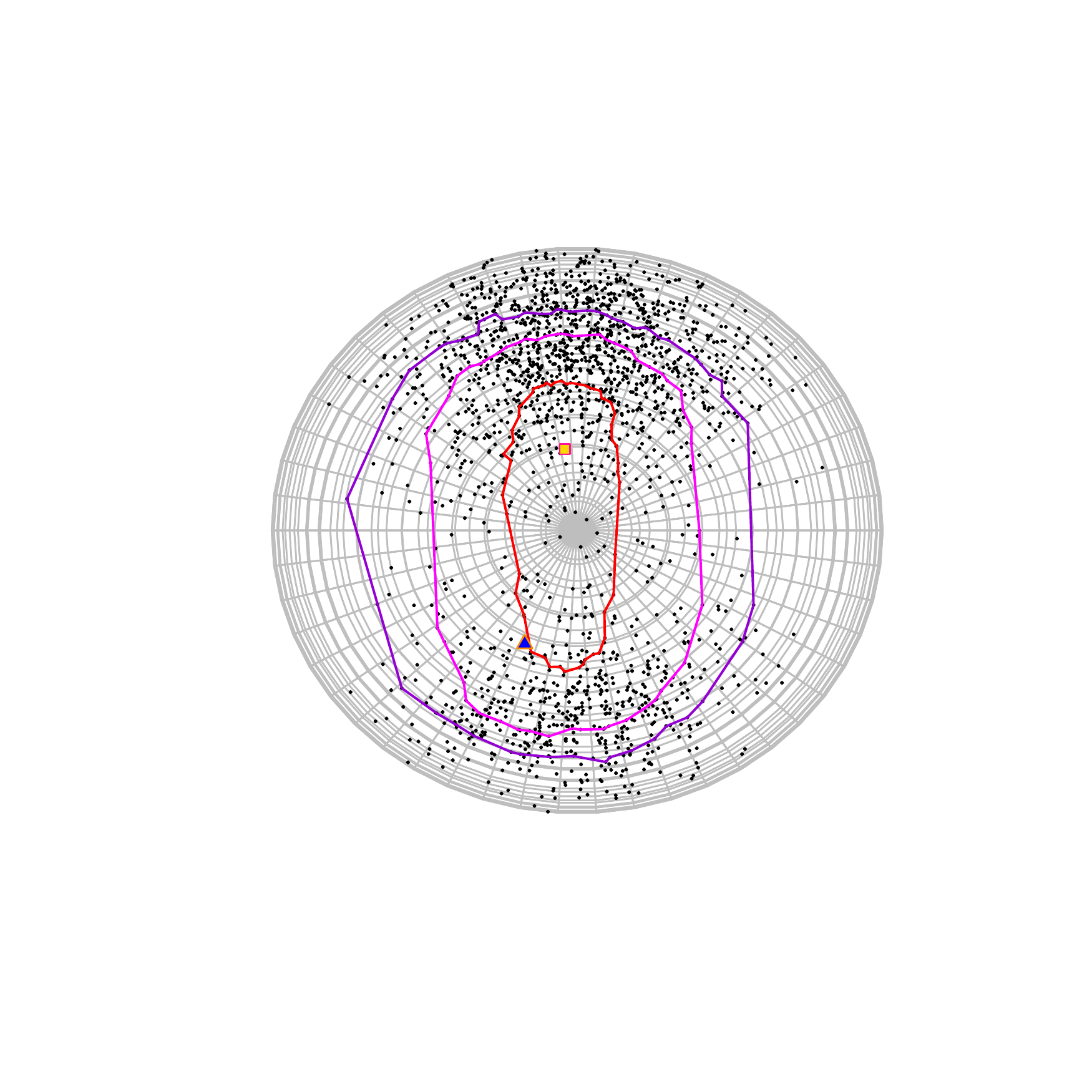} 
     \end{subfigure}
     \hfill  
  \begin{subfigure}[b]{0.28\textwidth}
         \hspace{-20mm}
       \vspace{-13mm}
         \includegraphics[trim=20 0 25 12,clip, width=1.53\textwidth]{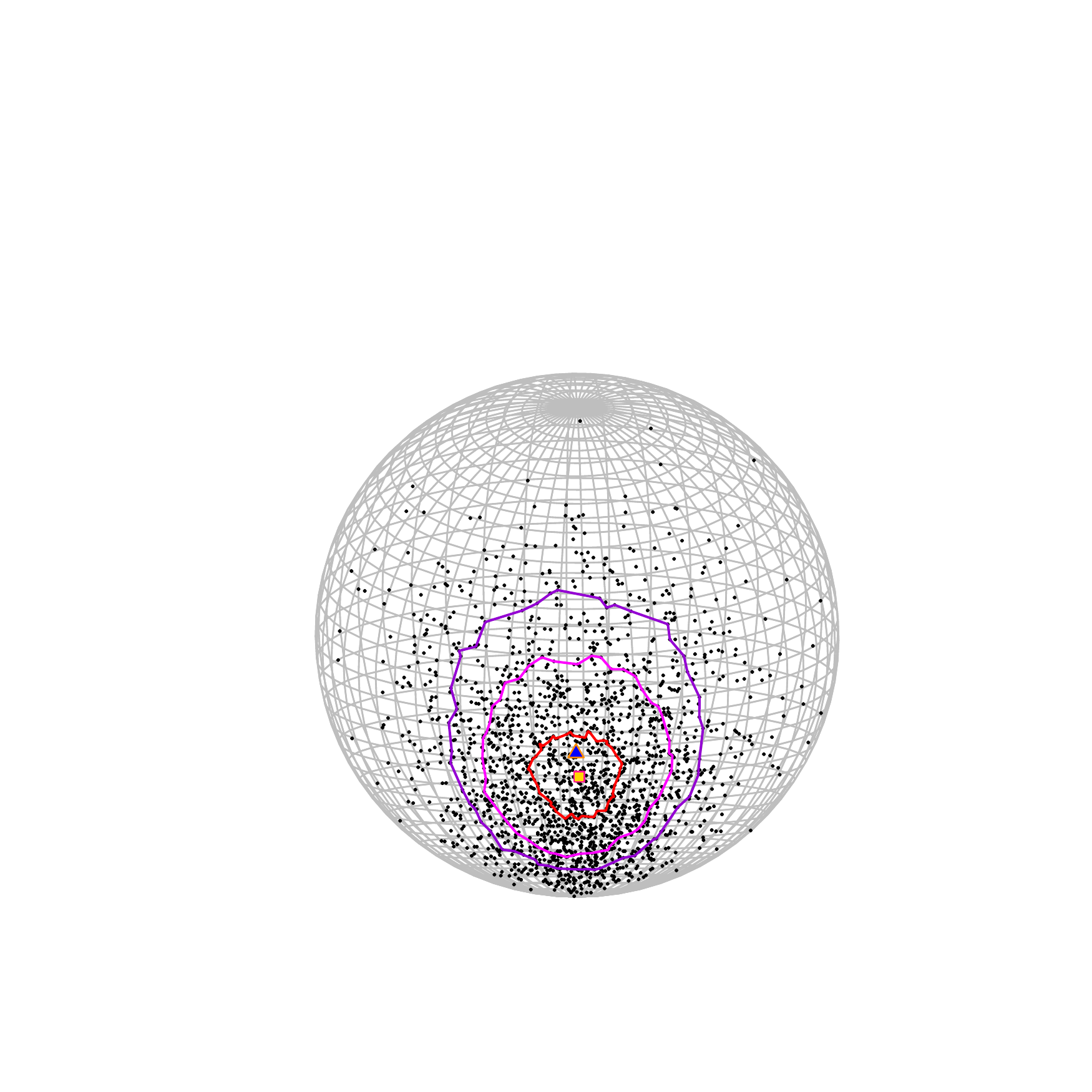} 
  \end{subfigure}  
          \begin{subfigure}[b]{0.31\textwidth}
         \vspace{-7mm}     
       \includegraphics[trim=20 10 25 50,clip, width=1.15 \textwidth ]{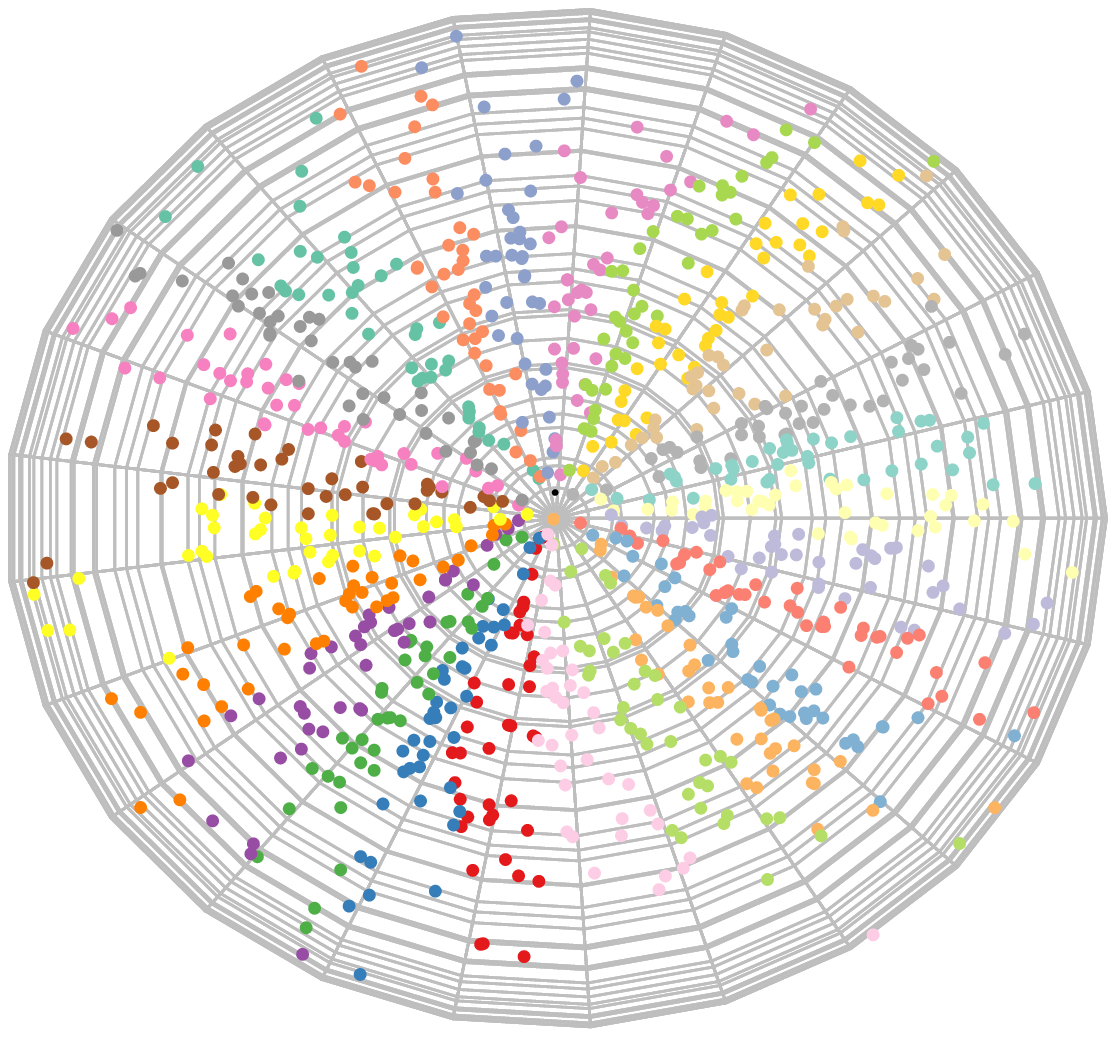}
     \end{subfigure}
     \hfill
     \begin{subfigure}[b]{0.31\textwidth}
     \hspace{-10mm}
       \vspace{-0mm}
         \includegraphics[trim=20 10 25 30,clip, width=1.12 \textwidth]{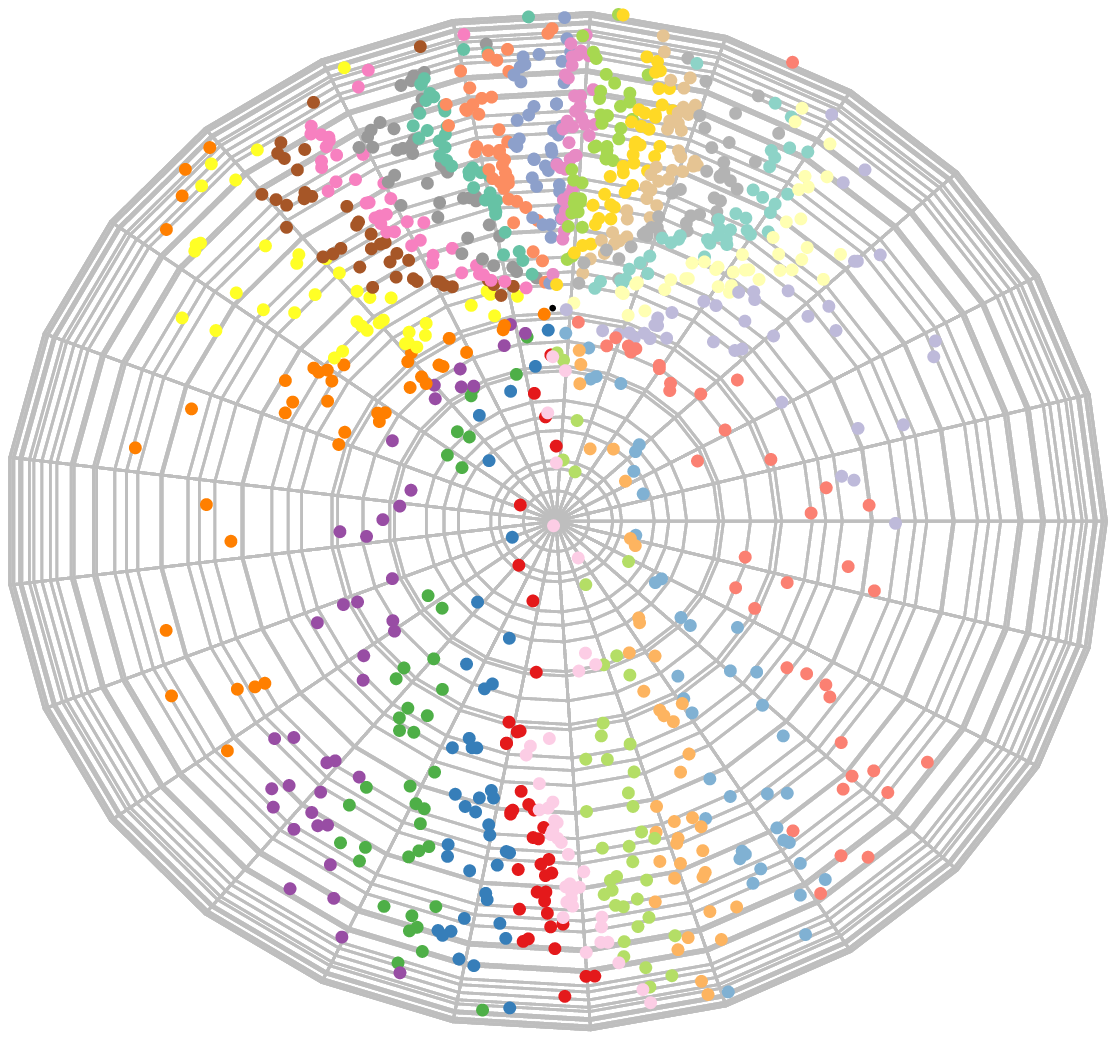}
     \end{subfigure}
     \begin{subfigure}[b]{0.3\textwidth}
      \hspace{-10.2mm}
       \vspace{0mm}
         \includegraphics[trim=30 10 25 78,clip, width=1.27 \textwidth]{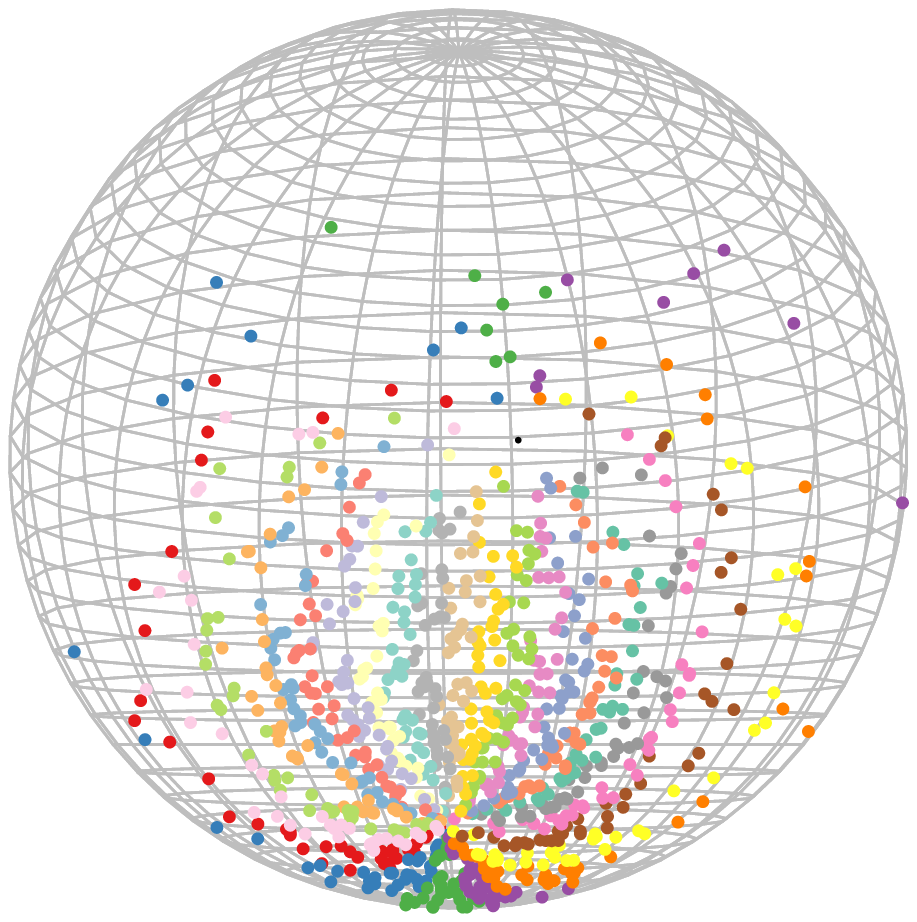}
     \end{subfigure}
     
 \caption{Upper panel: empirical quantile contours (probability contents 12.2\%, 48.8\%, and~70.7\%, respectively) computed from~$n=2001$ ($n_R = 40$, $n_S = 50$ and $n_0 =1$) points drawn from 
a von Mises-Fisher distribution (left), a mixture of two von Mises-Fisher distributions (middle) and a tangent von Mises-Fisher distribution (right). Lower panel: the corresponding empirical meridians; points with the same color have the same signs.}
 \label{Fig:Contours1}
\end{figure}

 For each simulation scheme, we computed the Fr\'echet mean via the {\tt rbase.robust} function from the R \linebreak package~{\tt RiemBase}. {The optimal coupling between the sample points and the regular grid has been obtained by using the fast network simplex algorithm (FNSA) of \cite{Bonneel2011} as implemented in the {\tt transport} function of the R package {\tt transport}. Many efficient algorithms have been proposed for solving the optimal transport problem, e.g.,  auction algorithm and its refined version \citep{Bertsekas1988, Bertsekas1992}. We used the FNSA which is reported in \cite{Bonneel2011} to reduce the complexity of the optimal transport problem from $O(n^3)$ to $O(n^2)$.} We factorized~$n$ into~ $n_R n_S + n_0$ with~$n_R = 40$, $n_S = 50$ and $n_0 =1$: each empirical quantile contour then consists in 50 points.  Plots of the median~$\thetat$, the empirical quantile contours for probability contents 12.2\% (rank~=~\!5), 48.8\% (rank = 20), and~70.7\% (rank = 29)  are shown in Figures~\ref{Fig:Contours1}  (left, middle, and right panels for the simulation schemes (i), (ii), and (iii),  respectively). We also provide in the same figures the corresponding empirical meridians (points with the same signs, \vspace{-0mm}
 same color in the figures).

 Under the vMF distribution (Figure~\ref{Fig:Contours1}, left panel),  the empirical quantile contours are symmetrically distributed around the median, as expected from the rotational symmetry of the   vMF distribution. Under the mixture distribution   (Figure~\ref{Fig:Contours1},  middle panel), quantile contours    adapt to the underlying multimodality, and the median is located in the mixture component with larger probability weight. Under the tangent vMF distribution, where skewness  is involved (Figure~\ref{Fig:Contours1},   right panel), the empirical quantile contours exhibit a distinctive skewed shape. Our empirical quantile contours thus  nicely pick up the shapes of the underlying distributions while controlling the probability contents of the corresponding quantile regions. \vspace{-0mm}

\section{Directional Goodness-of-fit} \label{sec:gof}\setcounter{equation}{0}

Classical  goodness-of-fit tests (Kolmogorov-Smirnov, Cram\' er-von~Mises, etc.) are based on distances between distribution functions. Our concepts of   population and empirical directional distribution functions quite naturally lead to the construction of  directional versions based on distances between ${\bf F}\n$ and ${\bf F}$. Assuming that a sample of \mbox{i.i.d.} random unit vectors $\Zb_1\n, \ldots, \Zb_n\n$ have common distribution ${\rm P}^\Zb$, the goodness-of-fit problem consists in testing the null hypothesis ${\cal H}_0: {\rm P}^\Zb= {\rm P}_0$ against~${\cal H}_1: {\rm P}^\Zb \neq {\rm P}_0$, where ${\rm P}_0$ is some specified distribution. This   has been studied extensively in directional statistics, with special attention to the problem of testing uniformity---in which case ${\rm P}_0$ is the uniform distribution $ {\rm P}^{\bf U}$ over $\mathcal{S}^{d-1}$. In view of Proposition \ref{characterize},   testing ${\cal H}_0\!:\! {\rm P}^\Zb\!=~\!\!{\rm P}_0$ against~${\cal H}_1\!:\! {\rm P}^\Zb\! \neq~\!\! {\rm P}_0$ is equivalent to testing ${\cal H}_0\!:\! {\bf F}\!=\! {\bf F}_0$ against~${\cal H}_1\!:~\!\!{\bf F}\! \neq~\!\!{\bf F}_0$ where ${\bf F}_0$ denotes the distribution function of ${\rm P}_0$---i.e., ${\bf F}_0(\zb)= \zb$  in case ${\rm P}_0= {\rm P}^{\bf U}$.

The test 
 we are proposing in this Section is a Cram\' er-von Mises-type test that rejects the null hypothesis for large values of the test statistic\vspace{-2.9mm}
\begin{equation}\label{CvM}
T_n\coloneqq  n^{-1} \sum_{i=1}^n \left\| {\bf F}\n (\Zb_i\n)- {\bf F}_0 (\Zb_i\n) \right\|^2,\vspace{-2mm}
\end{equation}
where ${\bf F}\n$ is obtained as in   Section~\ref{sec:emp} (the construction considered in Section~\ref{sec:emp1} here is sufficient). 
That test statistic is the empirical counterpart of the (squared) $L_2$ distance~${\rm E}\left[\|{\bf F}(\Zb_i\n) -{\bf F}_0 (\Zb_i\n)\|^2\right].$
As always in goodness-of-fit tests, the critical value $c_\alpha$ such that ${\rm P}[T_n> c_\alpha]= \alpha$ under ${\cal H}_0$ is easily approximated via Monte Carlo simulations since the null hypothesis is simple; the resulting test thus has exact size~$\alpha$. The following result establishes its  consistency.\vspace{-1mm} 

\begin{proposition}\label{Consist}
Assume that $\Zb_1\n, \ldots, \Zb_n\n$ are \mbox{i.i.d.} with common optimal transport map ${\bf F}$. Then,  provided that the uniform discrete distribution over the $n$-points  grid~${\mathfrak{G}}\n$ converges weakly,  as $n\to\infty$,  to the uniform distribution~${\rm P}^{\Ub}$ over~$ \mathcal{S}^{d-1}$, 
\begin{compactenum}
\item[(i)]  $T_n= o_{\rm P}(1)$ as $\ny$ under ${\cal H}_0$, while 
\item[(ii)] $T_n$ converges in probability to a strictly positive constant as $\ny$  if ${\bf F} \neq {\bf F}_0$.\vspace{-1mm} 
\end{compactenum}\end{proposition}
It follows from Proposition \ref{Consist} that the test $\phi_n$ is asymptotically consistent against any fixed alternative and therefore  qualifies as a universally consistent  {\it  omnibus} test. \smallskip

\textbf{Simulations for $d=3$.} In order to investigate the finite-sample performances of our test of uniformity based on~\eqref{CvM}, we performed Monte Carlo size and power  comparisons  of  our test and  the    projected Cram\' er-von Mises   (PCvM),   projected Anderson-Darling(PAD), and   projected Rothman  (PRt) tests 
recently proposed in \cite{Garcia-Portugues2020b} and the Rayleigh, Bingham, Ajne, Gin\'e and Bakshaev tests of uniformity \citep{Rayleigh1919, Bingham1974, Ajne1968, Gine1975, Bakshaev2010}; these tests were  implemented in the R package {\tt sphunif}. To do so, we generated $N=1000$ independent samples of~$n=400$ \mbox{i.i.d.} unit random vectors with the following distributions: 
\begin{compactenum}
\item[(i)] the uniform distribution;
\item[(ii)] vMF distributions, all with location parameter $\thetab = (0, 0, 1)^\top$, with   concentration parameters $\kappa = 0.05$,  $\kappa = 0.1$, and  $\kappa = 0.5$;
\item[(iii)] tangent vMF distributions as in \eqref{tangentvMF}, all with  angular distribution $F_{*}$, loca\-tion~$\thetab = (0, 0, 1)^\top$,   skewness direc\-tion~${\mbf \mu} = (0, 1)^\top$,  and  skewness intensities $\kappa = 0.05$,  $\kappa = 0.1$, and $\kappa = 0.2$;
\item[(iv)] mixtures of two vMFs, of the form 
 $I[U\leq 0.5] \Zb_1 + I[U> 0.5] \Zb_2$,  
where $U \sim {\rm U}[0, 1]$,  $\Zb_1 \sim \mathcal{M}_3(\thetab_1, \kappa_1)$, and~$\Zb_2 \sim \mathcal{M}_3(\thetab_2, \kappa_2)$ are mutually independent,  with $\thetab_1 = (0, -0.3, \sqrt{0.91})^\top$ and~$\thetab_2 = (0.3, \sqrt{0.66}, 0.5)^\top$, and  concentration parameters $\kappa_1 = \kappa_2 = 0.1$, $\kappa_1 = \kappa_2 = 0.2$, and~$\kappa_1 = \kappa_2 = 0.3$;
\item[(v)] mixtures of two vMFs and a tangent vMF, of the form   
 $I[U\leq 0.5] \Zb_1 + I[0.5<  U < 0.75] \Zb_2 + I[U \geq 0.75] \Zb_3$, 
where  $U \sim {\rm U}[0, 1]$,  $\Zb_1$, $\Zb_2$ and $\Zb_3$ are mutually independent;, with ${\Zb}_1$ following the tangent vMF distribution as in  (iii) with   skewness intensities $\kappa_1 = 0.07$,  $\kappa_1 = 0.1$, and $\kappa_1 = 0.2$,   and~$\Zb_2 \sim \mathcal{M}_p(\thetab_2, \kappa_2)$ and $\Zb_3 \sim \mathcal{M}_p(\thetab_3, \kappa_3)$ as in scheme (iv) with 
concentration parameters $\kappa_2 = \kappa_3 = 0.07$, $\kappa_2 = \kappa_3 = 0.1$; and $\kappa_2 = \kappa_3 = 0.2$.
\end{compactenum}

In Table~\ref{Tab.UnifTest}, we report  the rejection frequencies (at  nominal level $\alpha = 0.05$), out of~$N =~\!1000$ replications, of the optimal transport-based Cram\' er-von Mises test described in \eqref{CvM} (OT), the PCvM,   PAD, PRt, Rayleigh, Bingham, Ajne, Gin\'e, and Bakshaev tests of uniformity. The critical values are obtained through $2000$ Monte Carlo replications. Inspection of Table~\ref{Tab.UnifTest} reveals that all rejection frequencies under the null hypothesis of a uniform distribution are close to the nominal level $0.05$. Under the rotationally  symmetric  vMF alternatives, all tests show similar performances, except for the Bingham and Gin\'e tests, which have significantly lower powers. Under non-rotationally  symmetric alternatives such as  the tangent vMF or mixture distributions, our test uniformly outperforms all its competitors. 

\begin{table}[h!]
\caption{\small Rejection frequencies of the OT, PCvM, PAD, PRt, Rayleigh, Bingham, Ajne, Gin\'e and Bakshaev tests of uniformity over $\mathcal{S}^2$ under the uniform, vMF, tangent vMF, mixtures of two vMF distributions and mixtures of two vMFs and a tangent vMF, with various values of concentration or intensity parameter (simulation settings:~$N=1000$ replications, sample size $n=400$)}\label{Tab.UnifTest}
\centering
\small
\begin{tabular}{llllllllll} \hline \hline
                                                       & OT    & PCvM  & PAD   & PRt   & Rayleigh & Bingham & Ajne  & Gin\'e & Bakshaev \\ \hline
Uniform                                                & 0.045 & 0.054 & 0.054 & 0.054 & 0.053    & 0.060    & 0.051 & 0.059    & 0.054    \\
vMF ($\kappa = 0.05$)                                     & 0.078 & 0.076 & 0.073 & 0.072 & 0.070     & 0.052   & 0.072 & 0.053    & 0.076    \\
vMF ($\kappa = 0.1$)                                      & 0.134 & 0.137 & 0.139 & 0.140  & 0.139    & 0.053   & 0.137 & 0.055    & 0.137    \\
vMF ($\kappa = 0.5$)                                      & 0.990  & 0.998 & 0.998 & 0.998 & 0.998    & 0.082   & 0.998 & 0.082    & 0.998    \\
tangent vMF ($\kappa = 0.05$)                             & 0.073 & 0.042 & 0.042 & 0.044 & 0.050     & 0.031   & 0.051 & 0.03     & 0.042    \\
tangent vMF ($\kappa = 0.1$)                              & 0.202 & 0.138 & 0.133 & 0.138 & 0.137    & 0.027   & 0.143 & 0.029    & 0.138    \\
tangent vMF ($\kappa = 0.2$)                              & 0.669 & 0.578 & 0.575 & 0.573 & 0.580     & 0.030    & 0.586 & 0.028    & 0.578    \\
\begin{tabular}{@{}c@{}}Mixture of two vMFs \vspace{-1mm} \\ ($\kappa_1 = \kappa_2 = 0.1$) \end{tabular}  & 0.122 & 0.099 & 0.100   & 0.098 & 0.105    & 0.055   & 0.100   & 0.059    & 0.099    \\                       
\begin{tabular}{@{}c@{}}Mixture of two vMFs \vspace{-1mm} \\ ($\kappa_1 = \kappa_2 = 0.2$)  \end{tabular} & 0.325 & 0.303 & 0.302 & 0.297 & 0.295    & 0.058   & 0.297 & 0.061    & 0.303    \\
\begin{tabular}{@{}c@{}}Mixture of two vMFs \vspace{-1mm} \\ ($\kappa_1 = \kappa_2 = 0.3$)  \end{tabular}  & 0.628 & 0.615 & 0.611 & 0.615 & 0.621    & 0.042   & 0.624 & 0.046    & 0.615    \\
\begin{tabular}{@{}c@{}}Mixture of two vMFs \vspace{-1mm} \\  and a tangent vMF \vspace{-1mm} \\  ($\kappa_1 = \kappa_2 = \kappa_3 = 0.07$) \end{tabular} & 0.117 & 0.108 & 0.112 & 0.108 & 0.102    & 0.074   & 0.104 & 0.077    & 0.108    \\
\begin{tabular}{@{}c@{}}Mixture of two vMFs \vspace{-1mm} \\  and a tangent vMF \vspace{-1mm} \\  ($\kappa_1 = \kappa_2 = \kappa_3 = 0.1$) \end{tabular}  & 0.272 & 0.203 & 0.205 & 0.206 & 0.201    & 0.064   & 0.202 & 0.061    & 0.203    \\
\begin{tabular}{@{}c@{}}Mixture of two vMFs \vspace{-1mm} \\  and a tangent vMF \vspace{-1mm} \\  ($\kappa_1 = \kappa_2 = \kappa_3 = 0.2$) \end{tabular}  & 0.707 & 0.609 & 0.612 & 0.606 & 0.598    & 0.083   & 0.604 & 0.082    & 0.609 \\   \hline \hline
\end{tabular}
\end{table}

\textbf{Simulations for $d=2$.} We investigate, for $d=2$, the finite-sample performance 
of the same tests   by genera\-ting~$N=1000$ replications of $n=100$ i.i.d.~random vectors from
\begin{compactenum}
\item[(i)] the uniform distribution;
\item[(ii)] vMF distributions  with location parameter $\thetab = (0, 1)^\top$, and  concentration parameters $\kappa = 0.05$,  $\kappa = 0.1$, and  $\kappa = 0.5$;
\item[(iii)]  mixtures of two vMFs, of the form 
 $I[U\leq 0.7] \Zb_1 + I[U> 0.7] \Zb_2$,  
where $U \sim {\rm U}[0, 1]$,  $\Zb_1 \sim \mathcal{M}_2(\thetab_1, \kappa_1)$, and~$\Zb_2 \sim \mathcal{M}_2(\thetab_2, \kappa_2)$ are mutually independent,  with $\thetab_1 = (-0.3, \sqrt{0.91})^\top$ and~$\thetab_2 = (0.6, 0.8)^\top$, and   concentration parameters $\kappa_1 = \kappa_2 = 0.1$, $\kappa_1 = \kappa_2 = 0.25$, and~$\kappa_1 = \kappa_2 = 0.5$;
\item[(iv)] the sine-skew distribution proposed by \cite{umbach2009building} and \cite{abe2011sine}, the density of which takes form
 $\phi \mapsto f(\phi-\mu)(1+\lambda \sin(\phi-\mu), \quad \phi \in [-\pi, \pi),$ 
where $f$ is a density (over $ [-\pi, \pi)$) symmetric about zero, $\mu\in  [-\pi, \pi)$ is an angular location, and $\lambda \in (-1, 1)$  a skewness parameter. See \citet[Chapter~2]{Ley2017a} for the data-generating process. In the simulations, we chose $f$ to be a vMF density symmetric about $0$ with concentration parameters $\kappa = 0.1$, set $\mu = 0$ and $\lambda = 0.1$, 0.3, and 0.35.\vspace{-0mm}
\end{compactenum}

In Table~\ref{Tab.UnifTestS1}, we report the rejection frequencies (at  nominal level $\alpha = 0.05$) of the OT, PCvM, PAD, PRt, Rayleigh, Bingham, Ajne, Gin\'e, and Bakshaev tests of uniformity over $\mathcal{S}^1$. The results are consistent with those for $d=3$ in Table~\ref{Tab.UnifTest}. Specifically, all tests have rejection frequencies close to the nominal level under the null hypothesis of uniformity. Under the vMF alternatives, the  Bingham and Gin\'e tests have less power than all  other ones. Our test outperforms all its competitors whenever skewness or multi-modality   are present. {Further simulations yielding, in dimension $d = 5$, similar conclusions, are provided in the supplementary material}.

\begin{table}[t!]
\caption{\small Rejection frequencies of the OT, PCvM, PAD, PRt, Rayleigh, Bingham, Ajne, Gin\'e, and Bakshaev tests of uniformity over $\mathcal{S}^1$ under the uniform, vMF, mixtures of two vMF distributions and sine-skew distributions (simulation settings:~$N=1000$ replications, sample size $n=100$)}\label{Tab.UnifTestS1}
\centering
\small
\begin{tabular}{llllllllll} \hline \hline
                                                       & OT    & PCvM  & PAD   & PRt   & Rayleigh & Bingham & Ajne  & Gin\'e & Bakshaev \\ \hline
Uniform   & 0.058 & 0.053 & 0.055 & 0.052 & 0.054    & 0.039   & 0.054 & 0.042    & 0.054    \\
vMF ($\kappa = 0.05$)          & 0.080  & 0.069 & 0.068 & 0.068 & 0.067    & 0.054   & 0.067 & 0.051    & 0.070     \\
vMF ($\kappa = 0.1$)         & 0.105 & 0.098 & 0.097 & 0.096 & 0.096    & 0.045   & 0.098 & 0.048    & 0.097    \\
vMF ($\kappa = 0.5$)       & 0.877 & 0.873 & 0.867 & 0.876 & 0.875    & 0.064   & 0.873 & 0.069    & 0.875    \\ 
\begin{tabular}{@{}c@{}}Mixture of two vMFs \vspace{-1mm} \\ ($\kappa_1 = \kappa_2 = 0.1$) \end{tabular}   & 0.092 & 0.080  & 0.079 & 0.079 & 0.083    & 0.045   & 0.084 & 0.043    & 0.078    \\
\begin{tabular}{@{}c@{}}Mixture of two vMFs \vspace{-1mm} \\ ($\kappa_1 = \kappa_2 = 0.25$) \end{tabular}   & 0.288 & 0.256 & 0.249 & 0.256 & 0.264    & 0.055   & 0.262 & 0.053    & 0.260     \\
\begin{tabular}{@{}c@{}}Mixture of two vMFs \vspace{-1mm} \\ ($\kappa_1 = \kappa_2 = 0.5$) \end{tabular}   & 0.847 & 0.829 & 0.822 & 0.831 & 0.831    & 0.057   & 0.830  & 0.058    & 0.827    \\
Sine-skew ($\lambda=0.1$)  & 0.093 & 0.080  & 0.079 & 0.084 & 0.081    & 0.041   & 0.081 & 0.041    & 0.080     \\
Sine-skew ($\lambda=0.3$) & 0.498 & 0.465 & 0.457 & 0.465 & 0.473    & 0.035   & 0.470  & 0.038    & 0.471    \\
Sine-skew ($\lambda=0.35$) & 0.637 & 0.618 & 0.610  & 0.624 & 0.627    & 0.046   & 0.627 & 0.049    & 0.622   \\ \hline \hline
\end{tabular}
\end{table}

\color{black}

\section{Directional MANOVA} \label{MANOVAsec}\setcounter{equation}{0}

Let $\Xb_{i1}^{(n_i)}, \ldots, \Xb_{in_i}^{(n_i)}, i = 1, \ldots, m,$  denote $m (\geq 2)$ independent samples on the unit hypersphere~$\mathcal{S}^{d-1}$. For each $i$, we assume that $\Xb_{i1}^{(n_i)}, \ldots, \Xb_{in_i}^{(n_i)}$ are i.i.d.~with  common absolutely continuous distribution ${\rm P}_i$,  density~$f_{i}$, and  directional distribution function $\Fb_i$, $i=1,\ldots, m$. The null hypothesis of interest is the hypothesis $\mathcal{H}_0: \Fb_{1} = \ldots = \Fb_{m} =: \Fb$ (un\-spe\-cified~$\Fb$) of no treatment effect.

  Let $\Fb\n$ denote the empirical directional distribution function computed from the pooled sample 
  $$\{\Yb \n_1, \ldots, \Yb\n_n\} \coloneqq \{\Xb_{11}^{(n_1)}, \ldots, \Xb_{1n_1}^{(n_1)}, \ldots, \Xb_{m1}^{(n_m)}, \ldots, \Xb_{mn_m}^{(n_m)}\}.$$
The tests we propose here are based on a $(md_{\bf J})$-dimensional statistic  ${\bDelta}_{\Jb}\n\coloneqq (({\bDelta}_{1; \Jb}^{(n_1)})^\top, \ldots, ({\bDelta}_{m; \Jb}^{(n_m)})^\top)^\top$ 
where, for some 
 score function~$\Jb:\mathcal{S}^{d-1} \to \mathbb{R}^{d_{\bf J}}$, \vspace{-2mm}
\begin{equation*}
{\bDelta}_{i; \Jb}^{(n_i)} \coloneqq 
n_i^{-1/2} \sum_{j = 1}^{n_i} \Jb(\Fb\n(\Xb_{ij}^{(n_i)})) - n^{-1/2} \sum_{\ell =1}^n \frac{n_i^{1/2}}{n^{1/2}} \Jb(\Fb\n(\Yb \n_{\ell})) \label{def.tildeDeltai}
=
 n^{-1/2} \sum_{\ell =1}^n (a_{i\ell}\n - \bar{a}_i\n) \Jb(\Fb \n(\Yb \n_{\ell})),\ \  i=1,\ldots, m
\end{equation*}
with 
 $
a_{i\ell}\n\!\coloneqq (n^{1/2}/n_i^{1/2})I[{n_1+\ldots+n_{i-1}+1 \leq \ell \leq n_1+\ldots+n_{i-1}+n_i}]$ and $\bar{a}_i\n\! = n^{-1} \sum_{\ell =1}^n a_{i\ell}\n = n_i^{1/2}/n^{1/2}
$   ($I[{.}]$, as usual, denotes the indicator function). 
Below, we assume that  $r_i\n\coloneqq  n_i/n$ converges to a constant $0<r_i<1$ as~$n \rightarrow \infty$ (possibly, along some subsequence). Letting~{$\Db_{\Jb}\coloneqq  {\rm Var}(\Jb (\Ub))$ where~$\Ub$ is uniform over~$\mathcal{S}^{d-1}\!$, our tests reject 
 ${\cal H}_0$ 
 for large values of 
 $
{Q}_{\Jb}\n \coloneqq  ({\bDelta}_{\Jb}^{(n)})^\top ({\bf I}_m \otimes {\bf D}_{\Jb}^{-}) {\bDelta}_{\Jb}\n,$  
where~${\bf D}_{\Jb}^{-}$ denotes the Moore-Penrose inverse of ${\bf D}_{\Jb}$.}
The next result shows that ${Q}_{\Jb}\n$ is asymptotically chi-square under  ${\mathcal H}_0$ 
 and some mild regularity assumptions on the score function ${\bf J}$.
\begin{proposition} \label{MANOVAprop}
Assume that (i) $\Jb $ is continuous over $\mathcal{S}^{d-1}$ and (ii) $\Jb$ is  square-integrable, that is, 
$\int_{\mathcal{S}^{d-1}} \|{\bf J}(\ub) \|^2 {\rm dP}^{\Ub}(\ub) < \infty$, and (iii) for any sequence $\mathfrak{G}^{(n)} \coloneqq \{\mathfrak{G}_1^{(n)}, \ldots, \mathfrak{G}_n^{(n)}\}$ of $n$-tuples in  $\mathcal{S}^{d-1}$ such that the uniform discrete distribution over $\mathfrak{G}^{(n)}$ converges weakly to ${\rm P}^{\Ub}$ as $n \rightarrow \infty$, 
$
\underset{n\rightarrow\infty}{\lim} n^{-1} \sum_{{\ell}=1}^n 
 \| {\bf J} (\mathfrak{G}_{\ell}^{(n)}) \|^2 =    \int_{\mathcal{S}^{d-1}} \| {\bf J}(\ub) \|^2 {\rm d P}^{\Ub}(\ub).
$ 
Then, under ${\cal H}_0$ as $\ny$, ${Q}_{\Jb}\n$ is asymptotically chi-square  with {$(m-1) d^*$ degrees of freedom, where $d^*$ is the rank of the $d_{\bf J}\times d_{\bf J}$ matrix~${\bf D}_{\Jb}$.}
\end{proposition}
The MANOVA tests $\phi_{\Jb}\n$ we are proposing reject the  hypothesis of no treatment effect at  asymptotic level $\alpha$ when\-ever~
${Q}_{\Jb}\n > \chi^2_{(m-1) d^*; 1- \alpha}$ 
where $ \chi^2_{p; \tau}$ stands for the quantile of order $\tau$ of a chi-square distribution with $p$ degrees of freedom. Being based on the $\Fb\n$-measurable test statistic ${Q}_{\Jb}\n$, these tests are fully distribution-free. Below, we investigate their asymptotic properties  under local alternatives. To do so, let us consider a parametric framework where the underlying distributions are indexed by some finite-dimensional parameter $\boldsymbol\omega  \in \bOmega \subset \R^k$. More precisely, denote by 
 ${\rm P}^{(n)}_{\underline\tbomega}$,  where~$\underline{\boldsymbol\omega  }\coloneqq (\boldsymbol\omega  _1, \ldots, \boldsymbol\omega  _m)\in \bOmega^m \subset \R^{mk}$, the joint distribution of the pooled sample $(\Yb \n_1, \ldots, \Yb\n_n)$ when the~$i$th sample~$(\Xb_{i1}^{(n_i)}, \ldots, \Xb_{in_i}^{(n_i)})$ has distribution ${\rm P}^{(n_i)}_{\tbomega_i}$, $i = 1, \ldots, m$. Clearly, under the null hypothesis,  ${\rm P}^{(n)}_{\underline\tbomega}$ is of the form~${\rm P}\n_{(\tbomega_0, \ldots, \tbomega_0)}$ for some $\boldsymbol\omega _0 \in \bOmega$: write  $\underline{\boldsymbol\omega }_0$ for $( \boldsymbol\omega_0, \ldots,  \boldsymbol\omega_0)$.

In order to study the local power of our tests $\phi_{\Jb}\n$, consider local alternatives of the form $\underline{\boldsymbol\omega }_0+  n^{-1/2}\bnu\n \btau\n \in~\!\bOmega^m$, where  $\bnu\n \coloneqq {\rm diag}((r_1\n)^{-1/2} {\bf I}_{k}, \ldots, (r_m\n)^{-1/2} {\bf I}_{k})$ and $\btau\n \coloneqq ((\btau_1^{(n_1)})^\top, \ldots, (\btau_m^{(n_m)})^\top)^\top$ with $\btau_i^{(n_i)}$ a bounded sequence of~$\mathbb{R}^{k}$. We will assume that the underlying sequence of experiments is sufficiently regular in the sense that,  
under~${\rm P}^{(n)}_{{\underline \tbomega}_0}$, as $n \rightarrow \infty$,\vspace{-3mm} 
\begin{equation}\label{eq.jointLik}
\Lambda\n\coloneqq \log\frac{ {\rm d}{\rm P}^{(n)}_{{\underline \tbomega}_0 +   n^{-1/2}\tbnu\n \tbtau\n} }{ {\rm d}{\rm P}^{(n)}_{{\underline \tbomega}_0}} = (\btau^{(n)})^\top \bDelta^{(n)}_{\underline{\tbomega}_0} - \frac{1}{2} (\btau^{(n)})^\top {\boldsymbol{\mathcal I}}_{\underline{\tbomega}_0} \btau^{(n)} + o_{\rm P}(1),
\end{equation}
where $\bDelta^{(n)}_{{\underline \tbomega}_0} \coloneqq ((\bDelta^{(n_1)}_{1; \tbomega_0})^\top, \ldots, (\bDelta^{(n_m)}_{m; \tbomega_0})^\top)^\top$, with $\bDelta^{(n_i)}_{i; \tbomega_0}$ of the form 
$
\bDelta^{(n_i)}_{i; \tbomega_0}  \coloneqq n_i^{-1/2}   \sum_{j = 1}^{n_i} \bvp_{\tbomega_0}(\Xb_{ij}^{(n_i)})$,  $ i=1,\ldots, m
$ 
for some square-integrable  function $\bvp_{\tbomega_0}$, and ${\boldsymbol{\mathcal I}}_{{\underline \tbomega}_0} \coloneqq {\rm diag}({\boldsymbol{\mathcal I}}_{1; \tbomega_0}, \ldots, {\boldsymbol{\mathcal I}}_{m; \tbomega_0})$ are such that, still under ${\rm P}^{(n)}_{{\underline \tbomega}_0}$,  $\bDelta^{(n)}_{{\underline \tbomega}_0}$ is asymptotically normal with mean ${\boldsymbol 0}$ and covariance matrix ${\boldsymbol{\mathcal I}}_{{\underline \tbomega}_0}$. That is, we  assume that the underlying model is LAN with central sequence $\bDelta^{(n)}_{{\underline \tbomega}_0}$ and information matrix ${\boldsymbol{\mathcal I}}_{\underline{\tbomega}_0}$.
  Letting
 ${\bf K}_{{\bf J}, \tbomega_0} \coloneqq \int_{\mathcal{S}^{d-1}}   {\bf J} (\ub)  \bvp_{\tbomega_0}({\bf F}^{-1}(\ub)) {\rm dP}^{\Ub}(\ub),$ 
we have the following result.
\begin{proposition} \label{MANOVAproplocal}
Under ${\rm P}\n_{{\underline \tbomega}_0+   n^{-1/2}\tbnu\n \tbtau\n}$, ${Q}_{\Jb}\n$ is asymptotically non-central chi-square with $(m-1){d}^*$ \vspace{-2mm}degrees of freedom and, letting $ \btau_j\coloneqq  \lim_{n \rightarrow \infty} \btau_j\n$,  non-centrality parameter 
$\displaystyle{\sum_{j=1}^m (1- r_j) \btau_j^\top {\bf K}_{{\bf J}, \tbomega_0}^\top {\bf D}_{{\bf J}}^{-} {\bf K}_{{\bf J}, \tbomega_0} \sum_{j=1}^m (1- r_j) \btau_j .}$
\end{proposition}

{Proposition \ref{MANOVAproplocal} and Le Cam's third lemma readily yield integral expressions (involving the scores $\bf J$, the directional distribution function ${\bf F}^{(n)}$, and $\bvp_{\tbomega_0}$) for  local and asymptotic powers in this class of parametric models. They also help select a score function~$\bf J$: the test based on ${\bf J}= \bvp_{\tbomega_0}\circ\left( {\bf F}^{(n)}\right)^{-1}$, for instance, achieves optimality in the parametric LAN model just described.}  Rather than listing such theoretical expressions 
 for specific LAN models and specific alternatives (location, concentration, multimodality, skewness, ... ), we rather provide here a ``universal'' consistency result for 
  $\phi_{\Jb}\n$.
 \begin{proposition} \label{consistency}
 Let $\Xb_{i1}^{(n_i)}, \ldots, \Xb_{in_i}^{(n_i)}$, $i=1, \ldots, m$ be mutually independent samples of \mbox{i.i.d} random directions with directional distribution functions $\Fb_i$, $i=1,\ldots, m$ respectively.  Assume that there exists a couple $(i, j)$ such that  $\Fb_i \neq  \Fb_j$. Then~$\lim_{\ny} {\rm E}[\phi_{\Jb}\n]=1$.
 \end{proposition}
We conclude this Section by investigating the finite-sample performances of $\phi_{\Jb}\n$ for various scores: 
\begin{compactenum}  
\item[(i)] the uniform score ${\Jb}({\ub}) \coloneqq \ub$; 
\item[(ii)] the estimated vMF-location score $\Jb(\Fb\n(\Yb\n_{\ell})) \coloneqq \hat{\kappa} \sqrt{1- \Big( G_{\hat \kappa}^{-1}\Big(1-\frac{R\n_{\ell}}{n_R+1} \Big) \Big)^2}  \,\mathbf{S}_{\ell}\n$,  
where $\hat{\kappa}$ is the vMF maximum-likelihood estimator of the concentration parameter, $G_\kappa$  the distribution function of $\Zb^\top \thetab$ with $\Zb \sim  {\cal M}_d(\bth, \kappa)$,
 $R\n_{\ell} \coloneqq  R\n(\Yb_{\ell}\n) =  (n_R+1)\left [1 - F_{*}(({\bf F}\n(\Yb_{\ell}\n))^\top\thetat))\right]$,  and 
$\mathbf{S}_{\ell}\n \coloneqq  \mathbf{S}_{\thetatsub}(\Fb\n(\Yb\n_{\ell}))$;  the estimator~$\thetat$ is obtained as in Section~\ref{ranksandsigns},  step 1 of the construction of the grid;
\item[(iii)] the estimated vMF-concentration score 
$
\Jb(\Fb\n(\Yb\n_{\ell})) \coloneqq  G_{\hat \kappa}^{-1}\Big(1-\frac{R\n_{\ell}}{n_R+1} \Big),
$ 
and 
\item[(iv)] the estimated vMF-location-concentration score (a linear combination of the scores in (ii) and (iii))\vspace{-0mm}
\[
\Jb(\Fb\n(\Yb\n_{\ell})) \coloneqq \hat{\kappa} \Bigg\{ G_{\hat \kappa}^{-1}\Big(1-\frac{R\n_{\ell}}{n_R+1} \Big) \, \thetat  + \sqrt{1- \Big( G_{\hat \kappa}^{-1}\Big(1-\frac{R\n_{\ell}}{n_R+1} \Big) \Big)^2}  \,\mathbf{S}_{\ell}\n \Bigg\}.\vspace{-0mm}
\]
\end{compactenum}

In the simulation exercise below, we compare the resulting tests with their pseudo-vMF (pvMF) counterpart, which  rejects the null hypothesis at  asymptotic level $\alpha$ whenever $Q\n> \chi^2_{(m-1) (d-1); 1- \alpha}$, where\vspace{-1mm} 
\begin{align*}
Q\n & \coloneqq  (d-1) \Bigg(\sum_{i=1}^m \frac{n_i D_{i}}{E_{i}} (\bar{\Xb}_{i}^{(n_i)})^\top ({\bf I}_d - \hat{\bth} \hat{\bth}^\top) \bar{\Xb}_{i}^{(n_i)}- \sum_{i, j}^m \frac{n_i n_j}{n} \frac{D_{i} D_{j}}{H} (\bar{\Xb}_{i}^{(n_i)})^\top ({\bf I}_d - \hat{\bth} \hat{\bth}^\top) \bar{\Xb}_{j}^{(n_j)} \Bigg),
\end{align*}
with $\hat{\bth}$   the sample Fr\'echet mean,  $\bar{\Xb}_{i}^{(n_i)}\coloneqq  n_i^{-1} \sum_{j=1}^{n_i} \Xb_{ij}^{(n_i)}$, $E_{i}\coloneqq  n_i^{-1} \sum_{j=1}^{n_i} (\Xb_{ij}^{(n_i)})^\top  \hat{\bth}$, $H\coloneqq  \sum_{i=1}^m r_i\n D_{i}^2 B_{i}$, $D_{i} \coloneqq E_{i}/B_{i}$, and~$B_{i} \coloneqq  1- n_i^{-1} \sum_{j=1}^{n_i} ((\Xb_{ij}^{(n_i)})^\top  \hat{\bth})^2$; see \cite{Ley2017} for details.

 \begin{figure}[b!]
     \centering
     \begin{subfigure}[b]{0.45\textwidth}
         \centering
         \includegraphics[scale=0.27]{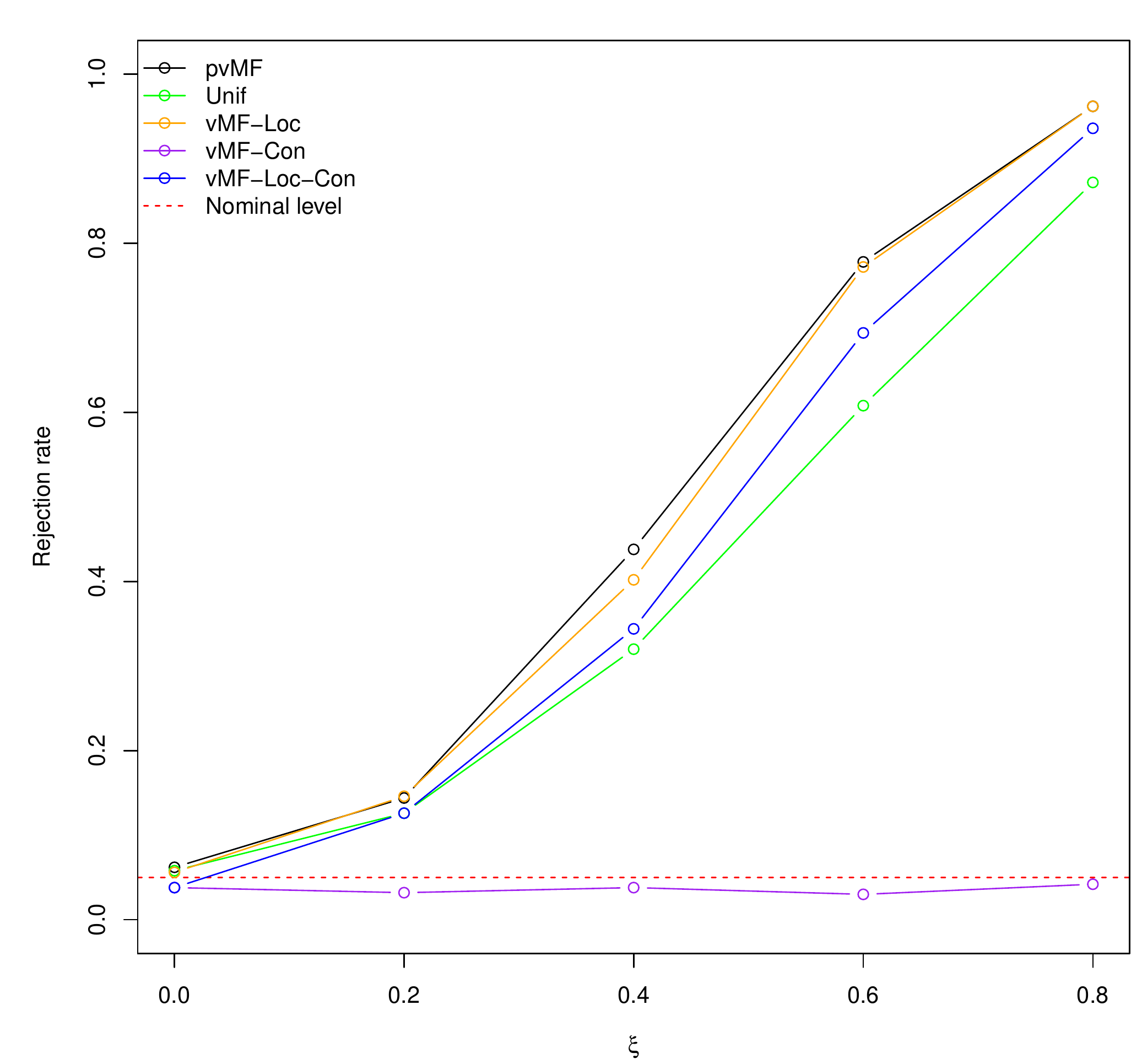} \vspace{-5mm}
     \end{subfigure}
     \begin{subfigure}[b]{0.45\textwidth}
         \centering
         \includegraphics[scale=0.27]{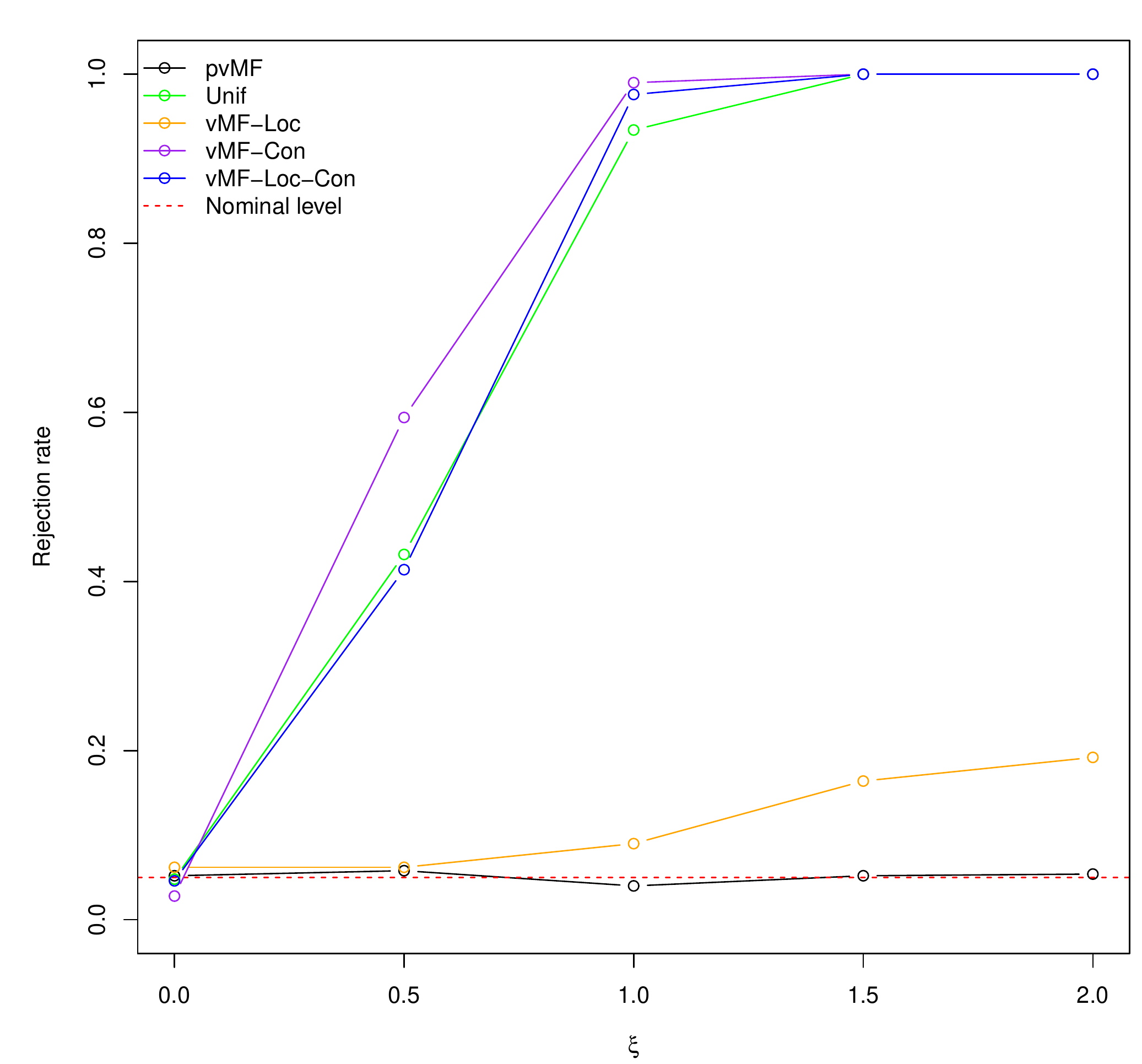} \vspace{-5mm}
     \end{subfigure}
     \begin{subfigure}[b]{0.45\textwidth}
         \centering
         \includegraphics[scale=0.27]{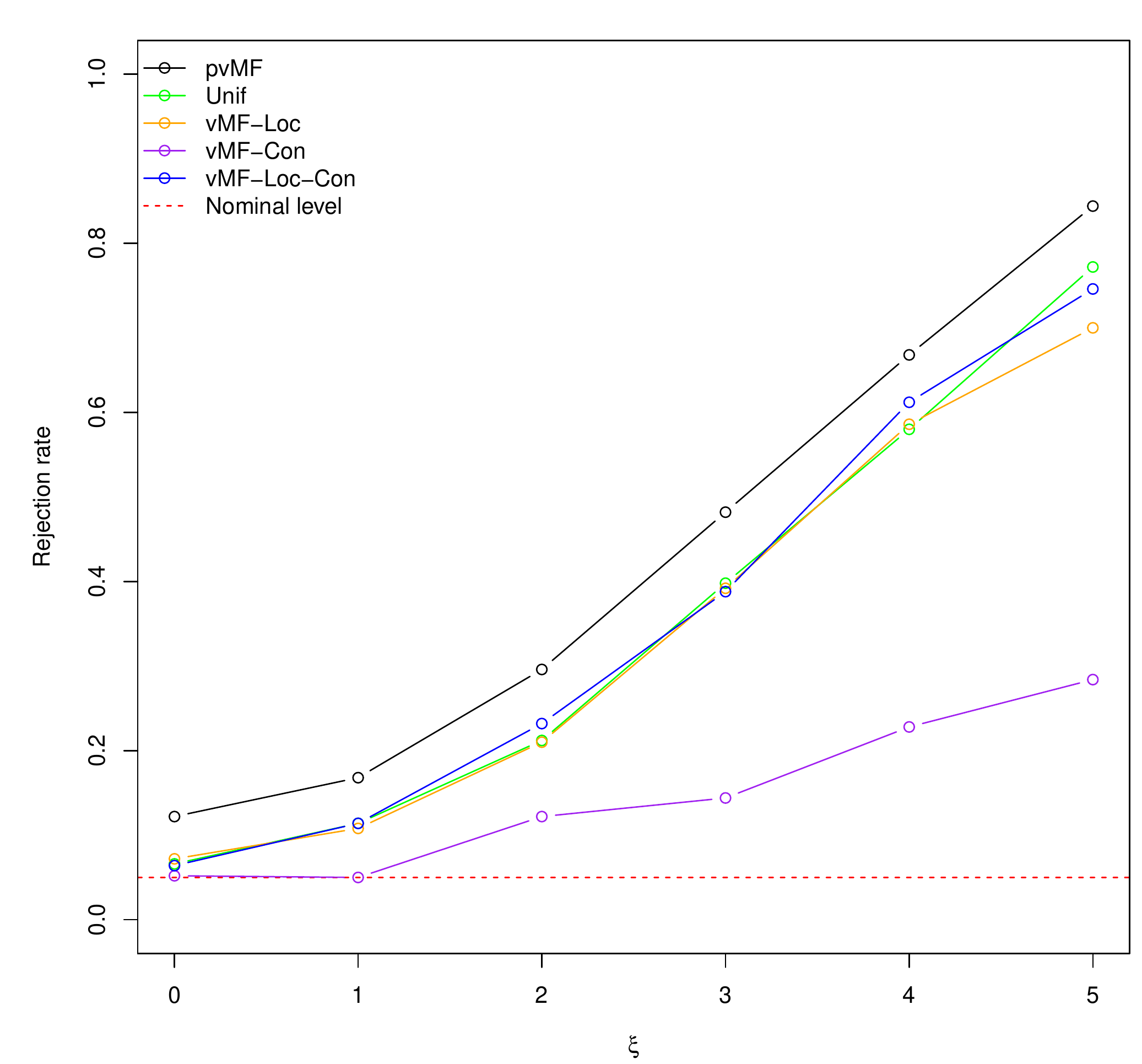}  
     \end{subfigure}
     \begin{subfigure}[b]{0.45\textwidth} 
         \centering
         \includegraphics[scale=0.27]{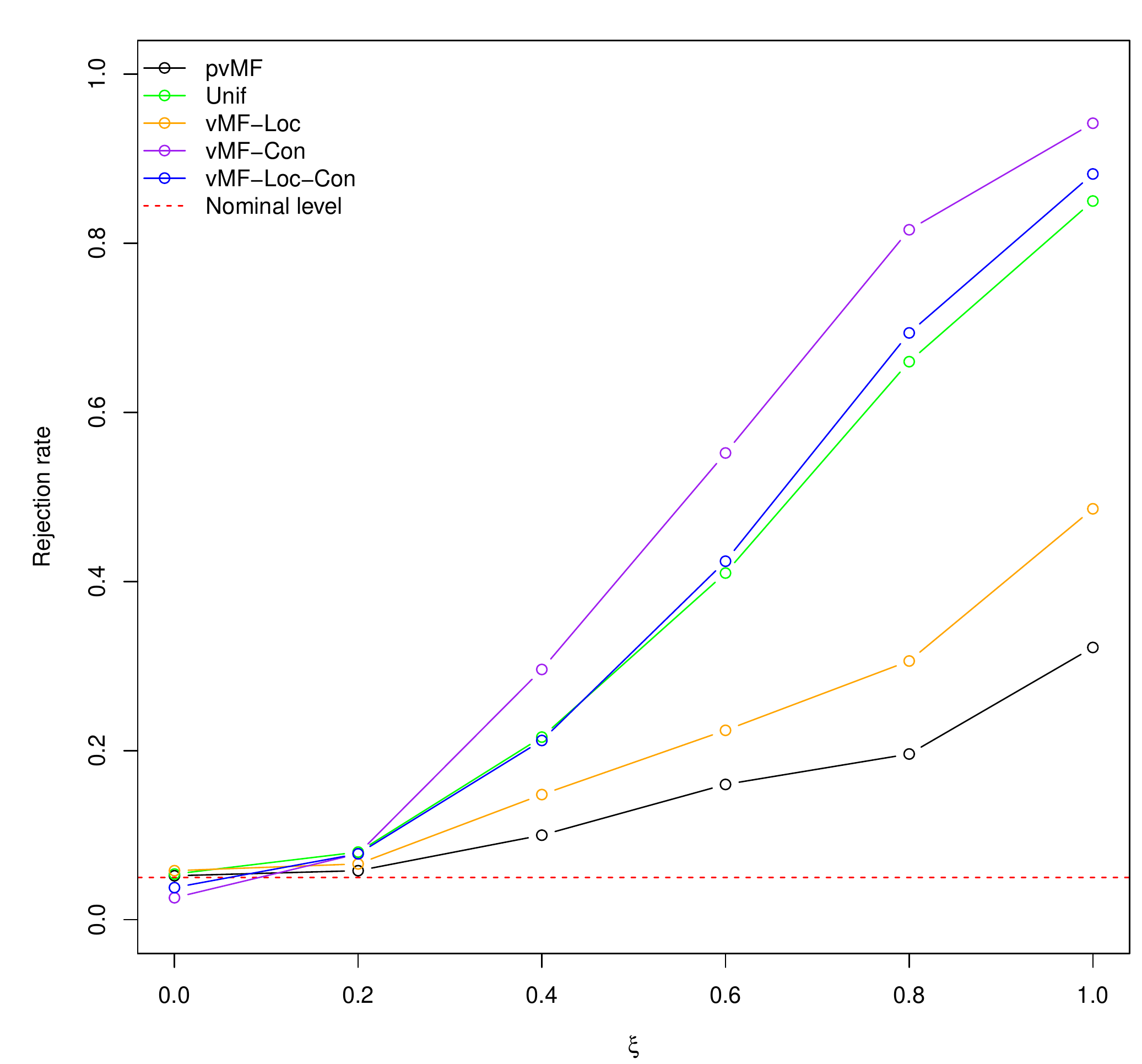} 
     \end{subfigure}
        \caption{\small Rejection rates (at $0.05$ nominal level) of the pseudo-vMF and rank-based  tests with uniform, vMF-location, vMF-concentration, and vMF-location-concentration scores,  for Cases (1) (top-left panel), (2) (top-right panel), (3) (bottom-left panel), and (4) (bottom-right panel).}
 \label{Fig:RejectRate}
\end{figure}

In the Monte-Carlo simulation, we set $d=3$, $m=2$, $n_1 = 500$, $n_2 = 600$, hence $n=1100$. The structured regular  grid for the computation of the vMF-location,  vMF-concentration, {and  vMF-location-concentration scores} was based on $n_R = 44$, $n_S = 25$, and $n_0 = 0$. 
The following data-generating processes were considered.  

\begin{compactenum}
\item[Case (1)]\hspace{-1mm} (location alternatives): $\Xb_{1j}^{(n_1)} \!\!\sim\! {\cal M}_d(\bth, \kappa)$, $j = 1, \ldots, n_1$ and $\Xb_{2\ell}^{(n_2)} \!\!\sim\! {\cal M}_d({\bf O}_{\xi} \bth, \kappa)$, $\ell = 1, \ldots, n_2$ where $\kappa=3$,\vspace{-1mm}
\begin{equation}\label{Oxi}
 \bth = \left(
\begin{array}{c}1\\ 0\\ 0\end{array}\right),
\quad\text{ and}\quad 
{\bf O}_{\xi} = 
\begin{pmatrix}
\cos(\pi\xi/15) & -\sin(\pi\xi/15) & 0 \\
\sin(\pi\xi/15) & \cos(\pi\xi/15) & 0  \\
0 & 0 & 1
\end{pmatrix}
\quad\text{ with } \xi = 0, \ 0.2, \  0.4,\  0.6,  \ 0.8;\vspace{-1mm}
\end{equation}
\item[Case (2)]\hspace{-1mm} (concentration alternatives): $\Xb_{1j}^{(n_1)} \!\!\sim\! {\cal M}_d(\bth, \kappa)$, $j = 1, \ldots, n_1$ and $\Xb_{2\ell}^{(n_2)} \!\!\sim\! {\cal M}_d(\bth, \xi + \kappa)$, $\ell = 1, \ldots, n_2$ where~$\kappa=3$ and $\bth = (1, 0, 0)^\top$, with~$\xi = 0, \ 0.5, \ 1, \ 1.5,  \ 2$;
\item[Case (3)]\hspace{-1mm} (multimodal  alternatives): $\Xb_{1j}^{(n_1)}$, $j = 1, \ldots, n_1$, is a mixture, with mixing probabilities $\frac{3}{8}$,  $\frac{3}{8}$, and $\frac{1}{4}$, of  ${\cal M}_d(\bth_1, \kappa_1)$, ${\cal M}_d(\bth_2, \kappa_2)$,  and~${\cal M}_d(\bth_3, \kappa_3)$,  where $\bth_1 = (1, 0, 0)^\top$, $\kappa_1 = 3$, $\bth_2 = (-0.8, 0.3, \sqrt{0.27})^\top$, $\kappa_2 = 2$, $\bth_3 = (0, -0.7, \sqrt{0.51})^\top$, $\kappa_3 = 3$, 
 and 
$\Xb_{2\ell}^{(n_2)} = {\bf O}_{\xi} \Zb_{\ell}^{(n_2)}$, $\ell = 1, \ldots, n_2$, where $\Zb_{\ell}^{(n_2)}$ has the same distribution as $\Xb_{1j}^{(n_1)}$ and ${\bf O}_{\xi}$ is as in \eqref{Oxi} with $\xi = 0, 1, \ldots, 5$;
\item[Case (4)]\hspace{-1mm} (skewed alternatives):  $\Xb_{1j}^{(n_1)}$, $j=1,\ldots, n_1$ and $\Xb_{1j\ell}^{(n_2)}$, $\ell=1,\ldots, n_1$ are generated from a {tangent vMF distribution} of the form \eqref{tangentvMF} with $\thetab = (0, 0, 1)^\top$,  ${\mbf \mu} = (0.7, \sqrt{0.51})^\top$, $\kappa = 1$,  with~$\widetilde{V}\sim\text{\rm Beta}(2, 5)$ for $\Xb_{1j}^{(n_1)}$ and~$\widetilde{V}\sim\text{\rm Beta}(2, 5+\xi)$, $\xi = 0, 0.2, \ldots, 1$ for $\Xb_{2\ell}^{(n_2)}$. 
\end{compactenum}

\vspace{3mm}

For each case, the experiments were repeated $500$ times for the directional rank-based tests and the pvMF test at $0.05$ nominal level,  for the various values of $\xi$ (where $\xi=0$ yields ${\mathcal H}_0$). 
Figure~\ref{Fig:RejectRate} shows a plot of the rejection frequencies of these tests against $\xi$, for Cases (1)--(4). Under the null ($\xi=0$), all directional rank-based tests yield rejection rates close to the nominal level   while the pvMF test exhibits a severe over-rejection frequency of $0.122$ under Case~(3), revealing a poor behavior under  multi-modality. In terms of power, for the vMF   location problem (Case (1)), {the vMF-location score, uniform score and vMF-location-concentration score rank-based tests  all have rejection rates very close to the pvMF test which is  the  optimal test in this case. For the vMF-concentration parameter problem (Case~(2)), the pvMF test has no power and is outperformed by all the directional rank-based tests, the power of which increases with~$\xi$.  For the mixture of three vMF distributions (Case (3)), the rejection frequencies of the pvMF test under the alternative are meaningless in view of its size problem;  the uniform, vMF-location, and  vMF-location-concentration score rank-based tests have greater power than the vMF-concentration score one.  When  skewness is present (Case~(4)), all the rank-based tests outperform the classical pvMF test. In general, the superiority of the vMF-location-concentration score and uniform score rank-based tests is clear, considering their consistency, sizeable power against multi-modal alternatives, skewness and vMF distributions (both for   location and concentration treatment effects).  Due to its simplicity, the uniform score is a good choice in terms of computational efficiency.}

\section{Application: MANOVA for Sunspots data}\label{sec:real}\setcounter{equation}{0}

Sunspots are regions on the Sun's photosphere that are darker and cooler than the surrounding areas. They are caused by concentrations of magnetic flux that inhibit convection and are temporary phenomena that experience continuous changes, lasting for hours to days. Their population first increases rapidly and then decreases slowly over a period of approximately 11 years, which is referred to as the {\it solar cycle}. Early in a solar cycle, sunspots appear at higher latitudes and then move towards the equator as the cycle approaches maximum, a phenomenon known as {\it Sp\"{o}rer's law}. Sunspots are widely used to
study and measure solar activity, whose effects may affect earth’s long-term climate \citep{Haigh2007}.

The data we analyze is based on the Debrecen Photoheliographic Data (DPD) sunspot catalogue, which contains locations of sunspots since 1974 and is a continuation of the Greenwich Photoheliographic Results (GPR) catalogue (spanning 1874-1976). The data is available from the R package {\tt rotasym}.  
\begin{figure}[t!]
\centering
\includegraphics[scale = 0.45, trim = 0 0 0 18mm, clip]{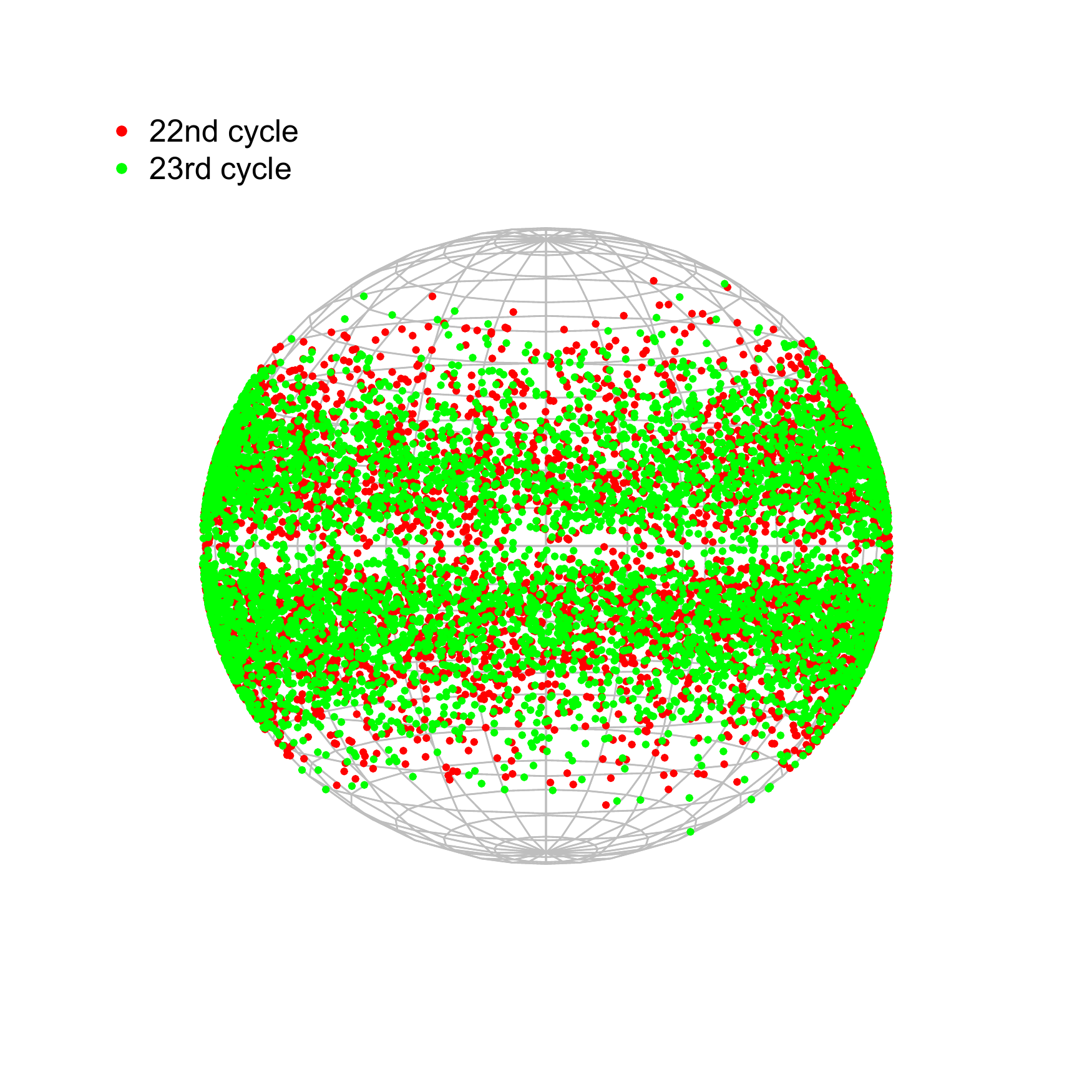}\vspace{-13mm}
\caption{\small Sunspots data: plot of the    sunspots of the 22nd solar cycle ($n_1 = 4551$ red points) and the 23rd ($n_2 = 5373$ green points) solar cycle.}\label{Fig:SolarCycles22and23}
\end{figure}

The sunspots of the 22nd (September 1986 to July 1996; $n_1 = 4551$ points, in red) and~23rd (August 1996 to November 2008; $n_2 = 5373$ points, in green) solar cycles are shown in in Figure~\ref{Fig:SolarCycles22and23}.  Visual inspection   hardly helps decide whether these two samples  are from the same  distribution or not. According to \cite{Garcia-Portugues2020},  various tests   suggest  rotational symmetry   around the north pole   for the 23rd cycle while the same hypothesis is to be rejected  ($p$-values smaller than $0.02$) for the 22nd cycle.

We performed  the pvMF  and the directional rank-based tests based on the uniform, vMF-location and vMF-concentration scores described in Section~\ref{MANOVAsec} for the null hypothesis of equal distributions in the two samples. For the vMF-location and vMF-concentration scores, we factorize $n= 9924$ into $n_S = 121$, $n_R = 82$, and $n_0 = 2$. The $p$-values of these tests  are shown in Table~\ref{Tab:SolarCycles22and23}. The pvMF test, with  $p$-value 0.14, does not reject the null of the equality of distributions, even  at significance level $0.10$. The rank-based test with vMF-concentration scores  rejects the same hypothesis at level~$0.10$, but not at level $0.05$. The rank-based tests with uniform, vMF-location and vMF-location-concentration scores, with  $p$-values at most 0.036, all suggest that the two samples have different distributions. This result is in line with the Monte~Carlo experiments (see, e.g., Case (4)), where the directional rank-based tests significantly outperform the~pvMF test when the underlying distribution are not rotationally symmetric.

\begin{table}[t!]
\centering
\begin{tabular}{lcccccc}
\hline \hline
scores & pvMF & Unif & vMF-location &  vMF-concentration & vMF-location-concentration \\ \hline
$p$-values & 0.140 &  0.036   &  0.005 & 0.052 &  0.005 \\ \hline 
 \end{tabular}
 \vspace{2mm}
 \caption{\small Sunspots data: $p$-values of the pvMF and directional rank-based tests of the equality of distributions of the sunspots of the 22nd  and 23rd solar cycles.}\label{Tab:SolarCycles22and23}
\end{table}

\section{Conclusions} 

In the present paper, we propose various nonparametric tools, based on optimal transport maps from the underlying distribution to the uniform distribution over the hypersphere, for the analysis of directional data. More precisely,  (i)~we propose concepts of directional distribution and quantile functions, distribution-free directional signs and ranks, (ii)~based  on empirical distribution functions, we construct a distribution-free and universally consistent Cram\' er-von Mises-type test of uniformity on $\mathcal{S}^{d-1}$,   and (iii) based on directional ranks and signs, we develop a class of fully distribution-free  MANOVA procedures. Contrary to the pseudo-von Mises methods proposed in the literature---the (asymptotic) validity of which is restricted to rotationally symmetric distributions---our procedures are unrestrictedly finite-sample valid, while  outperforming their  competitors (particularly so in the absence of rotational symmetry.)  \vspace{-2mm}

\section*{Acknowledgements}\vspace{-2mm}
{The authors express their warm thanks to Alberto Gonz\' alez-Sanz for insightful comments and suggestions on the convergence of empirical transport maps. They also  thank the referees and the AE for their their careful reading of the manuscript and  their valuable comments that helped improve the manuscript. 

\section*{Conflict of Interest  Statement}\vspace{-2mm}
 None of the authors are aware of any conflict of interest. \vspace{-2mm} 

\section*{Data availability}\vspace{-2mm}
The Sunspots data is available from the R package {\tt rotasym} \citep{rotasym}. 

\bibliographystyle{wiley-article}
\bibliography{DirectionalStatsMarc}  

\end{document}